\documentclass[10pt]{amsart}

\usepackage{amsmath,amssymb,amsthm}
\usepackage{amsmath,amssymb,amsthm,amscd}
\usepackage[frame,cmtip,arrow,matrix,line,graph,curve]{xy}

\numberwithin{equation}{section}

\DeclareMathOperator{\rank}{rank} 
 \DeclareMathOperator{\Aut}{Aut}

\DeclareMathOperator{\spec}{Spec}

\DeclareMathOperator{\Hom}{Hom}
\DeclareMathOperator{\Ext}{Ext}

\DeclareMathOperator{\codim}{Codim}
\DeclareMathOperator{\image}{Im}

\DeclareMathOperator{\supp}{Supp}
\DeclareMathOperator{\coker}{Coker}
\DeclareMathOperator{\pr}{pr}

\def\hom{\mathop{\cH om\skp}}
\def\endo{\mathop{\cE nd\skp}}
\DeclareMathOperator{\ext}{\cE\!{\it xt}}

\def\cD{\mathcal D}
\def\cE{{\mathcal E}}
\def\cF{{\mathcal F}}
\def\cG{{\mathcal G}}
\def\cH{{\mathcal H}}
\def\cI{{\mathcal I}}

\def\cK{{\mathcal K}}
\def\cL{{\mathcal L}}

\def\cN{{\mathcal N}}
\def\cO{{\mathcal O}}
\def\cP{{\mathcal P}}
\def\cQ{{\mathcal Q}}
\def\cR{{\mathcal R}}

\def\cU{{\mathcal U}}
\def\cV{{\mathcal V}}
\def\cW{{\mathcal W}}
\def\cX{{\mathcal X}}
\def\cY{{\mathcal Y}}
\def\cZ{{\mathcal Z}}

\def\bB{{\mathbf B}}

\def\bD{{\mathbf D}}

\def\bH{{\mathbf H}}
\def\bG{{\mathbf G}}

\def\bE{{\mathbf E}}
\def\bP{{\mathbf P}}

\def\bV{{\mathbf V}}
\def\bY{{\mathbf Y}}

\def\bA{{\mathbf A}}

\def\bM{{\mathbf M}}

\def\bU{{\mathbf U}}
\def\bV{{\mathbf V}}
\def\bp{{\mathbf p}}

\def\bZ{{\mathbf Z}}

\def\bq{{\mathbf q}}
\def\kk{{\rm k}}

\def\ZZ{{\mathbb Z}}
\def\QQ{{\mathbb Q}}
\def\RR{{\mathbb R}}
\def\CC{{\mathbb C}}
\def\II{{\mathbb I}}

\def\FF{{\mathfrak F}}
\def\EE{{\mathfrak E}}
\def\WW{{\mathfrak W}}
\def\LL{{\mathfrak L}}

\def\MM{{\mathfrak M}}

\def\XX{{\mathfrak X}}

\def\AUT{\mathfrak{Aut}}

\newtheorem{prop}{Proposition}[section]
\newtheorem{theo}[prop]{Theorem}
\newtheorem{lemm}[prop]{Lemma}
\newtheorem{coro}[prop]{Corollary}

\newtheorem{defi}[prop]{Definition}

\newtheorem{basi}[prop]{Basic Assumption}
\def\begeq{\begin{equation}}
\def\endeq{\end{equation}}

\def\and{\quad{\rm and}\quad}

\def\Ao{{{\mathbf A}^{\! 1}}}
\def\At{{\mathbf A}^{\! 2}}

\def\An{{{\mathbf A}^{\! n}}}

\def\Anpo{{\mathbf A}^{\!n+1}}

\def\Ampo{{\mathbf A}^{\! m+1}}

\def\alpm{{\alpha^-}}
\def\alpp{{\alpha^+}}
\def\ad{\alpha$-$\bd}

\def\bl{\bigl(}
\def\br{\bigr)}

\def\bd{{\mathbf d}}

\def\cn{{C[n]}}
\def\cl{{C[r]}}

\def\dual{^{\vee}}

\def\defeq{\triangleq}

\def\dag{\dagger}

\def\dpri{^{\prime\prime}}

\def\eps{\epsilon}

\def\gn{G[n]}
\def\gr{G[r]}

\def\Ga{{\mathbf G}^\alpha}

\def\half{\frac{1}{2}}

\def\kk{{\mathbf k}}

\let\lra=\longrightarrow
\def\lsta{_{\ast}}

\def\lalp{_{\alpha}}

\def\lrd{_{r,d}}

\def\lrc{_{r,\chi}}
\def\lrcp{_{r,\chi\pri}}

\def\lnn{_{n,m}}
\def\lnm{\lnn}

\def\lpot{_{p_1+p_2}}
\def\lot{}
\def\lthreeseven{}
\def\lthreefour{}

\def\mapright#1{\,\smash{\mathop{\lra}\limits^{#1}}\,}

\def\mrc{\MM\lrc}
\def\mh{\!:\!}

\def\mm{{\mathfrak m}}

\def\Ma{{\mathbf M}\ualp}

\def\Ma{\bM\ualp}

\def\pri{^{\prime}}
\def\ppri{^{\prime\prime}}

\def\pr{{\rm pr}}
\def\Po{{\mathbf P}^1}

\def\PFT{\PP\cF|_{p_2}\dual}

\def\rk{\mathbf{rk}}

\def\sub{\subset}
\def\sta{^{\ast}}

\def\skp{\hspace{1pt}}
\def\shar{^{\dag}}

\def\st{^{{\rm st}}}

\def\tilcv{{\tilde\cV}}

\def\tga{\tilde{\bG}\ualp}

\def\teo{\tilde E|_{q_0}}
\def\ten{\tilde E|_{q_n}}

\def\upmo{^{-1}}

\def\ualp{^{\alpha}}

\def\upcir{^{\circ}}

\def\upm{^{\pm}}
\def\ump{^\mp}
\def\utf{^{\text{t.f.}}}
\def\ugp{^{\text{G}}}
\def\uz{^{0}}
\def\uo{^{1/2}}
\let\uhalf=\uo

\def\xns{X_n\shar}

\def\wn{{W[n]}}
\def\wr{W[r]}
\def\wm{W[m]}

\def\lab{\label}

\def\solid{\line(1,0){16}}
\def\halfsolid{\line(1,0){8}}
\def\dotline{\qbezier[10](0,0)(8,0)(15,0)}
\def\halfdotline{\qbezier[5](0,0)(4,0)(7,0)}

\def\vsp{\vskip5pt}

\title{Vanishing of the top Chern classes of the moduli of vector bundles}

\author{Young-Hoon Kiem and Jun Li}

\address{Dept of Mathematics, Seoul National University,
Seoul, 151-747, Korea} \email{kiem@math.snu.ac.kr}
\address{Department of Mathematics, Stanford University, Stanford, CA 94305,
USA}\email{jli@math.stanford.edu}
\thanks{Young-Hoon Kiem was partially supported by KOSEF and SNU;
Jun Li was partially supported NSF grants} 

\date{}


\begin{document}
\maketitle

\setcounter{section}{-1}


\section{Introduction}

Let $Y$ be a smooth projective curve of genus $g\ge 2$ and let
$M_{r,d}(Y)$ be the moduli space of stable vector bundles of rank
$r$ and degree $d$ on $Y$.
In case $d$ and $r$ are relarively prime, $M_{r,d}(Y)$ is a smooth
projective variety of dimension $r^2(g-1)+1$. A classical
conjecture of Newstead and Ramanan states that
\begin{equation}
\lab{NRconj} c_i(M_{2,1}(Y))=0\ \ \ \ \text{for} \ i>2(g-1);
\end{equation}
i.e. the top $2g-1$ Chern classes vanish. The purpose of this
paper is to generalize this vanishing result to higher rank cases
by generalizing Gieseker's degeneration method.

In the rank 2 case, there are two proofs of \eqref{NRconj} due to
Gieseker \cite{Gie} and Zagier \cite{Zag}. Zagier's proof is
combinatorial based on the precise knowledge of the cohomology
ring of the moduli space: by the Grothendieck-Riemann-Roch
theorem, Zagier found an expression for the total Chern class
$c(M_{2,1}(Y))$ and then used Thaddeus's formula on intersection
pairing to show the desired vanishing. Because the computation is
extremely complicated even in the rank 2 case, it seems almost
impossible to generalize this approach to higher rank cases.

A more geometric proof of the vanishing \eqref{NRconj} was
provided by Gieseker via  induction on the genus $g$. Let $W\to C$
be a flat family of projective curves over a pointed smooth curve
$0\in C$ such that
\begin{enumerate}
\item $W$ is nonsingular,
\item the fibers $W_s$ over $s\ne 0$ are smooth projective curves of genus
$g$,
\item the central fiber $W_0$ is an irreducible stable curve $X_0$ with one node as its only singular point.
\end{enumerate}
Gieseker constructed a flat family of projective varieties
$\bM_{2,1}(\WW)\to C$ such that \begin{enumerate}
\item the total space $\bM_{2,1}(\WW)$ is nonsingular,
\item the fibers $\bM_{2,1}(\WW_s)$ over $s\ne 0$ are the moduli
spaces $M_{2,1}(W_s)$ of stable bundles over $W_s$,
\item the central fiber $\bM_{2,1}(\WW_0)$ over $0$ has only normal crossing
singularities.
\end{enumerate}
Recently, this construction was generalized to higher rank case by
Nagaraj and Seshadri in \cite{NS} by geometric invariant theory,
and the central fiber $\bM_{r,d}(\WW_0)$ of their construction
parameterizes certain vector bundles on semistable models of
$X_0$. In this paper, we will provide a different construction,
using the technique developed in \cite{Li}.

To prove the vanishing of Chern classes by induction on genus $g$,
Gieseker relates the central fiber $\bM_{2,1}(\WW_0)$ with the
moduli space $M_{2,1}(X)$ where $X$ is the normalization of the
nodal curve $W_0=X_0$. Let $\bM^0$ be the normalization of
$\bM_{2,1}(\WW_0)$, which is a smooth projective variety. Its
general points represent vector bundles on $X_0$ whose pull-back
to $X$ are stable bundles, and hence induces a rational map
\begin{equation}\label{rat}
\bM^0\dashrightarrow M_{2,1}(X).
\end{equation}
Gieseker then proves that the indeterminacy locus of this rational
map is precisely a projective bundle $\mathbb{P}  E^+$ of a vector
bundle $E^+$ over the product $B=Jac_0(X)\times Jac_1(X)$ of
Jacobians, and the normal bundle to $\mathbb{P} E^+$ is the
pull-back of a vector bundle $E^-\to B$ tensored with
$\cO_{\mathbb{P} E^+}(-1)$. This is a typical situation for flips:
we blow up $\bM^0$ along $\mathbb{P} E^+$ and then blow down along
the $\mathbb{P}  E^+$ direction in the exceptional divisor
$\mathbb{P} E^+\times_B\mathbb{P}  E^-$. Let $\bM^1$ be the result
of this flip. Then the rational map becomes a morphism $\bM^1\to
M_{2,1}(X)$, and is a fiber bundle with fiber $\overline{GL(2)}$
--- the wonderful compactification of $GL(2)$ --- that is constructed
as follows. We first compactify $GL(2)$ by embed it in
$$\mathbb{P}  \left(\mathrm{End}(\CC^2)\oplus \CC\right).$$
Its complement consists of two divisors $Z_0'$ (divisor at
infinity) and $Y_1'$ (zero locus of determinant). The wonderful
compactification is by blowing it up along $0\in
\mathrm{End}(\CC^2)$ and along $Y_1'\cap Z_0'$.
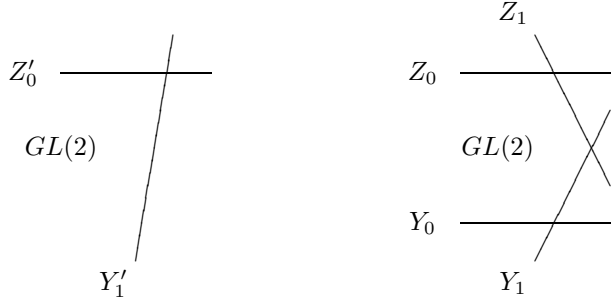
$$\xy
(0,20)*=0{};(20,20)*=0{} **@{-}, (15,25)*=0{};(10,-5)*=0{}
**@{-},(0,10)*=0{GL(2)}, (-5,20)*=0{Z_0'},(7,-8)*=0{Y_1'}
\endxy\qquad \qquad \qquad \qquad
\xy (0,20)*=0{};(20,20)*=0{} **@{-}, (0,0)*=0{};(20,0)*=0{}
**@{-}, (0,20)*=0{};(20,20)*=0{} **@{-}, (10,25)*=0{};(20,5)*=0{}
**@{-}, (20,15)*=0{};(10,-5)*=0{} **@{-}, (5,10)*=0{GL(2)},
(-5,20)*=0{Z_0}, (-5,0)*=0{Y_0}, (7,28)*=0{Z_1}, (7,-8)*=0{Y_1}
\endxy
$$\vspace{.5cm}

The wonderful compactification of $GL(r)$ for $r\ge 2$ is
similarly defined by blowing up $\mathbb{P}
\left(\mathrm{End}(\CC^r)\oplus \CC\right)$  along smooth
subvarieties $2(r-1)$ times and the complement of $GL(r)$ in it
consists of $2r$ smooth normal crossing divisors. This was
carefully studied by Kausz in \cite{Kau}. In summary, Gieseker
obtains the following diagram:
\begin{equation}\label{diagram}
\xymatrix{ &&&\widetilde{\bM}\ar[dl]_{blow-up}\ar[dr]^{blow-up}\\
&&\bM^0\ar[dl]_{normalization} &&\bM^1\ar[dd]^{\overline{GL(2)}}\\
M_{2,1}(Y)\ar@{~>}[r]_{degeneration}&\MM_{2,1}(W)_0&&\\
&&&&M_{2,1}(X)
  }
\end{equation} \normalsize
Afterwards, the proof of the vanishing result \eqref{NRconj} is
reduced to a series of very concrete Chern class computations.

To prove the vanishing of Chern classes for the higher rank case,
we first construct a diagram similar to \eqref{diagram}. In \S
\ref{section1}, we (re-)construct a degeneration $\bM_{r,d}(\WW)$
of Nagaraj and Seshadri, by using the stack of degeneration
defined in \cite{Li}. We take $\bM^0$ as the normalization of the
central fiber of the family $\bM_{r,d}(\WW)$. Next, we define
$\bM^1$ as a fiber bundle over $M_{r,d}(X)$ whose fiber is
$\overline{GL(r)}$ --- the wonderful compactification of $GL(r)$.
Explicitly, working with a universal bundle $\cU\to
M_{r,d}(X)\times X$, $\bM^1$ is the blow-up of
$$\mathbb{P} \left(\mathrm{Hom}(\cU|_{p_1},\cU|_{p_2})\oplus \cO
\right)$$ along suitable smooth subvarieties exactly as in the
construction of $\overline{GL(r)}$. The question is then how to
relate $\bM^0$ with $\bM^1$. Our strategy is to construct a family
of complex manifolds $\bM^{\alpha}$ for $0<\alpha<1$ and study
their variations as $\alpha$ moves from $1$ to $0$. We define a
suitable stability condition for each $\alpha$ (Definition
\ref{adef1.1}) and then show that the set of $\alpha$-stable
vector bundles on semistable models of $X_0$ admits the structure
of a proper separated smooth algebraic space. In particular,
$\bM^\alpha$ is a compact complex manifold.

To prove the vanishing result of Chern classes, we need a very
precise description of the variation of $\bM^\alpha$. We achieved
this for the case of rank 3. Since $M_{3,1}(Y)\cong M_{3,2}(Y)$ by
the morphism $[E]\to[E^*]$ and tensoring a line bundle of degree
1, we only need to consider the case when $r=3$ and $d=1$. By the
stability condition, the moduli spaces $\bM^{\alpha}$ vary only at
$1/3$ and $2/3$. We prove that $\bM^{1/2}$ is obtained from
$\bM^1$ as the consequence of two flips and similarly $\bM^0$ is
the consequence of two flips from $\bM^{1/2}$. The description is
quite explicit and we have the following diagram.
$$\xymatrix{
& &{\bM^0}\ar@{<..>}[r]_{flips}\ar[dl]_{normalization}&{\bM^{1/2}}
\ar@{<..>}[r]_{flips}&{\bM^1}\ar[dd]^{\overline{GL(3)}}\\
M_{3,1}(Y)\ar@{~>}[r]_{{\small degeneration}}&M_{3,1}(X_0)\\
& & & & {M_{3,1}(X)} }
$$
Now it is a matter of explicit but very involved Chern class
computations to verify the vanishing result by induction on genus
$g$. The vanishing result we prove in the end is the following.

\begin{theo}
$c_i(M_{3,1}(Y))=0$ for $i>6g-5$.
\end{theo}

In other words, the top $3g-3$ Chern classes vanish. It seems that
$c_{6g-5}$ should also vanish, but we haven't proved that. Notice
that we also have $c_i(M_{3,2}(Y))=0$ for $i>6g-5$.

The paper is organized as follows. In section one, we will
introduce and construct the moduli space of $\alpha$-stable
bundles on nodal curves; A special case of this construction is
the Gieseker's degeneration for high rank cases. The next section
is devoted to the study of the $\alpha$-stable bundles and the
generalized parabolic bundles on curves. The initial investigation
of the variation of the moduli spaces is carried out in section
three and the detailed study of the flips is achieved in section
four. The last section is about the Chern class calculation.

\section{$\alpha$-stable sheaves and Gieseker's
degeration}\label{section1}

In this section, we will introduce the notion of $\alpha$-stable
sheaves and prove their basic properties. We will then give an
alternative construction of Gieseker's moduli of stable sheaves on
nodal curves in high rank case. In the end, we will show that the
normalization of such moduli spaces can be realized as the moduli
spaces of $\alpha$-stable bundles over marked nodal curves.

\subsection{$\alpha$-stable vector bundles}

Let $g\geq 2$ be an integer and let $X_0$ be a reduced and
irreducible curve of arithmetic genus $g$ with exactly one node,
$q\in X_0$. For $n\geq 0$, we denote by $X_n$ the semistable model
of $X_0$ that contains a chain of $n$-rational curves (i.e.
$\Po$). In this paper, we will fix such an $X_0$ once and for all.
Let $X$ be the normalization of $X_0$, with $p_1$ and $p_2\in X$
the two liftings of the node of $X_0$ under the normalization
morphism. For $X_n$, we denote by $D$ the union of its rational
curves and denote by $D_1,\cdots, D_n$ its $n$ rational
components. We order $D_i$ so that $D_1\cap X=p_1$, $D_i\cap
D_{i+1}\ne\emptyset$ and $D_n\cap X=p_2$. We let
$X^0=X_0-\{q\}=X-\{p_1,p_2\}$, which is an open subset of $X_n$.
We define the based automorphisms of $X_n$ to be
$$\AUT_0(X_n)=\{\sigma\mh X_n \mapright{\cong} X_n\mid \sigma|_{X^0}\equiv \text{id}_{X^0} \}
\cong (\CC^\times)^n.
$$
(Namely, they are automorphisms of $X_n$ whose restrictions to
$X^0\sub X_n$ are the identity maps.)

Later, we need to study pairs $(X_n,q\shar)$, where $q\shar\in
X_n$ are nodes of $X_n$. In this paper, we will call
$(X_n,q\shar)$ based nodal curves, and denote them by $X_n\shar$
with $q\shar\in X_n$ implicitly understood. For $m\geq n$, we say
$\pi\mh X_m\to X_n$ is a contraction if $\pi|_{X^0}$ is the
identity and $\pi|_{D_k}$ is either an embedding or a constant
map. A contraction of $X_m\shar\to X_n\shar$ is a contraction of
the underlying spaces $X_m\to X_n$ that send the based node of
$X_m\shar$ to the based node of $X_n\shar$.

We now fix a pair of positive integers $r\geq 2$ and $\chi$. Let
$X_n\shar=(X_n,q\shar)$ be a based nodal curve and let $E$ be a
rank $r$ locally free sheaf of $\cO_{X_n}$-modules with
$\chi(E)=\chi$. We say $E$ is admissible if the restriction
$E|_{D_i}$ has no negative degree factor\footnote{By
Grothendieck's theorem, every vector bundle on $\Po$ is a direct
sum of line bundles. Each such line bundle is called a factor of
this vector bundle.} for each $i$. (Here and later, for a closed
subscheme $A\sub X_n$ we use $E|_{A}$ to mean
$E\otimes_{\cO_{X_n}}\cO_{A}$.)
Next, we pick an $n$-tuple $\bd=(d_i)\in\ZZ^{+n}$ and let $\eps$
be a sufficiently small positive rational. The pair $(\bd,\eps)$
defines a $\QQ$-polarization $\bd(\eps)$ on $X_n$ whose degree
along the component $\overline{X^0}$ (resp. $D_i$) is
$1-\eps|\bd|$ (resp. $d_i\eps$). (Here $|\bd|=\sum d_i$.)
Now let $F$ be any subsheaf of $E$. We define the rank of $F\sub
E$ at $q\shar\in X_n\shar$ to be
$$ r\shar(F)=\dim \image\{F|_{q\shar}\to E|_{q\shar}\}.
$$
For real $\alpha$ we define the $\ad$-slope (implicitly depending
on the choice of $\eps$) of $F\sub E$ to be
\begin{equation}
\lab{slope} \mu_\bd(F,\alpha)=(\chi(F)-\alpha r\shar(F))/\rk_\bd F
\in\QQ.
\end{equation}
Here the denominator is defined to be
$$\rk_\bd F=
(1-\eps|\bd|)\rank F|_{X^0}+\sum_{i=1}^n \eps d_i \rank F|_{D_i}.
$$
We define the automorphism group $\AUT_0(E)$ to be the group of
pairs $(\sigma,f)$ so that $\sigma\in\AUT_0(X_n)$ and $f$ is an
isomorphism $E\cong\sigma\sta E$.

\begin{defi}\lab{adef1.1}
Let $E$ be a rank $r$ locally free sheaf over $\xns$. Let
$\alpha\in [0,1)$ be any real number. We say $E$ is
$\ad$-semistable (resp. weakly $\ad$-stable) if for any proper
subsheaf $F\sub E$ we have
$$  \mu_\bd(F,\alpha)\leq \mu_\bd(E,\alpha)\quad ({\rm resp}. <\ )
$$
for $\epsilon$ sufficiently small. We say $E$ is $\ad$-stable if
$E$ is weakly $\ad$-stable and $\deg E|_{D_i}>0$ for all $i$.

In case $E$ is a vector bundle on $X_n$ (without the marked node),
we say $E$ is $\bd$-(semi)stable if the same condition hold with
$\alpha=0$.
\end{defi}


We remark that when $\alpha=0$ the $\ad$-stability defined here
coincides with the Simpson stability (compare also the stability
used in \cite{GL}).

We now collect a few facts about $\ad$-stable sheaves on
$X_n\shar$. To avoid complications arising from strictly
semistable sheaves, we will restrict ourselves to the case
$(r,\chi)=1$ and $\alpha\in [0,1)- \Lambda_r$: \begin{equation}
\label{Lambdadef}\begin{array}{ll} \Lambda_r=\{\alpha\in[0,1)\mid
&\alpha=\frac{r_0}{r_0-r\shar}\bl\frac{\chi}{r}-\frac{\chi_0}{r_0}\br,\,
\chi_0, r_0, r\shar\in\ZZ, \\ & 0<r_0<r,\, 0\leq r\shar\leq r,\,
2r_0-r\leq r\shar\leq 2r_0\}\end{array}
\end{equation}
Clearly, $\Lambda_r$ is a discrete subset of $[0,1)$. When
$(\chi,r)=1$, $0\not\in\Lambda_r$.

\begin{lemm}
\lab{1.2} Let $(r,\chi)=1$ and $\chi>r$ be as before. Let
$\alpha\in [0,1)-\Lambda_r$ be any real and let $\bd\in \ZZ^{+n}$
be any weight. Then for any rank $r$ $\ad$-semistable sheaf $E$ on
$X_n\shar$ of Euler characteristic $\chi(E)=\chi$, we have
\newline
(a) The restriction $E|_{D_i}$ has no negative degree factors and
there is no nontrivial section of $E$ which vanishes on $X$;
\newline
(b) For any (partial) contraction $\pi\mh X_m^{\dagger}\to
X_n^{\dagger}$ the pull back $\pi\sta E$ is weakly
$\alpha$-$\bd\pri$-stable for any weight $\bd\pri\in \ZZ^{+m}$;
\newline
(c) Suppose $E$ is $\ad$-stable. Then $n\leq r$ and
$\AUT_0(E)=\CC^\times$.

The similar statement hold for $\bd$-(semi)stable sheaves $E$ on
$X_n$.
\end{lemm}


\begin{proof}
We first prove (a). Suppose $E|_{D_i}$ has a negative degree
factor, say $\cO_{D_i}(-t)$. Let $F$ be the kernel of $E\to
E|_{D_i}\to \mathcal{O}_{D_i}(-t)$. Then $\rk_\bd F<\rk_\bd E$,
$r\shar(F)\leq r$ and
$$\chi(F)=\chi(E)-\chi(\mathcal{O}_{D_i}(-t))=\chi(E)+t-1\ge \chi(E).
$$
Hence,
$$\mu_\bd(F,\alpha)=\frac{\chi(F)-\alpha r\shar(F)}{\rk_\bd F}
>\frac{\chi(E)-\alpha r}{r}=\mu_\bd(E,\alpha).
$$
This is a contradiction. Similarly, suppose there is a section
$s\in H^0(E)$ so that its restriction to $X^0\sub X_n$ is trivial.
Let $L$ be the subsheaf of $E$ generated by this section. Then
since $\rk_\bd L=c\eps$, $c>0$ and $\chi(L)\geq 1$,
$\mu_\bd(L,\alpha)>(1-\alpha)/c\eps>\mu_\bd(E,\alpha)$. This is a
contradiction, which proves (a).

We now prove (b).
Let $F$ be any subsheaf of $E$. Since $E$ is $\ad$-semistable,
\begin{equation}
\lab{1.5} \mu_\bd(F,\alpha)\leq \mu_\bd(E,\alpha).
\end{equation}
Let $r_0(F)=\rank F|_{X^0}$. If $r_0(F)=0$, then
$\mu_\bd(F,\alpha)=c\eps\upmo$, for some $c\in \QQ$. Obviously
$c>0$ is impossible by (\ref{1.5}) and $\chi>r$. In case $c=0$,
then $\chi(F)-\alpha r\shar(F)=0$ and thus $\mu_{\bd'}(F,\alpha)=
0$ for all $\bd\pri$. Since $\chi(E)>r$, the strict inequality in
(\ref{1.5}) holds for all $\bd\pri$. When $c<0$, then
$\mu_{\bd\pri}(F,\alpha)=c\pri\eps\upmo$ for some $c\pri<0$ as
well. Hence the strict inequality in (\ref{1.5}) holds for
$\bd\pri$ too.

We next consider the case $r>\rank E|_{X^0}>0$. Then
$\mu_\bd(F,\alpha)- \mu_\bd(E,\alpha)$ is a continuous function of
$\epsilon$ whose value at $\epsilon=0$ is
$$\frac{\chi(F)-\alpha r\shar(F)}{r_0(F)}-\frac{\chi}{r}+\alpha.
$$
Because $\alpha\in [0,1)-\Lambda_r$, this is never zero. Hence for
sufficiently small $\epsilon$,
$$\mu_\bd(F,\alpha)\leq \mu_\bd(E,\alpha)
\Longleftrightarrow \mu_{\bd\pri}(F,\alpha)<
\mu_{\bd\pri}(E,\alpha).
$$

It remains to consider the case where $r_0(F)=r$. Obviously, we
may consider only proper subsheaves $F$ such that $F|_{X^0}\equiv
E|_{X^0}$ and thus
$E/F$ is a nonzero sheaf of $\cO_D$-modules. Because $E/F$ is a
quotient sheaf of $E|_D$ which is non-negative along each $D_i\sub
D$, $\chi(E/F)\geq r-r\shar(F)$. Since each vector $v\in
E|_{q^\dag}$ which lies in a subspace complementary to $\image
\{F|_{q^\dag}\to E|_{q^\dag}\}$ extends to a section of $E$. Hence
$\chi(E)-\chi(F)>\alpha r-\alpha r^\dag (F)$ and thus
$\mu_{\bd\pri}(F,\alpha)<\mu_{\bd\pri}(E,\alpha)$ for any
$\bd\pri$. This proves that $E$ is weakly $\bd'$-stable for all
$\bd'$.



Now we prove that the pull-back of $E$ to any $X_m$ is weakly
$\alpha$-$\bd\pri$-stable. We first consider the case $m=n+1$. For
simplicity we assume $\pi$ is the contraction of the last rational
component of $X_{n+1}$. We pick $\bd\pri$ so that $d'_i=d_i$ for
$i\le n$ and $d'_{n+1}=1$. Let $F$ be a subsheaf of $\pi^*E$.
Since $F|_{D_{n+1}}$ has no positive degree factor, $\pi\lsta F$
is torsion free and is a subsheaf of $E$. Further, since
$R^1\pi_*(F)$ is a skyscraper sheaf, by Leray spectral sequence we
have
$$\chi(F)=\chi(\pi_*F)-\chi(R^1\pi_*F) \le \chi(\pi_*F).
$$
We now investigate the case $r_0(F)>0$. (The other case can be
treated as in the proof of (b), and will be omitted.) Because $E$
is weakly $\ad$-stable and $\alpha\not\in\Lambda_r$, we must have
$$\frac{\chi(\pi\lsta F)-\alpha r\shar(\pi\lsta F)}{r_0(\pi\lsta F)}< \frac{\chi(E)}{r}-\alpha.
$$
Combined with the above inequality and $r\shar(F)\ge
r\shar(\pi\lsta F)$, we conclude
$$\mu_{\bd\pri}(F,\alpha)<\mu_{\bd\pri}(\pi\sta E,\alpha).
$$
This proves that $\pi\sta E$ is weakly $\alpha$-$\bd\pri$-stable.
The general case $m>n+1$ follows by induction. This completes the
proof of (b).

Now assume $E$ is $\ad$-stable. Consider the vector space
$$V=\{s\in H^0(E|_D)\mid s(p_2)=0\}.
$$
We know that its dimension is no less than $n$ by Riemann-Roch
since $\mathrm{degree}(E|_D)\ge n$. By part (b), the evaluation
map
$$V\lra E|_{p_1}; \quad \text{via}\quad
s\mapsto s(p_1)\in E|_{p_1}
$$
is injective. Because $\dim E|_{p_1}=r$, we have $n\le r$. This
proves (c).

The proof of the second part of (d) is based on the notion of GPB,
and will be proved in section 2 (Corollary \ref{aut}).
\end{proof}

In the light of this lemma, the $\ad$-stability is independent of
the choice of $\bd$ as long as $\alpha\not\in\Lambda_r$. In the
remainder of this paper, we will restrict ourselves to the case
where the following is satisfied.

\begin{basi}
\lab{bas} We assume $(r,\chi)=1$, $\chi>r$ and $\alpha\in
[0,1)-\Lambda_r$
\end{basi}
Henceforth, we can and will call $\ad$-stable simply
$\alpha$-stable with some choice of $\bd$ understood.

\subsection{Gieseker's degeneration of moduli spaces}

In this subsection, we will give an alternative construction of
Gieseker's degeneration of moduli of bundles in high rank case.
The first such construction was obtained by Nagaraj and Seshadri
\cite{NS}.

Let $(r,\chi)=1$ be as before and fix $\bd=(1,\cdots,1)\in \ZZ^n$
for any $n$. Let
$$\cV\lrc(\XX)=\{E\mid E \text{\, is $\bd$-stable on $X_n$ for some
$n$, $\rank E=n$ and $\chi(E)=\chi$}\}/\sim
$$
Here two $E$ and $E\pri$ over $X_n $ and $X_m $ are isomorphic if
there is a based isomorphism $\sigma\mh X_n \to X_m $ so that
$\sigma\sta E\pri\cong E$.

\begin{lemm}
The set $\cV\lrc(\XX )$ is bounded.
\end{lemm}

\begin{proof}
It follows immediately from the bound $n\leq r$ in Lemma
\ref{1.2}.
\end{proof}

Let $0\in C$ be a pointed smooth curve and let $\pi\mh W\to C$ be
a projective family of curves all of whose fibers $W_s$ except
$W_0$ are smooth and the central fiber $W_0$ is the nodal curve
$X_0$ chosen before. Let $C\upcir=C-0$ and $W\upcir=W-W_0$. For
$s\in C\upcir$ we let $\mrc(W_s)$ be the moduli space of rank $r$
and Euler characteristic $\chi$ (namely $\chi(E)=\chi$) semistable
vector bundles on $W_s$. We then let $\mrc(W\upcir/C\upcir)$ be
the associated relative moduli space, namely for $s\in C\upcir$ we
have
$$\mrc(W\upcir/C\upcir)\times_{C\upcir}s=\mrc(W_s).
$$
The goal of this section is to construct the degeneration of the
family $\mrc(W\upcir/C\upcir)$
by filling the central fiber of the family.

Our construction
is aided by the construction of an Artin stack $\WW$
parameterizing all semi-stable models of $W/C$.
Without loss of generality, we can assume $C\sub \Ao$ is a Zariski
open subset with $0\in C$ is the origin of $\Ao$. Let $W[0]=W$ and
let $W[1]$ be a small resolution of
$$W[0]\times_{\Ao}\At,\quad\text{where}\
\At\to\Ao \ \hbox{is via}\ (t_1,t_2)\mapsto t_1t_2.
$$
The small resolution is chosen so that the fiber of $W[1]$ over
$0\in\At$ is $X_1$, and the fibers of $W[1]$ over the first
(second) coordinate line $\Ao\sub\At$ is a smoothing of the first
(resp. second) node of $X_1$. Next $W[2]$ is constructed as a
small resolution of
$$W[1]\times_{\At}\bA^{\!3},\quad\text{where}\
\bA^{\!3}\to\At\ \hbox{is via}\ (t_1,t_2,t_3)\mapsto (t_1,t_2t_3).
$$
The small resolution is chosen so that the fiber of $W[2]$ over
$0\in\bA^{\!3}$ is $X_2$ and the fibers of $W[2]$ over the $i$-th
coordinate line $\Ao\sub\bA^{\!3}$ is a smoothing of the $i$-th
node of $X_2$. The degeneration $W[n]$ is defined inductively. For
the details of this construction please see \cite[\S 1]{Li}.

Let $C[n]$ be $C\times_\Ao \Anpo$. Then $\wn$ is a projective
family of curves over $\cn$ whose fibers are isomorphic to one of
$X_0,\cdots,X_n$. As shown in \cite{Li}, there is a canonical
$\gn\equiv (\CC^\times)^n$ action on $\wn$ defined as follows. Let
$\sigma=(\sigma_1,\cdots,\sigma_n)$ be a general element in $\gn$.
Then the $\gn$ action on $\Anpo$ is
$$(t_1,\cdots, t_{n+1})^\sigma=(\sigma_1 t_1, \sigma_1\upmo\sigma_2 t_2,\cdots,
\sigma_{n-1}\upmo\sigma_n t_n,\sigma_n\upmo t_{n+1}).
$$
The $\gn$ action on $\wn$ is the unique lifting of the trivial
action on $W[0]$ and the above action on $\Anpo$. Consequently, we
can view $\wn/\gn$ as an Artin stack. It is easy to see that
$\wn/\gn$ is an open substack $W[n+1]/G[n+1]$ via the embedding
$\wn/\An\sub W[n+1]/\Anpo$ with $\An\sub\Anpo$ the embedding via
$(t_i)\mapsto(t_i,1)$.

We define the groupoid $\WW$ to be the category of pairs $(W_S/S,
\pi)$, where $S$ are $C$-schemes, $W_S$ are families of projective
curves over $S$ and $\pi\mh W_S\to W$ are $C$-projections, such
that there are open coverings $U\lalp$ of $S$ and $\rho\lalp\mh
U\lalp\to C[n\lalp]$ so that
$$W_S\times_S U\lalp\cong W[n\lalp]\times_{C[n\lalp]}U\lalp,
$$
compatible with the projection $W_S\to W$. Two families $W_S$ and
$W_S\pri$ are isomorphic if there is an $S$-isomorphism $f\mh
W_S\to W_S\pri$ compatible to the tautological projections $W_S\to
W$ and $W_S\pri\to W$. By our construction, $\WW$ is indeed a
stack.

We next define the family of stable sheaves over $\WW$. Let $S/C$
be any scheme over $C$. An $S$-family of locally free sheaves over
$\WW/C$ consists of a member $W_S$ in $\WW$ over $S$ and a flat
family of locally free sheaves $E$ over $W_S/S$. We say $E$ is
admissible (resp. (semi)stable) if for each closed $s\in S$ the
restriction of $E$ to the fiber $W_s=W_S\times_S s$ is admissible
(resp. (semi)stable). Let $W_S$ and $W\pri_S$ be two families in
$\WW(S)$ and let $E$ and $E\pri$ be two families of sheaves over
$W_S$ and $W\pri_S$. We say $E\sim E\pri$ if there is an
isomorphism $f\mh W_S\to W_S\pri$ in $\WW(S)$ and a line bundle
$L$ on $S$ so that $f\sta E\pri\cong E\otimes\pr_S\sta L$.

We now define the groupoid of our moduli problem. For any
$C$-scheme $S$, let $\FF_{r,\chi}(S)$ be the set of equivalence
classes of all pairs $(E, W_S)$, where $W_S$ are members in
$\WW(S)$ and $E$ are flat $S$-families of rank $r$ Euler
characteristic $\chi$ stable vector bundles on $W_S$.


We continue to assume the basic assumption \ref{bas}.

\begin{prop}
\lab{existence} The functor $\FF\lrc(\WW)$ is represented by an
Artin stack $\mrc(\WW)$.
\end{prop}

\begin{proof}
The proof is straightforward and will be omitted.
\end{proof}

Since all stable sheaves have automorphism groups isomorphic to
$\CC^{\times}$, the coarse moduli space $\bM\lrc(\WW)$ of
$\mrc(\WW)$ exists as an algebraic space.

\begin{theo}
The coarse moduli space $\bM\lrc(\WW)$ is separated and proper
over $C$. Further, it is smooth and its central fiber (over $0\in
C$) has normal crossing singularities.
\end{theo}

We divide the proof into several lemmas.

We let $\FF\lrc(\wr)$ be the functor that associates to each
$\cl$-scheme $S/\cl$ the set of equivalence classes\footnote{The
equivalence relation for sheaves over $\wr$ is the usual
equivalence relation. Namely $E\sim E\pri$ if $E\cong
E\pri\otimes\pr_S\sta L$.} of stable sheaves $E$ over
$\wr\times_\cl S$ of rank $r$ and Euler characteristic $\chi$.
Since $(r,\chi)=1$, the stability is independent of the choice of
the polarizations $\bd(\eps)$. By \cite{Sim}, $\FF\lrc(\wr)$ is
represented by an Artin stack $\mrc(\wr)$ and its coarse moduli
space $\bM\lrc(\wr)$ is quasi-projective over $\cl$.

\begin{lemm}
The moduli stack $\mrc(\wr)$ and the moduli space $\bM\lrc(\wr)$
are smooth over $\cl$.
\end{lemm}

\begin{proof}
This is a simple consequence of the deformation theory of sheaves
on $\wr$. Let $[E]\in\mrc(\wr)$ be a closed point, represented by
the sheaf $E$, as a sheaf of $\cO_{\wr}$-modules. Let $Z$ be the
support of $E$, which is a fiber of $\wr$ over $\xi\in\cl$. We
denote by $\iota$ the inclusion $Z\to\wr$. It follows from the
deformation theory of sheaves that the first order deformation of
$E$ is given by the extension group $\Ext^1_{\wr}(E,E)$ which fits
into the exact sequence
$$0\lra H^1(Z,\endo(\iota\sta E))\lra
\Ext^1_{\wr}(E,E)\lra H^0(Z, \endo(\iota\sta E))\otimes
T_\xi\cl\lra 0.
$$
Since $E$ is stable, $H^0(Z, \endo(\iota\sta E))\equiv \CC$.
Further, it is a direct check that the homomorphism of the tangent
spaces
\begin{equation}
\lab{tangent} T_{[E]}\mrc(\wr)\lra T_\xi \cl
\end{equation}
is the next to the last arrow in the above exact sequence.

To show that $\mrc(\wr)$ is smooth, we need to show that there is
no obstruction to deforming $E$. Since the support $Z$ is a closed
fiber of $\wr\to \cl$ and since $E$ is a locally free sheaf of
$\cO_Z$-modules, the obstruction to deforming $E$ lies in
$$H^2(Z,\endo(\iota\sta E))\otimes T_\xi\cl=0.
$$
This shows that $\mrc(\wr)$ is smooth. Since (\ref{tangent}) is
surjective, $\mrc(\wr)$ is smooth over $\cl$. Finally, since the
automorphism group of each $[E]\in\mrc(\wr)$ is isomorphic to
$\CC^\times$, the coarse moduli space $\bM\lrc(\wr)$ is also
smooth over $\cl$. This proves the Lemma.
\end{proof}

\begin{coro}
The coarse moduli space $\bM\lrc(\WW)$ is smooth and is flat over
$C$.
\end{coro}

\begin{proof}
Clearly the $\gr$-action on $\wr$ naturally lifts to an action on
$\mrc(\wr)$ and $\bM\lrc(\wr)$. By Lemma \ref{1.2}, the
stabilizers of the $\gr$-action at all points of $\bM\lrc(\wr)$
are trivial. Hence the quotient $\bM\lrc(\wr)/\gr$ is an algebraic
space. Further, because $\bM\lrc(\wr)$ is smooth,
$\bM\lrc(\wr)/\gr$ is also smooth.

Now let
$$\Phi: \bM\lrc(\wr)/\gr\lra \bM\lrc(\WW)
$$
be the induced morphism. To prove that $\bM\lrc(\WW)$ is smooth it
suffices to show that $\Phi$ is surjective and is \'etale.
Since $\bM\lrc(\wr)$ is the coarse moduli space of the stack
$\mrc(\wr)$, its quotient $\bM\lrc(\wr)/\gr$ is the coarse moduli
space of the stack $\mrc(\wr)/\gr$. Because the nature morphism
$$\mrc(\wr)/\gr\lra \mrc(\WW)
$$
is \'etale, the morphism $\Phi$ is also \'etale. Further, $\Phi$
is surjective because of Lemma \ref{1.2}. This proves that
$\bM\lrc(\WW)$ is smooth. Finally, since $\bM\lrc(\wr)$ is smooth
over $C[r]$ and $C[r]\to C$ is flat, $\bM\lrc(\wr)\to C$ is flat.
Hence $\bM\lrc(\wr)/\gr\to C$ is flat and so does $\bM\lrc(\WW)\to
C$.
\end{proof}

\begin{lemm}
\lab{prop} The algebraic space $\bM\lrc(\WW)$ is separated and
proper over $C$.
\end{lemm}

\begin{proof}
We first check that $\bM\lrc(\WW)$ is proper over $C$, using the
valuation criterion. Let $\xi\in S$ be a closed point in a smooth
curve over $C$. We let $S\upcir=S-\xi$ and let $E\upcir$ be a
family of stable sheaves over $W_{S\upcir}$ for a $C$-morphism
$S\upcir\to \cn$ for some $n$. We need to check that $S\upcir\to
\bM\lrc(\WW)$ extends to $S\to \bM\lrc(\WW)$. Since $\bM\lrc(\WW)$
is flat over $C$, it suffices to check those $S\upcir\to
\bM\lrc(\WW)$ that are flat over $C$. In case $\xi$ lies over
$C\upcir$, then $W_{S\upcir}$ extends to $W_S=W\times_C S$ with
smooth special fiber $W_\xi$. Hence there is an extension of
$E\upcir$ to a family of stable sheaves over $W_S$. We now assume
$\xi$ lies over $0\in C$. Since $S\upcir$ is flat over $C$, $S\to
C$ does not factor through $0\in C$. Without loss of generality we
can assume $S\upcir\to C$ factor through $C\upcir\sub C$. Then
$E\upcir$ is a family of stable sheaves on $W\times_C S\upcir$. We
consider $W_S=W\times_CS$, which possibly has singularity along
the node of $W_\xi$. We let $\tilde W_S$ be the canonical
desingularization of $W_S$. The central fiber $\tilde W_\xi$ (of
$\tilde W_S$ over $\xi$) is isomorphic to $X_m$ for some $m$. We
next fix a polarization $\bd(\eps)$ on $\tilde W_S$ so that its
degree along the irreducible components of $\tilde W_\xi$ are
$1-m\eps, \eps,\cdots,\eps$. Then by \cite{Sim}, there is an
extension of $E\upcir$ to an $\cO_S$-flat sheaf of $\cO_{\tilde
W_S}$-modules $\tilde E$ so that $\tilde E|_{\tilde W_\xi}$ is
semistable with respect to a polarization $\bd(\eps)$. Therefore
$\tilde E$ is locally free since $\tilde W_S$ is smooth, $\tilde
E$ is $\cO_S$-flat and $\tilde E|_{\tilde W_\xi}$ has no torsion
elements supported at points. By Lemma \ref{1.2}, for any rational
curve $D_i\sub \tilde W_\xi$ the restriction $\tilde E|_{D_i}$ is
either trivial ($\cong \cO_{D_i}^{\oplus r}$) or is admissible
along $D_i$. We let $\bar W_S$ be the contraction of $\tilde W_S$
along those $D_i$ so that $\tilde E|_{D_i}$ is trivial. Let
$\pi\mh \tilde W_S\to \bar W_S$ be the contraction morphism and
let $\bar E=\pi\lsta \tilde E$. It is direct to check that $\bar
E$ is locally free and its restriction to $W_\xi$ is admissible
and stable.

It remains to show that the family $\bar W_S$ can be derived from
a $C$-morphism $S\to\cl$. Let $t\in\Gamma(\cO_{\Ao})$ be the
standard coordinate function. Since the exceptional divisor of
$\tilde W_S\to W_S$ has $m$ irreducible components, $\rho\sta t\in
\mm_\xi^{m+1}-\mm_\xi^{m+2}$, where $\rho\mh S\to\Ao$ is the
composite $S\to C\to\Ao$. Hence possibly after an \'etale base
change, we can assume $\rho\sta(t)=\tilde t_0\cdots\tilde t_{m}$
for $\tilde t_i\in\Gamma(\mm_\xi)$. Let $S\to
C[m]=C\times_\Ao\Ampo$ be induced by $S\to C$ and $(\tilde
t_0,\cdots,\tilde t_m)\mh S\to\Ampo$. Then one checks directly
that the fiber product $\wm\times_{C[m]}S$ is isomorphic to
$\tilde W_S$. As to $\bar W_S$, we consider the morphism $S\to
C[m']$ defined as follows. Let
$\tilde\zeta_0,\cdots,\tilde\zeta_m$ be the nodes of $\tilde
W_\xi$ and let $\bar\zeta_0,\cdots,\bar\zeta_{m'}$ be the nodes of
$\bar W_\xi$, ordered according to our convention. Then the
restriction of the contraction $\pi_\xi\mh \tilde W_\xi\to\bar
W_\xi$ induces a surjective map
$$\phi: \{\tilde\zeta_0,\cdots,\tilde\zeta_m\}
\lra\{\bar\zeta_0,\cdots,\bar\zeta_{m'}\}.
$$
We then define $\bar t_k=\prod\{\tilde t_i\mid
\phi(\tilde\zeta_i)=\bar\zeta_k\}$. It is direct to check that
$W[m']\times_{C[m']}S$ is isomorphic to $\bar W_S$. Finally, we
have $m'\leq r$ by Lemma \ref{1.2}. Hence $W[m]\times_{C[m]}S$ can
be realized as a product $W[r]\times_{\cl}S$ for some $S\to\cl$.
This verifies the valuation criterion for properness.

We next show that $\bM\lrc(\WW)$ is separated, using the valuation
criterion. Let $\xi\in S$ be a closed point in a smooth curve over
$C$, let $W_{S\upcir}$ be a family associated to a $C$-morphism
$S\upcir\to \cl$ and let $E\upcir$ be a family of stable sheaves
on $W_{S\upcir}$, as before. To verify the separatedness, we need
to show that there is at most one extension of $E\upcir$ to
families of stable sheaves over $S$. Suppose there are two
extensions $W_S\pri$ and $W_S\dpri$ of $W_{S\upcir}$ and two
extensions $E\pri$ on $W_S\pri$ and $E\dpri$ on $W_S\dpri$ of
$E\upcir$. We need to show that there is a based isomorphism
$\sigma\mh W_S\pri\to W_S\dpri$ and an isomorphism
$E\pri\cong\sigma\sta E\dpri$. Clearly, if $\xi$ lies over a point
in $C\upcir$, this follows from that the moduli of stable sheaves
on smooth curves are separated. We now assume $\xi$ lies over
$0\in C$. Again, because $\bM\lrc(\WW)$ is flat over $C$, we only
need to check those so that $S\upcir$ is flat over $C$. As before,
we let $W_S=W\times_CS$ and let $\tilde W_S$ be the canonical
desingularization of $W_S$. Because both $W_S\pri$ and $W_S\dpri$
are fiber products $W[n\pri]\times_{C[n\pri]}S$ and
$W[n\dpri]\times_{C[n\dpri]}S$ for some $C$-morphisms $S\to
C[n\pri]$ and $S\to C[n\dpri]$, the canonical desingularizations
of both $W_S\pri$ and $W_S\dpri$ are isomorphic to $\tilde W_S$.
Let
$$\pi\pri\mh \tilde W_S\to W_S\pri
\quad{\rm and}\quad \pi\dpri\mh \tilde W_S\to W_S\dpri
$$
be the projections. Then by Lemma \ref{1.2}, both
$\pi^{\prime\ast}E\pri$ and $\pi^{\prime\prime\ast}E\dpri$ are
families of locally free stable vector bundles, extending
$E\upcir$. Hence they are isomorphic. Further, we observe that a
rational curve $D$ is contracted by $\pi\pri$ if and only if
$\pi^{\prime\ast}E\pri|_D$ is trivial. Because
$\pi^{\prime\ast}E\pri \cong\pi^{\prime\prime\ast}E\dpri$, both
$\pi\pri$ and $\pi\dpri$ contract the same set of rational curves
in $\tilde W_S$. Hence $W_S\pri\cong W_S\dpri$ and under this
isomorphism
$$E\pri\equiv \pi\pri\lsta \pi^{\prime\ast}E\pri\cong
\pi\dpri\lsta \pi^{\prime\prime\ast}E\dpri\equiv E\dpri.
$$
This proves the separatedness.
\end{proof}

We now complete the proof of the theorem.

\begin{proof}
It remains to show that the central fiber of $\bM\lrc(\WW)/C$ has
normal crossing singularities. We consider the following Cartesian
square
$$\begin{CD}
\bM\lrc(\wr) @>>> \bM\lrc(\wr)/\gr\\
@VVV @VVV\\
\cl @>>> C=\cl/\gr
\end{CD}
$$
We know that the upper and left arrows are smooth morphisms while
fibers of $\cl\to C$ are smooth except the central fiber, which
has normal crossing singularities. Thus the central fiber of the
right arrow has normal crossing singularities. Finally, since
$\bM\lrc(\wr)/\gr\to \bM\lrc(\WW)$ is \'etale and the arrow
commutes with the two projections to $C$, the fibers of
$\bM\lrc(\WW)$ over $C$ are smooth except the central fiber which
has normal crossing singularities.
\end{proof}

We conclude this subsection by remarking that the set of closed
points of $\bM\lrc(\WW)$ over $0\in C$ is exactly the set
$\cV\lrc(\XX)$ defined in the beginning of this subsection. This
set is exactly the set used by Gieseker and Nagaraj - Seshadri in
their construction degeneration of moduli spaces \cite{Gie, NS}.

\subsection{Normalization of the central fiber}

We close this section by constructing the normalization of the
central fiber
$$\bM\lrc(\WW_0)=\bM\lrc(\WW)\times_C 0.
$$
Let $[E]\in\mrc(\WW_0)$ be any closed point associated to a stable
vector bundle $E$ on $X_n$. We consider the coarse moduli space
$\bM\lrc(\wn)$ of stable sheaves on $\wn/\cn$. As we argued
before, the induced morphism $\bM\lrc(\wn)/\gn\to \bM\lrc(\WW)$ is
\'etale and its image contains $[E]$. For the same reason,
$$\bM\lrc(\wn)\times_{\cn}(\cn\times_C0)/\gn\lra \bM\lrc(\WW_0)
$$
is \'etale. Now let $\bH^k\sub\Anpo$ be the coordinate hyperplane
transversal to the $k$-th coordinate axis (we agree that the
coordinate axes of $\Anpo$ are indexed from $0$ to $n$) and let
$\bH^k\to \cn$ be the induced immersion. Then
$$\coprod_{k=0}^n \bH^k\lra \cn\times_C 0
$$
is the normalization morphism and the induced morphism
$$\coprod_{k=0}^n(\bM\lrc(\wn)\times_{\cn}\bH^k)/\gn\lra
\bM\lrc(\wn)\times_{\cn}(\cn\times_C0)
$$
is the normalization morphism. It follows that the closed points
of the normalization of $\bM\lrc(\WW_0)$ consists of the
equivalence classes of triples $(E,X_n,q\shar)$ for some $n$ where
$E$ are stable vector bundles over $X_n$ and $q\shar$ are nodes of
$X_n$.

The normalization of $\bM\lrc(\WW_0)$ is indeed itself a moduli
space. Let $\XX\shar$ be the Artin stack (groupoid) of pointed
semistable models of $X_0$. Namely, $\XX\shar(S)$ consists of
pairs $(W_S, q\shar)$ where $W_S$ are members in $\WW(S)$ with $S$
understood as $C$-schemes via the trivial morphisms $S\to 0$ with
$0\in C$, and $q\shar$ are sections of nodes of the fibers of
$W_S/S$. An isomorphism between two families $(W_S,q\shar)$ and
$(W_S\pri,q^{\prime\dagger})$ is an isomorphism $\sigma\mh W_S\to
W_S\pri$ in $\WW(S)$ that preserves the sections $q\shar$ and
$q^{\prime\dagger}$.  We denote the pair $(W_S,q\shar)$ by
$W_S\shar$. We define the moduli groupoid $\FF\ualp\lrc(\XX\shar)$
to be the category of all pairs $(E, W_S\shar)$, where
$W_S\in\XX\shar(S)$ and $E$ is a family of $\alpha$-stable sheaves
on $W_S\shar$ whose members (namely the restriction to closed
fibers of $W_S/S$) lie in the set  $\cV\ualp\lrc(\XX\shar)$ of
isomorphism classes of $\alpha$-stable rank $r$ locally free
sheaves $E$ on $X_n^{\dag}$ for some $n$ with $\chi(E)=\chi$. Two
families $(E,W_S\shar)$ and $(E\pri, W_S^{\prime\dagger})$ are
equivalent if there is an isomorphism $\sigma\mh W_S\shar\to
W_S^{\pri\dagger}$ in $\WW(S)$
 and a line bundle $L$ on $S$ so that $f\sta
E\pri\cong E\otimes\pr_S\sta L$. Following the proof of Theorem
\ref{existence}, one easily shows that the functor
$\FF\lrc\ualp(\XX\shar)$ is represented by a smooth Artin stack.
Further, because of Corollary \ref{aut} its coarse moduli space
$\bM\lrc\ualp(\XX\shar)$ is a smooth algebraic space. Again,
similar to the proof of Lemma \ref{prop}, the coarse moduli space
is proper and separated. Finally, the previous argument shows that
$\bM\lrc^0(\XX\shar)$ is a normalization of $\bM\lrc(\WW_0)$. We
summarize this in the following proposition.

\begin{prop}
The groupoid $\FF\lrc\ualp(\XX\shar)$ forms a smooth Artin stack
$\mrc\ualp(\XX\shar)$ and its coarse moduli space
$\bM\lrc\ualp(\XX\shar)$ is a proper smooth and separated
algebraic space. Further, the canonical morphism
$\bM\lrc^0(\XX\shar)\to \bM\lrc(\WW_0)$ induced by forgetting the
marked section of nodes is the normalization morphism.
\end{prop}

In case $r=2$, the moduli space $\bM_{2,\chi}(\XX\shar)$ is
exactly the normalization constructed in \cite{Gie}.


\def\a{\alpha}
\def\PP{{\mathbb P}}
\def\gpb{\mathcal{GPB}}

\section{$\alpha$-stable bundles and generalized parabolic bundles}

In this section, we will first relate $\alpha$-stable sheaves on
$X\shar$ to Generalized-Parabolic-Bundle (in short GPB) on
$X^+=(X,p_1+p_2)$. We will then show that the moduli of
$\alpha$-stable sheaves on $X\shar$ is a blow-up of the moduli of
$\alpha$-stable GPB on $X^+$. This will be used to study how
$\bM^0\lrc(\XX\shar)$ is related to $\bM^{1^-}\lrc(\XX\shar)$ in
the next section. We will continue to assume $(r,\chi)=1$,
$\chi>r$ and $\alpha\not\in \Lambda_r$ throughout this section.

\subsection{GPB}

We begin with the notion of GPB. Let $X^+$ be the pair
$(X,p_1+p_2)$. A rank $r$ GPB on $X^+$ is a pair $V^G=(V,V^0)$ of
a rank $r$ vector bundle $V$ on $X$ and an $r$-dimensional
subspace $V^0\sub V|\lpot$. In this paper, we will use the
convention that for any sheaf $F$ and closed $p\in X$ we denote by
$F|_p$ the vector space $F\otimes\kk(p)$ and denote by
$V|_{p_1+p_2}$ the vector space $V|_{p_1}\oplus V|_{p_2}$. For any
subsheaf $F\sub V$ we denote by $F|_{p_1+p_2}\cap V^0$ the
subspace $\image\{F|\lpot\to V|\lpot\}\cap V^0$ and define
$r^+(F)=\dim F|_{p_1+p_2}\cap V^0$.

We begin with the investigation of locally free sheaves on a chain
of rational curves. Let $R$ be a chain of $n$ $\Po$'s coupled with
two end points $q_0$ and $q_n$ in the first and the last
components of $R$. We order the $n$-rational curves into
$R_1,\cdots,R_n$ so that $R_i\cap R_{i+1}=\{q_i\}$ for $1\le i\le
n-1$, $q_0\ne q_1\in R_1$ and $q_n\ne q_{n-1}\in R_n$. In the
following we will call such $R$ with $q_0$ and $q_n$ understood an
end-pointed chain of rational curves. Now let $F$ be any
admissible locally free sheaf on $R$. Inductively, we define
vector spaces $W_i\sub F|_{q_i}$ inductively by $W_0=\{0\}$ and
$$W_i=\{s(q_i)\mid s\in H^0(R_i,F), s(q_{i-1})\in W_{i-1} \}
\sub F|_{q_i}.
$$
We call $T_\to=W_n$ the \emph{transfer} of $0\in F|_{q_0}$ along
$R$. Note that
\begin{equation}
\lab{tran} T_\to=\{s(q_n)\mid s(q_0)=0\ \text{and}\  s\in
H^0(R,F)\}.
\end{equation}
If we reverse the order of $R$ by putting $\tilde{R}_i=R_{n-i+1}$,
we call the resulting transfer $T_\leftarrow\sub F|_{q_0}$ the
\emph{reverse transfer} of $0\in F|_{q_n}$. Notice that we have a
well-defined homomorphism $F_{q_n}\to F|_{q_0}/T_\leftarrow$ by
assigning to each element $c\in F|_{q_n}$ the class of
$[s(q_0)]\in F|_{q_0}/T_\leftarrow$ for some $s\in H^0(R,F)$ such
that $s(q_n)=c$. The kernel of this homomorphism is precisely the
transfer $T_\to$. Hence we have a canonical isomorphism
\begin{equation}
\lab{1.12} \xi: F|_{q_n}/T_\to\mapright{\cong}
F|_{q_0}/T_\leftarrow.
\end{equation}
There is another way to see this isomorphism. Let
$H^0(R,F\dual)\otimes\cO_R\to F\dual$ and
$$\varphi: F\lra H^0(R,F\dual)\dual\otimes \cO_R
$$
be the canonical homomorphism. Then
$\ker(\varphi(q_0))=T_\leftarrow$, $\ker(\varphi(q_n))=T_\to$ and
the isomorphism (\ref{1.12}) is induced by
$$F|_{q_0}\lra H^0(R,F\dual)\dual\longleftarrow F|_{q_n}.
$$

\begin{defi}
\lab{def} We say a locally free sheaf $F$ on $R$ is regular if
there are integers $a_i$ so that to each $i\in [1,n]$,
$$F|_{R_i}=\cO_{R_i}^{r-a_i}\oplus \cO_{R_i}(1)^{a_i}
\quad\text{and}\quad \dim W_i=\dim W_{i-1}+a_i.
$$
\end{defi}

Note that $F$ is regular if and only if each restriction
$F|_{R_i}$ has only degree $0$ and $1$ factors and
\begin{equation}\label{reg2}
\{s\in H^0(R,F)\mid s(q_0)=s(q_n)=0\}=0.
\end{equation}
Also, if $F$ is regular then $\deg F=\dim T_\to=\dim
T_\leftarrow$. We now prove a lemma concerning regular bundles on
rational chains that will be useful later.

\begin{lemm}\lab{1.15}
Let the notation be as before and let $\sigma_0\mh F|_{q_0}\to
F|_{q_0}/T_\leftarrow$ and $\sigma_n\mh F|_{q_n}\to
F|_{q_n}/T_\to$ be the projections. Suppose $S_0\sub F|_{q_0}$ and
$S_n\sub F|_{q_n}$ are two subspaces so that
$\xi(\sigma_n(S_n))\sub \sigma_0(S_0)$. Then there is a subsheaf
$\cF\sub F$ so that $\image\{\cF|_{q_0}\to F|_{q_0}\}=S_0$,
$\image\{\cF|_{q_n}\to F|_{q_n}\}=S_n$ and
$$\chi(\cF)\geq \dim S_n+\dim S_0-\dim \sigma_0(S_0).
$$
\end{lemm}

\begin{proof}
The proof is straightforward, which is based on the following easy
observation. Suppose $A\sub F|_{q_0}$ and $B\sub F|_{q_n}$ are one
dimensional subspaces so that $\sigma_0(A)=\xi(\sigma_n(B))\ne
\{0\}$. Then there is a unique subsheaf $\cL\cong\cO_R\sub F$ so
that $\image\{\cL|_{q_0}\to F|_{q_0}\}=A$ and
$\image\{\cL|_{q_n}\to F|_{q_n}\}=B$.

We now construct the subsheaf $\cF$. First, we write
$S_n=A_1\oplus A_2$ so that $A_2=\ker\sigma_n\cap S_n$. Then by
assumption, there is a subspace $A_1\pri\sub S_0$ so that
$A\pri_1$ is isomorphic to $A_1$ under
$\sigma_n\upmo\circ\xi\upmo\circ\sigma_0$. We let
$A_3=\ker\sigma_0\cap S_0$ and $A_4\sub S_0$ be the compliment of
$A\pri_1\oplus A_3$. Then by the observation just stated, there is
a subsheaf $\cF_1\cong\cO_R^{\oplus a_1} \sub F$, where $a_1=\dim
A_1$, so that $\cF_1|_{q_0}\sub F|_{q_0}$ is $A_1$ and
$\cF_1|_{q_n}\sub F|_{q_n}$ is $A\pri_1$. We then pick $\cF_4\cong
\cO_R^{\oplus a_4}\sub F$, where $a_4=\dim A_4$, so that
$\cF_4|_{q_0}\sub F|_{q_0}$ is $A_4$. We then define $\cF_2\sub F$
to be the subsheaf spanned by $\{s\in H^0(F(-q_n))\mid s(q_0)\in
A_2\}$. Similarly, we defined $\cF_3\sub F$ to be the subsheaf
spanned by sections in $H^0(F(-q_0))$ whose values at $q_n$ are in
$B_2$. Clearly, $\chi(\cF_i)=a_i$, where $a_i=\dim A_i$. Finally,
let $\cF\sub F$ be the image sheaf of
$$\varphi:
\cF_1\oplus\cF_2\oplus\cF_3\oplus\cF_4(-q_n)\to F.
$$
By our construction, we certainly have $\image\{\cF|_{q_i}\to
F|_{q_i}\}=S_i$ for $i=0$ and $n$. As $\chi(\cF)$, it is equal to
$\sum\chi(\cF_i)-a_4-\chi(\ker \varphi)$. Because the restriction
of $\varphi$ to $q_0$ and $q_n\in R$ are injective, the structures
of $\cF_i$ guarantee that $H^0(R,\ker\varphi)=0$. Hence
$\chi(\ker\varphi)\leq 0$. This proves the inequality of
$\chi(\cF)$.
\end{proof}

Now we pick a pair of non-negative integers $(n,m)$, and form the
two end-pointed rational chains $R$ and $R\pri$ of length $n$ and
$m$, respectively. Let $q_0,q_n\in R$ and $q_0\pri,q_m\pri\in
R\pri$ be the respective end points with their nodes $q_i$ and
$q_i\pri$, respectively. We then form the $2$-pointed curve
$X_{n,m}$ by gluing $p_1$ and $p_2$ in $X$ with the $q_0\in R$ and
the $q_0\pri\in R\pri$, respectively, with $q_n$, $q_m\pri$ its
two marked points. Namely, $X_{n,m}$ has two tails of rational
curves, the left tail $R\ni q_n$ and the right tail $R\pri\ni
q_m\pri$.
If we identify $q_n$ and $q_m\pri$, we obtain a pointed curve
$X_{n+m}\shar$.

Now let $E$ be a rank $r$ vector bundle over $X_{n+m}\shar$ so
that its restriction to the chain of rational curves $D\sub
X_{n+m}$ is regular. Let
$$\rho: X\lnn\lra X_{n+m}\shar
\quad\text{and}\quad \pi:X\lnn\lra X
$$
be the tautological projections. First, $\tilde E\defeq\rho\sta E$
is a locally free sheaf on $X\lnn$ and $E$ can be reconstructed
from $\rho\sta E$ and the isomorphism $\phi: \tilde E|_{q_n}\equiv
E|_{q\shar}\equiv \tilde E|_{q_m\pri}$ via
$$0\lra E\lra \rho\lsta\tilde E\lra (\tilde E|_{q_n}\oplus \tilde E|_{q_m\pri})/\Gamma_{\phi}\lra 0
$$
where $\Gamma_\phi\sub \tilde E|_{q_n}\oplus \tilde E|_{q_m\pri}$
is the graph of $\phi$. Here we view the last non-zero term in the
sequence as a $\kk(q\shar)$ vector space, which is naturally a
sheaf of $\cO_{X_{n+m}}$-modules. In this way the vector bundle
$E$ on $X_{n+m}\shar$ is equivalent to the GPB $(\tilde
E,\Gamma_\phi)$ over $X\lnn$.

We now show how the GPB $(\tilde E,\Gamma_\phi)$ naturally
associates to a GPB over $X^+$. First let $V=(\pi\lsta \tilde
E\dual)\dual$. Since the restriction of $\tilde E$ to $D$ is
regular, $\rho\lsta\tilde E\dual$ and hence $V$ are locally free
sheaves on $X$. Clearly, we have $\chi(V)=\chi(\tilde
E)=\chi(E)+r$. Next, by its construction we have canonical
$\pi\sta(\pi\lsta\tilde E\dual)\to\tilde E\dual$ and its dual
\begin{equation}
\lab{a} \tilde E\lra \pi\sta V\equiv \pi\sta(\pi\lsta\tilde
E\dual)\dual.
\end{equation}
Restricting to $q_n$ and $q_m\pri$, we obtain
$$h_1: \tilde E|_{q_n}\lra \pi\sta V|_{q_n}\equiv V|_{p_1}
\quad\text{and}\quad h_2: \tilde E|_{q_m\pri}\lra \pi\sta
V|_{q_m\pri}\equiv V|_{p_2}.
$$
We then define $V^0\sub V|_{p_1}\oplus V|_{p_2}$ to be the image
of $\Gamma_\phi\sub \tilde E|_{q_n}\oplus\tilde E|_{q_m\pri}$
under the homomorphism $h_1\oplus h_2$. We claim that $\dim
V^0=r$. Suppose $\dim V^0<r$.  Then there is $v\in \Gamma_\phi$ so
that $h_1(v)=h_2(v)=0$. Since $v\in \Gamma_\phi$, its image in
$\tilde E|_{q_n}$ and in $\tilde E|_{q\pri_m}$ are identical,
using $\tilde E|_{q_n}\equiv\tilde E|_{q\pri_m}$. By our previous
discussion, $h_1(v)=0$ implies that $v$ lies in the kernel $\tilde
E|_{q_n}\to \tilde E|_{q_0}/T_\leftarrow$, which implies that
there is a section $s\in H^0(\tilde E|_R(-q_0))$ so that
$s(q_n)=v$. Similarly, $h_2(v)=0$ implies that there is a section
$s\pri\in H^0(\tilde E|_{R\pri}(-q_0\pri))$ so that
$s\pri(q_m\pri)=v$. The pair $(s,s\pri)$ then glue together to
form a section of $E$ that vanishes along $X^0$. Since $E|_{D}$ is
regular, such section must be trivial, and hence $v=0$. This
proves that $\dim V^0=r$ and hence the pair $V^G=(V,V^0)$ is a GPB
on $X$.

\subsection{$\alpha$-stable bundles and $\alpha$-stable GPB}

In this subsection we will show that the correspondence
constructed in the previous subsection relates $\alpha$-stable
bundles on $X\shar_{n+m}$ to $\alpha$-stable GPBs on $X^+$.

\begin{defi}
Let $\a\in[0,1)$. A GPB $V^G=(V,V^0)$ is $\alpha$-stable if for
any proper subbundle $F\sub V$ we have
$\mu^G(F,\alpha)<\mu^G(V,\alpha)$, where
$$\mu^G(F,\alpha)=(\chi(F)+(1-\alpha)r^+(F))/r(F).\footnote{The
$\alpha$-stability is the $(1-\alpha)$-stability of GBP introduced
in \cite{Bho}.}
$$
\end{defi}

Let $\cG\ualp\lrcp(X^+)$ be the set of all isomorphism classes of
$\alpha$-stable rank $r$ GPBs on $X^+$ of Euler characteristics
$\chi\pri=\chi+r$. In this subsection, we will show that the
previous correspondence defines a map
$$\cV\ualp\lrc(\XX\shar)\lra \cG\ualp\lrcp(X^+).
$$
This map will be useful in studying the moduli of $\alpha$-stable
sheaves on $\XX\shar$.

We now prove the following equivalence of the stable GPB and
stable vector bundles on $X\shar_{n+m}$.

\begin{prop}
\lab{2.16} Let $E$ be a rank $r$ vector bundle of Eular
characteristic $\chi$ on $X\shar_{n+m}$. Suppose $E|_{D}$ is
regular.\footnote{By Lemma \ref{1.2}, $E|_D$ is regular whenever
$E$ is $\alpha$-stable.} Then $E$ is $\a$-stable if and only if
its associated GPB $V^G=(V,V^0)$ is $\a$-stable.
\end{prop}

\begin{proof}
We first prove that $E$ being $\a$-stable implies that $V^G$ is
$\alpha$-stable. Let $U\sub V$ be any proper subbundle. We need to
show that $\mu^G(U,\alpha)<\mu^G(V,\alpha)$. We will prove this
inequality by constructing a subsheaf $F\sub E$ so that
$r_0(F)=r(U)$, $\chi(F)\geq\chi(U)-2r(U)+r^+(U)$ and
$r\shar(F)=r^+(U)$. Then we have the inequality because
$$\mu^G(U,\alpha)=\frac{\chi(U)+(1-\alpha)r^+(U)}{r(U)}
\leq\frac{(\chi(F)+2r(U)-r^+(U))+(1-\alpha)r^+(U)}{r(U)}
$$
$$
=\frac{\chi(F)-\alpha r\shar(F)}{r_0(F)}+2<\frac{\chi(E)-\alpha r
}{r}+2 =\frac{\chi(V)+(1-\alpha)r}{r}=\mu^G(V,\alpha).
$$

We now construct such $F$. First let $\tilde U=\tilde E|_X\cap U$,
which makes sense since by construction $V$ is just the result of
a elementary (Hecke) modification of $\tilde{E}|_X$ at $p_1$ and
$p_2$ so that
 $\tilde E|_X\sub V$ is a subsheaf. Let $A_1=\image\{\tilde U|_{p_1}\to \tilde E|_{q_0}\}$ and
$A_2=\image\{\tilde U|_{p_2}\to \tilde E|_{q_0\pri}\}$. We then
let $\tilde B\sub \Gamma_\phi$ be the preimage of $U^0=V^0\cap
(U|_{p_1}\oplus U|_{p_2})$ under the isomorphism
$V^0\cong\Gamma_\phi$, and let $B_1\sub \tilde E|_{q_n}$ and
$B_2\sub \tilde E|_{q_m\pri}$ be the image of $\tilde B$ under the
obvious projections. Let $\sigma_0\mh \teo\to \teo/T_\leftarrow$
and $\sigma_n\mh \ten\to\ten/T_\to$ be the projections and $\xi\mh
\ten/T_\to\cong\teo/T_\leftarrow$ be the tautological
isomorphisms, constructed in (\ref{1.12}). We claim that
$\xi(\sigma_n(B_1))\sub \sigma_0(A_1)$. Indeed, let $(u_1,u_2)\in
\tilde B$ be any element with $(v_1,v_2)\in U^0$ its preimage. By
definition, $\xi(\sigma_n(v_1))=\sigma_0(u_1)$. Hence
$\xi(\sigma_n(B_1))\sub \sigma_0(A_1)$. Similarly, if we let
$\sigma_0\pri$, $\sigma_m\pri$ and $\xi\pri$ be similar
homomorphisms associated to $q_0\pri$, $q_m\pri$ and $p_2$, we
have $\xi\pri(\sigma_m\pri(B_2))\sub\sigma_0\pri(A_2)$.


We now apply Lemma \ref{1.15} to conclude that there is a subsheaf
$F_1\sub \tilde E|_R$ so that
$$\image\{F_1|_{q_0}\to \tilde E|_{q_0}\}=A_1,\ \image\{F_1|_{q_n}\to \tilde E|_{q_n}\}
=B_1,\ \chi(F_1)\geq \dim B_1+e_1,
$$
where $e_1=\dim \ker \{\tilde U|_{p_1}\to U|_{p_1}\}$. Similarly,
we have a subsheaf $F_2\sub \tilde E|_{R\pri}$ having
$$\image\{F_2|_{q_0\pri}\to \tilde E|_{q_0\pri}\}=A_2,\ \image\{F_2|_{q_m\pri}\to \tilde E|_{q_m\pri}\}
=B_2, \ \chi(F_2)\geq \dim B_2+e_2,
$$
where $e_2=\dim \ker \{\tilde U|_{p_2}\to U|_{p_2}\}$. Note that
$B_1=B_2$ under the identification $\tilde E|_{q_n}\equiv \tilde
E|_{q'_m}$. Therefore, the subsheaves $\tilde U\sub \tilde
E|_{X}$, $F_1\sub \tilde E|_{R}$ and $F_2\sub \tilde E|_{R\pri}$
glue together to form a subsheaf $F\sub E$ such that
$r_0(F)=r(U)$, $r\shar(F)=\dim B_1=\dim B_2=\dim U^0=r^+(U)$ and
$$\chi(F)=\chi(\tilde U)+\chi(F_1)+\chi(F_2)-2r(U)-\dim C
\geq\chi(U)-2r(U)+\dim U^0.
$$
Here we used the fact that $\chi(\tilde U)=\chi(U)-e_1-e_2$. This
$F$ is the desired subsheaf. This prove the first part of the
lemma.

We postpone the other part of the proof until we give a more
precise description of the destabilizing subsheaf of $E$.
\end{proof}

\subsection{$\alpha$-stable GPB and $\alpha$-stable bundles}

We will complete the other half of Proposition \ref{2.16} in this
subsection.

We begin with a characterization of the destabilizing subsheaf of
$E$ on $X\shar_{n+m}$. Let $E$ be a vector bundle on
$X\shar_{n+m}$ as in Proposition \ref{2.16}. Let $F\sub E$ be an
$\alpha$-destabilizing subsheaf. Namely,
$\mu(F,\alpha)>\mu(E,\alpha)$ and is the largest possible among
all subsheaves of $E$. Since $E|_D$ is regular, as mentioned
before $r_0(F)=\rank F|_{X^0}>0$. In the following, we say that a
sheaf $F\pri$ with $F\sub F\pri\sub E$ is a small extension of
$F\sub E$ if $F\pri/F$ is a sheaf of $\cO_D$-modules. We say
$F|_D=F\otimes_{\cO_{X_n}}\cO_D$ is non-negative if the torsion
free part of the restriction of $F$ to each rational curves
$D_i\sub D$ has no-negative factors.

\begin{lemm}
\lab{3.16} Let $E$ be as in Proposition \ref{2.16} and $F\sub E$
be an $\alpha$-destabilizing subsheaf. Then $r_0(F)>0$, $F|_D$ is
non-negative and there are no non-negative small extension $F\sub
F\pri\sub E$ of $F\sub E$.
\end{lemm}

\begin{proof}
First, $r_0(F)=0$ is impossible because of our assumption
$\chi>r$. Suppose there is an irreducible component $D_i\sub D$ so
that $F|_{D_i}$ is not non-negative. Namely, there is a $t>0$ so
that $\cO_{D_i}(-t)$ is a quotient sheaf of $F|_{D_i}$. We let
$F\pri$ be the kernel of the composite $F\to F|_{D_i}\to
\cO_{D_i}(-t)$. Clearly, $\chi(F\pri)\geq \chi(F)$ and
$r\shar(F\pri)\leq r\shar(F)$ while $\rk_\bd F\pri<\rk_\bd F$.
Hence $\mu_\bd(F\pri,\alpha)>\mu_\bd(F,\alpha)$, violating that
$F$ is an $\alpha$-destabilizing subsheaf of $E$.

Now suppose $F\sub F\pri\sub E$ is a small extension of $F$ so
that $F\pri$ is also non-negative on $D$. We claim that
$\mu_\bd(F\pri,\alpha)>\mu_\bd(F,\alpha)$. We look at the quotient
sheaf $F\pri/F$. Since $F\pri|_D$ is non-negative, $F\pri/F$ is
also non-negative. Further, since $F\pri/F$ is a sheaf of
$\cO_D$-modules, it is easy to see that
$$\chi(F\pri/F)\geq r\shar(F\pri)-r\shar(F).
$$
Hence because $\alpha<1$, the slope $\mu_\bd(F\pri,\alpha)$ is
$$\frac{(\chi(F)-\alpha r\shar(F))+(\chi(F\pri/F)-\alpha(r\shar(F\pri)-r\shar(F)))}
{\rk_\bd F+O(\eps)}
>\frac{\chi(F)-\alpha r\shar(F)}{\rk_\bd F}=\mu_\bd(F,\alpha).$$
Here we used the fact that $r_0(F)>0$ and that $\eps$ is
sufficiently small. This violates the assumption that $F$ is an
$\a$-destabilizing subsheaf of $E$.
\end{proof}


We now give a more precise description of the destabilizing
subsheaf $F\sub E$. We begin with more notation. Let
$D_1,\cdots,D_{n+m}$ be the ordered rational curves of $D\sub
X\shar_{n+m}$ with nodes $q_0,\cdots,q_{n+m}$ and the marked node
$q\shar=q_n$. We let $D_{[i,j]}=\cup_{k=i+1}^j D_k$ be the
subchain of $D$. There are several (possible) subsheaves of $E|_D$
that are important to our later study. The first is the subsheaf
$\cO_{[i,j]}\sub E|_D$, which as a sheaf is isomorphic to
$\cO_{D_{[i,j]}}$ and such that the induced homomorphisms
$\sigma(q_i)\mh \cO_{[i,j]}\otimes\kk(q_i)\to E|_{q_i}$ and
$\sigma(q_j)\mh \cO_{[i,j]}\otimes\kk(q_j)\to E|_{q_j}$ are
non-zero. If we impose the condition $\sigma(q_i)\ne 0$ instead of
$\sigma(q_i)=0$, we denote the resulting subsheaf by
$\cO_{(i,j]}$. The sheaves $\cO_{[i,j)}$ is defined similarly in
the obvious way. The other subsheaf is $\cO_{[0,n+m]}^{[k]}\sub
E|_D$, which as a sheaf is an invertible sheaf of $\cO_D$-modules,
its degree on $D_k$ is $1$ and its degrees on other components are
all $0$, and the induced homomorphisms $\sigma(q_0)$ and
$\sigma(q_{n+m})$ are both non-zero.

\begin{lemm}
\lab{3.15} Let $E$ be as in Proposition  \ref{2.16} and $F\sub E$
be its $\alpha$-destabilizing subsheaf. Then the image subsheaf of
$F|_D\to E|_D$ is a direct sum of subsheaves from the list
$$
\cO_{[0,i)},\ \cO_{(i,n+m]},\ \cO_{[0,n+m]},\ \cO_{[0,n+m]}^{[i]}.
$$
\end{lemm}

\begin{proof}
Since $E|_D$ is regular, all restrictions $F|_{D_k}$ have only
degree $0$ and $1$ factors. We first assume there is a component
$D_k$ so that $F|_{D_k}$ has a factor $\cO_{D_k}(1)$.
Then there is a subchain $D_{[i,j]}$ containing $D_k$ so that this
factor $\cO_{D_k}(1)\sub F|_{D_k}$ extends to a subsheaf $\cL\sub
F|_D$ so that $\cL$ is an invertible sheaf of
$\cO_{D_{[i,j]}}$-modules and the degree of $\cL$ along each
$D_l\sub D_{[i,j]}$ is non-negative. However, if there is another
$l\ne k$ so that the degree of $\cL$ on $D_l$ is $1$, then we can
find a section of $\cL$ that vanishes at $q_i$ and $q_j$. Using
$\cL\to E|_D$, this section induces a section of $E$ that vanishes
on $X^0$, violating the fact that $E|_D$ is regular. Hence $\deg
\cL|_{D_l}=0$ for all $k\ne l\in[i+1,j]$. For the same reason we
conclude that $\sigma(q_i)\mh \cL\otimes \kk(q_i)\to E|_{q_i}$ can
not be zero since otherwise we can find a section of $E$ that
vanishes on $X^0$. Since $E|_D$ is locally free, this is possible
only if $i=0$. For the same reason we have $j=n+m$ and
$\sigma(q_{n+m})\ne 0$. Hence $\cL\sub E|_D$ is a subsheaf of the
type $\cO_{[0,n+m]}^{[k]}\sub E|_D$ described before.

Since $F|_D$ has only degree $0$ and $1$ factors, it is direct to
see that $\cO_{[0,n+m]}^{[k]}\sub F|_D$ must be a direct summand.
Let $F|_D=\cO_{[0,n+m]}^{[k]}\oplus \cF\pri$ be a decomposition.
By repeating this procedure to the sheaf $\cF\pri$, we conclude
that $F|_D$ is a direct sum of sheaves of type
$\cO_{[0,n+m]}^{[k]}$ with a sheaf $\cF$ whose restriction to each
$D_i$ has only degree $0$ factors. We now show that $\cF$ must be
a direct sum of a torsion sheaf and sheaves of the first three
kinds in the list of the lemma: $\cO_{[0,i]}$, $\cO_{(i,n+m]}$ and
$\cO_{[0,n+m]}$. Suppose $\cF$ is not supported at points. Then
there is a subchain $D_{[i,j]}\sub D$ so that $\cO_{D_{[i,j]}}\sub
\cF$ is a subsheaf. Let $\sigma(q_i)$ and $\sigma(q_j)$ be the
induced homomorphisms $\cO_{D_{[i,j]}}\otimes\kk(q_i)\to E|_{q_i}$
and $\cO_{D_{[i,j]}}\otimes\kk(q_j)\to E|_{q_j}$. As before, we
can show that $\sigma(q_i)$ and $\sigma(q_j)$ can not be
simultaneously zero. Further, $\sigma(q_i)\ne 0$ is possible only
if $i=0$. Hence $\cF$ has a factor from the first three kinds in
the above list. By repeating this argument, we conclude that
$F|_D$ is a direct sum of sheaves supported on $q_0$ and
$q_{n+m}$, and sheaves from the list.
\end{proof}

We now complete the proof of Proposition \ref{2.16}

\begin{proof}
We need to show that if $E$ is $\a$-unstable then $V^G$ is
$\a$-unstable as well. Let $F\sub E$ be the $\a$-destabilizing
subsheaf of $E$. We let $R$ and $R\pri$ be the left and the right
rational tails of $X_{n,m}$, let $\tilde E=\rho\sta E$ be the pull
back vector bundle on $X_{n,m}$ and let $T_\leftarrow\sub \tilde
E|_{q_0}$ and $T_\leftarrow\sub \tilde E|_{q_{0}\pri}$ be the
reverse transfer of $\tilde E|_R$ and $\tilde E|_{R\pri}$. We
follow the notation introduced before (\ref{a}) with $q_i$ and
$q_j\pri$ the nodes of $R$ and $R\pri$ so that $q_0=p_1$ and
$q_0\pri=p_2$. Then $\tilde E|_X$ and $V$ fits into the exact
sequence
$$0\lra \tilde E|_X\lra V\lra T_\leftarrow\otimes\kk(p_1)\oplus T\pri_\leftarrow\otimes\kk(p_2)
\lra 0.
$$
Let $\cF\sub \tilde E|_X$ be the image subsheaf of $\rho\sta
F|_X\to \tilde E|_X$, let $\cF\pri\supset\cF$ be the largest
subsheaf of $\tilde E|_X$ so that $\cF\pri/ \cF$ is torsion and
let $F\pri\sub V$ be the largest subsheaf so that $\cF\sub F\pri$
and $F\pri/\cF$ is torsion. Here the inclusion $\cF\sub F\pri$ is
understood in terms of the inclusion of sheaves $\tilde E|_X\sub
V$. Since $\cF\pri$ is the largest possible such subsheaf, it is a
subbundle of $\tilde E|_X$ and $\cF\pri/\cF$ is contained in
$\tilde E|_{p_1}\oplus \tilde E|_{p_2}$. We claim that
$\cF\pri/\cF|_{p_1}\cap T_\leftarrow=\cF\pri/\cF|_{p_2}\cap
T_\leftarrow\pri=\{0\}$. Since otherwise we can find a small
extension of $F\sub E$ so that its $\mu(\cdot,\alpha)$ degree is
larger than $\mu(F,\alpha)$, violating the maximality of the
latter. Combined with the maximality of $\mu(F,\alpha)$, we
conclude that $F\pri/\cF$ is a torsion sheaf supported at $p_1$
and $p_2$ and is indeed a direct sum of a $\kk(p_1)$-module and a
$\kk(p_2)$-module.

We now write the torsion free part of $F|_D$ as a direct sum of
sheaves in the list of Lemma \ref{3.15}. Let $a_-$ (resp. $a_+$)
be the number of summands of type $\cO_{[0,i)}$ with $i\leq n$
(resp. $i> n$) in the decomposition; let $b_-$ (resp. $b_+$) be
the number of $\cO_{(i,n+m]}$ in the summand with $i<n$ (resp.
$i\geq n$); let $c$ be the number of summand $\cO_{[0,n+m]}$ and
let $d_-$ (resp. $d_+$) be the number of summands
$\cO_{[0,n+m]}^{[i]}$ with $i\leq n$ (resp. $i>n$). A direct check
via the construction of $V$ shows that
$$\dim F\pri/\cF|_{p_1}=r_0(F)-(a_++c+d_+)
\quad\text{and}\quad \dim F\pri/\cF|_{p_2}=r_0(F)-(b_-+c+d_-).
$$
Hence
$$\chi(F\pri)=\chi(\cF)+2 r_0(F)-(a_+ + b_-+2c+d_-+d_+).
$$
Further, by a direct check we have $\chi(\cF)=\chi(F)+c$ and
$r\shar(F)=a_++b_-+c+d_-+d_+$. Hence
$$\chi(F\pri)=\chi(F)+2r_0(F)-r\shar(F).
$$
A similar argument shows that $r^+(F\pri)=a_++b_-+c+d_-+d_+$ and
hence $r\shar(F)=r^+(F\pri)$. This implies that
$$\frac{\chi(F\pri)+(1-\a)r^+(F\pri)}{r(F\pri)}=
\frac{\chi(F)-\a r\shar(F)}{r(F\pri)}+2
>\frac{\chi-\alpha r}{r}+2=\frac{\chi\pri+(1-\a)r}{r},
$$
violating the $\a$-stability of $V^G$. This completes the proof of
Proposition \ref{2.16}.
\end{proof}

\begin{coro}\lab{aut}
Each $E\in \cV\ualp\lrc(\XX\shar)$ has $\AUT_0(E)=\CC^{\times}$.
\end{coro}

\begin{proof}
Let $E\in \cV\ualp\lrc(\XX\shar)$ be a sheaf over $X_{n+m}\shar$
and let $(\sigma,f)$ be an automorphism of $E$. Namely, $\sigma\mh
X_{n+m}\to X_{n+m}$ is a based automorphism and $f\mh \sigma\sta
E\cong E$ is an isomorphism. Let $\sigma\pri\mh X_{n,m}\cong
X_{n,m}$ and $f\pri\mh \sigma^{\prime\ast}\tilde E \cong \tilde E$
be the induced isomorphisms. Since $\sigma|_X=\text{id}$, $f$
induces an isomorphism $\tilde f\mh \tilde E|_X\to \tilde E|_X$.
Clearly, $f'$ preserves the subspaces $T_\leftarrow\sub \tilde
E|_{q_0}$ and $T\pri_\leftarrow\sub \tilde E|_{q_0\pri}$. Hence if
we denote by $V\ugp=(V,V^0)$ the associated GPB of $E$, $\tilde f$
extends to an isomorphism $\bar f$ as shown
$$\begin{CD}
0 @>>> \tilde E|_X @>>> V @>>>
T_\leftarrow\otimes \kk(p_1)\oplus T\pri_\leftarrow\otimes \kk(p_2) @>>>0\\
@. @VV{\tilde f}V @VV{\bar f}V @VVV\\
0 @>>> \tilde E|_X @>>> V @>>>
T_\leftarrow\otimes \kk(p_1)\oplus T\pri_\leftarrow\otimes \kk(p_2) @>>>0\\
\end{CD}
$$
We claim that the GPB structure $V^0\sub V|_{p_1}\oplus V|_{p_2}$
is preserved under $\bar f$. Indeed, because $\Gamma_\phi\sub
\tilde E|_{q_n}\oplus \tilde E|_{q_m\pri}$ is the graph of the
tautological isomorphism $\tilde E|_{q_n}\cong \tilde
E|_{q_m\pri}$, and this isomorphism is preserved by $f\pri$,
$\Gamma_\phi$ is preserved by $f\pri$. Next, we look at the
homomorphism
$$\tilde E|_{q_n}\lra \tilde E|_{q_n}/T_\to\mapright{\cong} \tilde E|_{q_0}/T_\leftarrow.
$$
Obviously, the first arrow is canonical. The second arrow is
induced by $\cO_R^{\oplus a}\sub \tilde E|_R$, which is also
preserved by $f\pri$. Hence composite of the above arrows is
invariant under $f\pri$. Therefore, the image of $\Gamma_\phi$ in
$V|_{p_1}\oplus V|_{p_2}$ will be preserved under $f\pri$. This
shows that $(f,\sigma)$ induces an isomorphism $\bar f$ of $V\ugp$
whose restriction to $X^0\sub X$ is exactly $f|_{X^0}$.

Since $E$ is $\alpha$-stable, $V^G$ is $\alpha$-stable by
Proposition \ref{2.16} and hence by the usual argument, the
automorphism group of $V^G$ is $\mathbb{C}^{\times} $. Hence $\bar
f$ is a multiple of the identity map. In particular, after
replacing $f$ by a multiple of itself, we can assume
$f|_{X^0}=\text{id}$. We now show that $\sigma=\text{id}$. As
before, let $q_0,\cdots,q_{n+m}$ be the nodes of $X_{n+m}$. Since
$q_i$ are fixed points of $\sigma$, the isomorphism $f$ induces
automorphisms $f|_{q_i}\mh E|_{q_i}\to E|_{q_i}$. We first claim
that all $f|_{q_i}=\text{id}$. Suppose not, say
$f|_{q_k}\ne\text{id}$. Then $0< k< n+m$ since
$f|_{X^0}=\text{id}$. Since $E|_D$ is non-negative, there is a
section $s\in H^0(D, E)$ so that $s(q_k)$ is not fixed by
$f|_{q_k}$. Hence $s-f^{-1\ast}\sigma\sta s$ is a section of
$H^0(D,E)$ that vanishes on $q_0$ and $q_{n+m}$ but non-zero at
$q_k$. This violates the $\alpha$-stability of $E$. Hence all
$f|_{q_i}$ are the identities.

Next, we claim that $\sigma_k=\sigma|_{D_k}$ are identities for
all $D_k$.
 Indeed,
since $E|_{D_k}\cong \cO_{D_k}^{\oplus r-a}\oplus
\cO_{D_k}(1)^{\oplus a}$ for some $a>0$, there is no isomorphism
of $E|_{D_k}$ with $\sigma_k\sta(E|_{D_k})$ whose restrictions to
$q_{k-1}$ and $q_k$ are the identity maps unless $\sigma_k=id$.
Finally, since $\sigma=\text{id}$ and the restrictions of $f$ to
$X^0$ and all nodes are the identity maps, $f$ must be an identity
map since $E|_{D_k}$ has only degree $0$ and $1$ factors. This
proves that $\AUT_0(E)=\CC^{\times}$.
\end{proof}


The association from $E\in\cV\lrc^\alpha(\XX\shar)$ to
$V^G\in\cG\lrcp^\alpha(X^+)$ constructed above defines a map
\begin{equation}
\lab{4.11} \cV\ualp\lrc(\XX\shar)\lra \cG\ualp\lrcp(X\shar).
\end{equation}
On the other hand, by \cite{Bho} the moduli space of
$\alpha$-stable GPBs $(V,V^0)$ on $X^+$ of rank $r$ and
$\chi(V)=\chi\pri=\chi+r$ form a fine moduli space
$\bG\ualp\lrcp(X^+)$. We now show that the above correspondence
induces a morphism
\begin{equation}
\lab{4.12} \bM\ualp\lrc(\XX\shar) \lra \bG\ualp\lrcp(X^+).
\end{equation}
Later we will show that in case $r=3$  this is the composition of
two blow-ups along smooth subvarieties.

Let $\FF\ualp\lrc(\XX\shar)$ and $\FF\ualp\lrcp(X^+)$ be the
moduli functors of the sets $\cV\ualp\lrc(\XX\shar)$ and
$\cG\ualp\lrcp(X^+)$. To prove the statement, it suffices to show
that the map (\ref{4.11}) defines a transformation of functors
$\FF\ualp\lrc(\XX\shar)\Rightarrow\FF\ualp\lrcp(X^+)$. Namely, to
any scheme $S$ and a family $\cE\in\FF\ualp\lrc(\XX\shar)(S)$ it
associates a unique family $\cV^G\in\FF\ualp\lrcp(X^+)(S)$,
compatible to (\ref{4.11}) and satisfies the base change property.
Let $\cE\in\FF\ualp\lrc(\XX\shar)(S)$ be any family over
$(W_S,q\shar)$. Let $\rho\mh \tilde W_S\to W_S$ be the
normalization along $q\shar(S)\sub W_S$ and let $\pi\mh \tilde
W_S\to X\times S$ the contraction of all rational curves on the
fibers. Let $\bq_-$ and $\bq_+\sub \tilde W_S$ be the two sections
of $\tilde W_S/S$ that are the pre-images of $q\shar(S)$. We then
denote $\bp_i=p_i\times S\sub X\times S$. As usual, we index
$\bq_\pm$ so that $\pi(\bq_-)=\bp_1$. We define
$\cV=\bl\pi\lsta\rho\sta\cE\dual\br\dual$. Because $\cE$ is a
family of $\alpha$-stable sheaves, $R^i\pi\lsta\rho\sta\cE\dual=0$
for $i>0$. Hence $\pi\lsta\rho\sta\cE\dual$ and $\cV$ are locally
free sheaves on $X\times S$. Next, there are canonical
homomorphisms
$$\pi\lsta\rho\sta\cE\dual|_{\bp_1}\lra \rho\sta\cE\dual|_{\bq_-}
\quad\text{and}\quad \pi\lsta\rho\sta\cE\dual|_{\bp_2}\lra
\rho\sta\cE\dual|_{\bq_+}.
$$
Coupled with the identity
$\rho\sta\cE|_{\bq_-}\equiv\rho\sta\cE|_{\bq_+}$, we obtain
homomorphisms
$$\cE|_{q\shar(S)}\mapright{\text{diag}}
\rho\sta\cE|_{\bq_-+\bq_+}\lra
\bl\pi\lsta\rho\sta\cE\dual\br\dual|_{\bp_1+\bp_2} \equiv
\cV|_{\bp_1+\bp_2}.
$$
Here $\rho\sta\cE|_{\bq_-+\bq_+}$ is
$\rho\sta\cE|_{\bq_-}\oplus\rho\sta\cE|_{\bq_+}$, considered as a
sheaf of $\cO_S$-modules. We define $\cV^0$ to be the image sheaf
of the composition of the above arrows. It is direct to check that
this construction $\cE\Rightarrow (\cV,\cV^0)$ satisfies the base
change property. Hence as we argued before (in constructing
$V^G=(V,V^0)$) for each closed $\xi\in S$ the induced
$$\bl\cE|_{q\shar(S)}\br\otimes\kk(\xi)\lra
\cV|_{\bp_1+\bp_2}\otimes\kk(\xi)
$$
is injective. Hence $\cV^0$ is a subvector bundle of
$\cV|_{\bp_1+\bp_2}$. Consequently, by Lemma \ref{2.16} and the
base change property the pair $\cV^G=(\cV,\cV^0)$ is a family of
$\alpha$-stable GPB in $\FF\ualp\lrcp(X^+)$. This defines the
desired transformation of the functors.

\subsection{$\bM\lrc\ualp(\XX\shar)$ as a blow-up of $\bG\lrcp\ualp(X^+)$}

In this subsection, we will restrict ourselves to the case $r=3$
and $\chi=4$. We will show that $\bM_{3,4}\ualp(\XX\shar)$ is a
blow up of $\bG_{3,7}\ualp(X^+)$. We will prove this by looking at
the inverse of (\ref{4.12}):
\begin{equation}
\Phi: \bG_{3,7}\ualp(X^+) -\!-\!\to \bM_{3,4}\ualp(\XX\shar),
\end{equation}
and prove that after two blow-ups of the domain we can resolve the
indeterminacy and the resulting morphism is an isomorphism.

For simplicity, in the remainder of this paper, we will denote
$\bG_{3,7}\ualp(X^+)$ by $\bG\ualp$ and denote
$\bM_{3,4}\ualp(\XX\shar)$ by $\bM\ualp$. By \cite[Theorem
2]{Bho}, we know that there is a universal family of GPBs
$(\cV,\cV^0)$ on $\Ga\times X$, where $\cV$ is a rank 3 vector
bundle over $\Ga\times X$ and $\cV^0$ is a rank 3 subbundle of
$\cV|_{p_1+p_2}\defeq \cV|_{\Ga\times p_1}\oplus \cV|_{\Ga\times
p_2}$ over $\Ga$. On the open dense subset $U\sub \Ga$ where the
induced $\cV^0\to \cV|_{p_1}$ and $\cV^0\to\cV|_{p_2}$ are
isomorphisms, we get a family of vector bundles over $U\times X_0$
by taking the kernel of
$$(1\times\rho)_*\cV\lra (1\times\rho)\lsta \cV|_{q}/\cV^0\equiv  \bl\cV|_{p_1+p_2}\br/\cV^0,
$$
where $\rho\mh X\to X_0$ is the normalization map and $q\in X_0$
is its node. By Lemma \ref{2.16}, the resulting family is a family
of $\alpha$-stable sheaves on $X_0$.

The resulting $U$-family of $\alpha$-stable sheaves defines a
morphism $U\to \Ma$, which defines a rational map $\Ga-\!\to\Ma$,
inverse to the given $\Ma\to \Ga$. We now show how to eliminate
indeterminacy by blowing up the domain $\Ga$ twice.

Let $\bY_i$ (resp. $\bZ_i$) be the subvariety of $\Ga$ consisting
of GPB $(V,V^0)$ such that $V^0\to V|_{p_1}$ (resp. $V^0\to
V|_{p_2}$) have ranks at most $i$. Clearly
$\bY_0\sub\bY_1\sub\bY_2$ and $\bZ_0\sub \bZ_1\sub\bZ_2$ are
chains of subvarieties with $\bY_0$ and $\bZ_0$ smooth. Further,
because $\dim V^0=3$, we know that $\bZ_0\cap\bY_2=\bY_0\cap \bZ_2
=\bY_1\cap \bZ_1=\emptyset$ and $\bY_1$ intersects $\bZ_2$ and
$\bZ_1$ intersects $\bY_2$ transversally. We now blow up $\Ga$
along $\bY_0\cup \bZ_0$. We denote the blow-up of $\Ga$ by $\Ga_1$
and denote the proper transforms of $\bY_i$ and $\bZ_i$ by
$\bY_{i,1}$ and $\bZ_{i,1}$. Because of the intersection property
mentioned, the proper transforms $\bY_{1,1}$ and $\bZ_{1,1}$ are
smooth, satisfying similar intersection properties. We next blow
up $\Ga_1$ along $\bY_{1,1}\cup\bZ_{1,1}$. We denote the blown-up
by $\tilde\bG\ualp$ and denote the corresponding total transforms
by $\tilde \bY_i$ and $\tilde\bZ_i$. This time, all $\tilde \bY_i$
and $\tilde \bZ_i$ are smooth normal crossing divisors.

We now show that $\Ga-\!\to \Ma$ lifts to a morphism $\tga\to
\Ma$. Such morphism will be induced by a family of $\alpha$-stable
sheaves parameterized by  $\tga$. We now construct such a family.
First, we blow up the codimension $2$ subvarieties
$\tilde\bY_{0}\times p_1$ and $\tilde\bZ_{0}\times
p_2\sub\tga\times X$. We denote the resulting family (the
blown-up) by $ W_1$. Let
$$\pi_1:  W_1\to  \tga,\qquad
\Phi_1: W_1\lra\tga\times X \quad\text{and}\quad \phi:\tga\lra \Ga
$$
be the obvious projections. Note that the fibers of $\pi_1$ are
one of $X$, $X_{0,1}$ and $X_{1,0}$. Let $\bq_1\sub  W_1$ (resp.
$\bq_1\pri$) be the proper transform of $\tga\times p_1$ (resp.
$\tga\times p_2$) and let $\bD_1$ (resp. $\bD_1\pri$) be the
exceptional divisor over $\tilde\bY_{0}\times p_1$ (resp.
$\tilde\bZ_{0}\times p_2$).

Next let $(\cV,\cV^0)$ be the universal bundle of $\Ga$ given by a
vector bundle $\cV$ over $\Ga\times X$ and a subbundle $\cV^0$ of
$\cV|_{\bp_1+\bp_2}$, where $\bp_i=\Ga\times p_i\sub\Ga\times X$.
We introduce a new locally free sheaf on $ W_1$:
$$\cV_1\defeq \ker\{\Phi_1\sta\cV\lra \Phi_1\sta\cV|_{\bD_1}\oplus\Phi_1\sta\cV|_{\bD_1\pri}\}.
$$
Because $\tilde\bY_{0}$ is the locus where $V^0\to V|_{p_1}$ are
zeros, and likewise for $\tilde\bZ_{0}$, the pull back
$\phi_1\sta\cV^0\to \Phi_1\sta\cV|_{\bq_1+\bq_1\pri}$ factor
through
\begin{equation}
\lab{5.1} \cV_1^0\defeq \phi\sta\cV^0\lra
\cV_1|_{\bq_1+\bq_1\pri}.
\end{equation}
The pair $(\cV_1,\cV_1^0)$ is a family of GPBs on $(
W_1,\bq_1,\bq_1\pri)$, parameterized by $\tga$. As argued in
\cite{Gie}, for each $\xi\in \tga$ the homomorphisms
$$\cV_1^0|_{ W_{1,\xi}}\lra \cV_1|_{ W_{1,\xi}\cap \bq_1}
\quad\text{and}\quad \cV_1^0|_{ W_{1,\xi}}\lra\cV_1|_{
W_{1,\xi}\cap \bq_1\pri}
$$
have ranks at least $1$.

We now modify this family along the rank $1$ degeneracy loci
$\tilde\bY_1$ and $\tilde\bZ_1$. The construction is similar. We
first blow up $ W_1$ along the disjoint union of
$\pi_1\upmo(\bY_1)\cap\bq_1$ and $\pi_1\upmo(\bZ_1)\cap\bq_1\pri$.
Let $ W_2$ be the blown-up, let $\bq_2$ and $\bq_2\pri$ be the
proper transforms of $\bq_1$ and $\bq_1\pri$ and let $\bD_2$ and
$\bD_2\pri$ be the exceptional divisors of $\Phi_2\mh  W_2\to
W_1$.
As argued in \cite{Gie}, the cokernel
$$\cL_1\defeq\coker\{\cV_1^0|_{\pi_1\upmo(\tilde\bY_{1})}\lra\cV_1|_{\pi_1\upmo(\tilde\bY_{1})\cap\bq_1}\}
$$
is a rank two locally free sheaf on
$\pi_1\upmo(\tilde\bY_{1})\cap\bq_1$. Similarly, let $\cL_1\pri$
be defined with $\tilde\bY_{1}$ replaced by $\tilde\bZ_{1}$ and
with $\bq_1$ replaced by $\bq_1\pri$. It is also a rank two
locally free sheaf on $\pi_1\upmo(\tilde\bZ_{1})\cap\bq_1\pri$.
Similarly to the rank $0$ case, we define $\cV_2$ to be the kernel
of $\Phi_2\sta\cV_1\to \Phi_2\sta \cL_1\oplus\Phi_2\sta\cL_1\pri$.
It is a locally free sheaf over $ W_2$. Further the homomorphism
(\ref{5.1}) induces a homomorphism
\begin{equation}
\lab{5.2} \cV_2^0\defeq \phi\sta\cV^0\lra
\cV_2|_{\bq_2+\bq_2\pri}.
\end{equation}
The pair $(\cV_2,\cV_2^0)$ is a family of GPBs on $ W_2$ over
$\tga$.

Lastly, we resolve the rank $2$ degeneracy. Let $\pi_2\mh
W_2\to\tga$ be the projection and let $ W_3$ be the blow up of $
W_2$ along the disjoint union of
$\pi_2\upmo(\tilde\bY_{2})\cap\bq_2$ and
$\pi_2\upmo(\tilde\bZ_{2})\cap\bq_2\pri$. Let $\bq_3$, $\bq_3\pri$
be the proper transforms of $\bq_2$ and $\bq_2\pri$ and let
$\bD_3$ and $\bD_3\pri$ be the exceptional divisors of $\Phi_3\mh
W_3\to  W_2$.
Again the cokernel $\cL_2$ of
$\cV_2^0|_{\pi_2\upmo(\tilde\bY_{2})}\to\cV_2|_{\pi_2\upmo(\tilde\bY_{2})\cap\bq_2}$
and the similarly defined $\cL_2\pri$ are rank one locally free
sheaves on $\pi_2\upmo(\tilde\bY_2)\cap\bq_2$ and
$\pi_2\upmo(\tilde \bZ_2)\cap \bq_2\pri$. We define $\cV_3$ to be
the kernel of $\Phi_3\sta\cV_2\to
\Phi_3\sta\cL_2\oplus\Phi_3\sta\cL_2\pri$. Then we have the
canonical
\begin{equation}
\lab{5.4} \cV_3^0\defeq \phi\sta\cV^0\lra
\cV_3|_{\bq_3+\bq_3\pri}.
\end{equation}

\begin{lemm}
The family of GPBs $(\cV_3,\cV_3^0)$ on $ W_3$ has the property
that the induced homomorphisms $\cV^0\to \cV_3|_{\bq_3}$ and
$\cV_3^0\to \cV_3|_{\bq_3\pri}$ are both isomorphisms.
\end{lemm}

\begin{proof}
We omit the proof since it is similar to that in \cite{Gie} and
follows directly from the construction.
\end{proof}

For simplicity, we denote $\bq_3$ and $\bq_3\pri\sub W_3$ by $\bq$
and $\bq\pri\sub\tilde W$, and denote $(\cV_3,\cV_3^0)$ by
$(\tilcv,\tilcv^0)$, respectively. We then glue the two sections
$\bq$ and $\bq\pri$ of $\tilde W$ to obtain a family $ W\shar$
over $\tga$ with the marked section $\bq\shar$, the gluing locus.
Clearly, $ W\shar$ is a family of based semistable model of $X_0$
and the tautological projection $\rho\mh\tilde W\to W\shar$ is the
normalization of $ W\shar$ along $\bq\shar$. Over $ W\shar$, we
define $\cE$ via the exact sequence
$$0\lra \cE\lra \rho_*\tilde \cV\lra
\bl\tilde \cV|_{\bq}\oplus \tilde\cV|_{\bq\pri}\br/\tilde
\cV^0\lra 0.
$$

We now show that the family $\cE$ is a family of $\a$-stable
vector bundles over $ W\shar/\tga$. We begin with a closed
$\tilde\xi\in\tga$, with $E\to X_{n+m}\shar$ the restriction of
$\cE\to W\shar$ to the fiber over $\tilde\xi$. We will follow the
notation introduced before. In particular, let $D$ be the chain of
rational curves in $X_{n+m}\shar$ with $q\shar=q_n$ its based
node, let $X_{n,m}$ be its normalization along $q\shar$ and let
$R$ and $R\pri$ be its two rational tails. By our construction,
$E$ is the gluing of a GPB $(\tilde V,\tilde V^0)$ on $X_{n,m}$.
Let $\xi\in\Ga$ be the image of $\tilde\xi$ with $(V,V^0)$ the
corresponding GPB.

\begin{lemm}
The vector bundle $\tilde V|_{R}$ is a regular vector bundle, as
defined in Definition \ref{def}. Further, if we let
$A=\ker\{V^0\to V|_{p_1}\}$, then $\image\{A\to V^0\cong \tilde
V^0 \to \tilde V|_{q_n}\}$ is exactly the transfer $T_\to\sub
\tilde V|_{q_n}$ defined in (\ref{tran}).
\end{lemm}

\begin{proof}
We define the type of the left tail $\tilde V|_{R}$ to be the
triple $(i_0,i_1,i_2)$ defined by $i_j=1$ if $\tilde \xi\in \tilde
\bY_j$ and $i_j=0$ otherwise. Clearly, $n=i_0+i_1+i_3$. We first
study the case $(i_0,i_1,i_2)=(0,1,1)$. Since we only want to
understand $\tilde V|_{R}$, we can assume without loss of
generality that $m=0$. We begin with an explicit description of
the construction of $\tilde V\to X_{2,0}$. First, we let
$B_1\cong\CC^2$ be the cokernel of $V^0\to V|_{p_1}$ and let
$V\pri= \ker\{V\to B_1\otimes\kk(p_1)\}$. The definition of
$V\pri$ induces a canonical filtration $B_1\sub V\pri|_{p_1}$.
Next, let $U_2$ be $\cO_{D_1}(1)^{\oplus 2}\oplus\cO_{D_1}$. Again
the canonical inclusion $\cO_{D_1}(1)^{\oplus 2}\sub U_2$ defines
a filtration $\kk(q_0)^{\oplus 2}\sub U_2|_{q_0}$. We fix an
isomorphism $U_2|_{q_0}\cong V\pri|_{p_1}$ so that it preserves
the two subspaces $\CC^2$ just mentioned. We then define $V_2$ by
the induced exact sequence on $X_{1,0}$:
$$0\lra V_2\lra j\lsta V\pri\oplus j\pri\lsta U_2\lra V\pri|_{p_1}\otimes\kk(q_0)\lra 0.
$$
By our construction of $\cV_2$, $V_2$ is the restriction of
$\cV_2$ to $ W_{2,\tilde\xi}$. The restriction of
$\cV_2^0\to\cV_2|_{\bq_2}$ induces a homomorphism $V^0\to
V_2|_{q_1}$. Since $i_2=1$, its cokernel $B_2$ has dimension $1$.

We define $V_3$ similarly. Let $V_2\pri$ be the kernel of $V_2\to
B_2$. Note that $V_2\pri|_{q_1}$ has a filtration $\CC\sub\CC^3$.
Let $U_3=\cO_{D_2}(1)\oplus \cO_{D_2}^{\oplus 2}$. Then
$U_3|_{q_1}$ also has a filtration $\CC\sub\CC^3$. We fix an
isomorphism $V\pri_2|_{q_1}\cong U_3|_{q_1}$, preserving the two
filtrations. We then define $V_3$ on $X_{2,0}$ by the exact
sequence
$$0\lra V_3\lra \bar j\lsta V_2\oplus \bar j\pri\lsta U_3\lra V_2\pri|_{q_1}\otimes\kk(q_1)\lra 0.
$$
Here $\bar j\mh X_{1,0}\to X_{2,0}$ and $\bar j\pri\mh D_2\to
X_{2,0}$ are the obvious inclusions. Again $V_3$ is the
restriction of $\cV_3$ to the fiber of $ W_3=\tilde W$ over
$\tilde\xi$. Also, the restriction of $\cV_3^0\to\cV_3|_{\bq_3}$
gives us $V^0\to V_3|_{q_2}$ which must be an isomorphism.

We now check that $V_3|_{R_2}$ is regular. First, $V_3|_{D_2}$ has
one degree 1 factor and two trivial factors. We claim that
$V_3|_{D_1}$ also has one degree 1 factor and two trivial factors.
By our construction, this will be true if $\image\{V^0\to
V_2|_{q_1}\}\sub V_2|_{q_1}$ is different from the $\CC^2\sub
V_2|_{q_2}$ induced by the canonical $\cO_{D_1}(1)^{\oplus 2} \sub
V_2|_{D_1}$. Indeed, if they are identical, then $V^0\to V|_{p_1}$
has rank $0$, a contradiction to $i_0=0$. Hence $V_3|_{D_1}$ has
only one degree $1$ factor.

It remains to show that $\dim T_\to\sub \dim V_3|_{q_2}=2$. We
claim that $T_\to=\image\{A\to V^0\to V_3|_{q_2}\}$, where
$A=\ker\{V^0\to V|_{p_1}\}\cong\CC^2$. But this can be checked
directly based on our explicit construction, and will be left to
the readers.

The other cases are trivial except when the type of $\tilde
V|_{R}$ is of type $(1,1,1)$. The study of this case is parallel
to the case studied, and will be omitted. This proves the lemma.
\end{proof}

\begin{lemm}
Let the notation be as before. Then $E$ on $X_{n+m}\shar$ is
$\alpha$-stable.
\end{lemm}

\begin{proof}
We first check that $E|_{D}$ is regular. By the previous lemma,
the restriction of $E$ to each rational curve has only degree $1$
and $0$ factors. Hence to show $E|_D$ is regular we only need to
show that there are no non-trivial section $s\in H^0(E)$ that
vanishes on $X^0\sub X_{n+m}\shar$. Let $s$ be any such section
and let $\tilde s$ be its lift in $H^0(X_{n,m},\tilde V)$. Then
$\tilde s(q_n)\in T_\to$ and $\tilde s(q_m\pri)\in T_\to\pri$. On
the other hand, if we let $v\in V^0$ be the lift of $s(q\shar)$
via the canonical $V^0\cong E|_{q\shar}$, $\tilde s(q_n)=\tilde
s(q_m\pri)=v$, under the canonical $\tilde V|_{q_n}\cong\tilde
V|_{q_m\pri}\cong V^0$. By the previous lemma, $\tilde s(q_n)$
lies in the kernel of $V^0\to V|_{p_1}$ and $\tilde s(q_m\pri)$
lies in the kernel of $V^0\to V|_{p_2}$. This is impossible unless
$v=0$ since $V^0\to V|_{p_1+p_2}$ is injective. This shows that
$s(q\shar)=0$. Then $s=0$ since $\tilde V|_{R_n}$ and $\tilde
V|_{R_m\pri}$ are both regular. This proves that $E|_R$ is
regular.

Once we proved that $E|_D$ is regular, we can apply Lemma
\ref{2.16} to conclude that $E$ is $\alpha$-stable. This completes
the proof.
\end{proof}

\begin{coro}\label{stablem}
The family of locally free sheaves $\cE$ on $ W\shar$ over $\tga$
is a family of $\a$-stable vector bundles.
\end{coro}

As a consequence, we get a morphism $\lambda: \tga\to \Ma$ over
${\bG}\ualp$. Now let $\bU\sub\tga$ be the largest open subset so
that $\lambda|_{\bU}$ is one-one. Since both $\tga$ and $\bM\ualp$
are smooth, $\lambda$ will be an isomorphism if
$\codim(\tga-\bU)\geq 2$. The complement of the 6 divisors $\bY_i,
\bZ_i$ ($i=0,1,2$) represents GPBs $(V,V^0)$ such that $V^0\to
V|_{p_1}$ and $V^0\to V|_{p_2}$ are isomorphisms. Obviously this
is mapped isomorphically by $\lambda$ onto the open subset in
$\bM\ualp$ whose points represent $\alpha$-stable bundles over
$X_0$. By construction, the complement of this open set in $\Ma$
consists of 6 divisors whose generic points are bundles $E$ over
$X_1^\dag$ (2 choices for $q^\dag$) such that the restriction of
$E$ to the rational component is $\cO^a\oplus \cO(1)^{3-a}$ (3
choices for $a$). From our construction of the family of
$\alpha$-stable bundles over $\tga$, it is easy to see that
$\lambda$ maps the 6 divisors of $\tga$ to the 6 divisors of
$\Ma$. Hence the restriction of $\lambda$ to any of the 6 divisors
$\bY_i, \bZ_i$ is generically finite. Since $\lambda$ is injective
on the complement of the 6 divisors, $\lambda$ is a local
homeomorphism at generic points of the divisors and hence
$\lambda$ is injective on an open set $U$ whose complement has
codimension $\ge 2$. Therefore we proved
\begin{coro} $\tga\cong\Ma$.\end{coro}

%

\section{Variations of $\bM\lrd\ualp(\XX\shar)$ in $\alpha$}

The goal of this section is to investigate how
$\bM\lrd\ualp(\XX\shar)$ varies when $\alpha$ varies in $[0,1)$.
Following the work of \cite{Hu, Tha}, it is expected that there is
a finite set $A\sub (0,1)$ so that $\bM\lrd\ualp(\XX\shar)$ is a
constant family when $\alpha$ varies in a connected component of
$(0,1)-A$. Further for $\alpha\in A$ the two moduli spaces
\begin{equation}
\lab{3.4} \bM\lrd^{\alpha^-}(\XX\shar) \leftarrow \!\!- -\!
\!\rightarrow \bM\lrd^{\alpha^+}(\XX\shar)
\end{equation}
are birational and differ by a series of flips. In this section we
will give detailed description of the flips of (\ref{3.4}).

\subsection{The jumping loci}

We begin with determining the set $A$. We continue to assume
$(r,\chi)=1$ and $\chi>r$ throughout this section. Let $\alpha\in
[0,1)-\Lambda_r$ be any real number, let $\delta>0$ be a
sufficiently small number and let $\alpha^\pm=\alpha\pm\delta$. We
suppose $\bM\lrc^{\alpm}(\XX\shar)$ is different from
$\bM\lrc^{\alpp}(\XX\shar)$. Then there is a locally free sheaf
$E$ on $X_n\shar$ that is in $\bM\lrc^{\alpm}(\XX\shar)$ but not
in $\bM\lrc^{\alpp}(\XX\shar)$. Namely, there is a proper subsheaf
$F\sub E$ so that
\begin{equation}
\lab{3.1} \mu_\bd(F,\alpm)<\mu_\bd(E,\alpm) \quad\text{while}\quad
\mu_\bd(F,\alpp)\geq \mu_\bd(E,\alpp).
\end{equation}
Since $E|_D$ is regular, $r_0(F)>0$. Then (\ref{3.1}) implies that
$$\frac{\chi(F)-\alpha r\shar(F)}{r_0(F)}=\frac{\chi(E)-\alpha r}{r}.
$$
This means that $\alpha\in\Lambda_r$. Hence the set $A$ can be
chosen to be $\Lambda_r$. We summerize it as a lemma.

\begin{lemm}
For any two $\alpha_1,\alpha_2$ in a connected component of
$[0,1)-\Lambda_r$, the birational map (\ref{3.4}) is an
isomorphism.
\end{lemm}


It is direct to check that $\Lambda_3=\{1/3, 2/3\}$.

In the remainder of this section, we will restrict ourselves to
the case where $r=3$.
Since for $\alpha\not\in \Lambda_3$ the moduli space
$\Ma_{3,\chi}$ is a blow-up of $\Ga_{3,\chi\pri}$, the previous
lemma suggests the following lemma.

\begin{lemm}
When $\alpha$ varies in a connected component of $[0,1)-\Lambda_3$
the moduli spaces $\Ga_{3,\chi}$ are all isomorphic.
\end{lemm}

\begin{proof}
The proof is straightforward and will be omitted.
\end{proof}

As in \cite{Gie}, it is relatively easy to prove a vanishing
result of the top Chern classes of a certain vector bundle on
$\bM^{1^-}_{3,\chi}$. What we need is the vanishing result on
$\bM^0_{3,\chi}$. One strategy to achieve this is to give an
explicit description of the flips involved in the birational maps
$$ \bM^0_{3,\chi}\leftarrow\!-\!\to \bM^{1/2}_{3,\chi}\leftarrow\!-\!\to \bM^{1^-}_{3,\chi}.
$$
It turns out the two arrows are similar. So we only need to study
the first arrow in detail.

\subsection{Variation of $\bG^{\alpha}\lrcp$}
Since $\bM\ualp\lrc$ is a blow-up of $\bG\ualp\lrcp$, it is
natural to study the variation of $\bG\ualp$ in detail, which we
will do now.

As we will see, we need to study GPB $(V,V^0)$ with $\dim
V^0\ne\rank V$. Here is our convention. We denote by
$\Ga_{r,\chi,a}$ the moduli space of $\a$-stable GPBs $(V,V^0)$ of
rank $r$ vector bundles $V$ with $\chi(V)=\chi$ and
$a$-dimensional subspaces $V^0\subset V|_{p_1+p_2}$. We will still
use $\Ga\lrc$ to denote $\Ga_{r,\chi,r}$, i.e. when $\dim
V^0=\rank V$. Also, in the remainder of this paper we will mostly
interested in the case $r=3$ and $\chi=4$, for convenience we will
abbreviate $\Ga_{3,7}$ to $\Ga$ and abbreviate $\Ma_{3,4}$ to
$\Ma$.

We first investigate how $\Ga\lthreeseven$ varies when $\alpha$
varies. Recall that a GPB $(V, V^0)\in\bG^\alpha_{r,\chi,a}$ on
$X$ is $\alpha$-stable ($\a$-semistable) if for any nontrivial
proper subsheaf $F\sub V$, we have
$$\frac{\chi(F)+(1-\alpha) \dim V^0\cap F|_{p_1+p_2}}
{r(F)}<\frac{\chi(V)+(1-\alpha) a}{r}\ \ \ \ (\le).
$$
Since both sides of the above inequality are linear, if a GPB
$(V,V^0)$ on $X$ is $\alpha_1$-stable but $\alpha_2$-unstable for
some $\alpha_1< \alpha_2$ in the interval $[0,1)$, then we get the
equality
\begin{equation}\lab{ssseq}
\frac{\chi(F)+(1-\alpha) \dim V^0\cap
F|_{p_1+p_2}}{\mathrm{rank}(F)}=\frac{7+(1-\alpha)3 }{3}
\end{equation}
for some $\alpha$ between $\alpha_1$ and $\alpha_2$. It is
elementary to check that the equality can hold only when
$\alpha=1/3$ or $2/3$. Hence $\Ga\lthreeseven$ varies only at
$\alpha=1/3$ and $2/3$ and thus it suffices to consider only the
moduli spaces $\bG^{0^+}\lthreeseven$, $\bG^{1/2}\lthreeseven$ and
$\bG^{1^-}\lthreeseven$.

The variation of $\bG^{\alpha}\lthreeseven$ near $\alpha=1/3$ can
be described as follows: the equation (\ref{ssseq}) holds at
$\alpha=1/3$ only if we have a subbundle $F$ such that
\begin{equation}\lab{gpb2+}
\mathrm{rank}(F)=2,\ \ \ \ \chi(F)=4,\ \ \ \ \dim V^0\cap
F|_{p_1+p_2}=3
\end{equation}
or a subbundle $L$ such that
\begin{equation}\lab{gpb2-}
\mathrm{rank}(L)=1,\ \ \ \ \chi(L)=3,\ \ \ \ \dim V^0\cap
L_{p_1+p_2}=0.
\end{equation}

Suppose a GPB $(V,V^0)$ is $0^+$-stable but $1/2$-unstable. Then
$V$ has a subbundle $L$ satisfying (\ref{gpb2-}).
The quotient bundle $F=V/L$ is equipped with a 3-dimensional
subspace $F^0$ of $F|_{p_1+p_2}$ that is the image of $V^0$. Let
$L^0=0$. Then both GPBs $(L, L^0)$ and $(F,F^0)$ are $1/3$-stable
with the same parabolic slopes. Notice that the
$1/3$-semistability is equivalent to the $1/3$-stability for
$\bG^{1/3}_{2,4,3}$ and $\bG^{1/3}_{1,3,0}$.



We now let $\bA\lot=\bG^{1/3}_{2,4,3}\times \bG^{1/3}_{1,3,0}$.
The previous argument shows that there are maps
$$
\begin{CD}
\bG\uz\lthreeseven-\bG\uo\lthreeseven@>{\eta^-}>> \bA\lot
@<{\eta^+}<< \bG\uo\lthreeseven-\bG\uz\lthreeseven
\end{CD}
$$
that send $(V,V^0)$ to pairs $((F,F^0),(L,L^0))$. We now show that
there are two vector bundles $ W\lot^-$ and $ W\lot^+$ over
$\bA\lot$ so that
\begin{equation}
\lab{fiber} \bG\uz\lthreeseven-\bG\uo\lthreeseven=\PP W\lot^- \and
\bG\uo\lthreeseven-\bG\uz\lthreeseven=\PP W\lot^+.
\end{equation}

Let $(F^G,L^G)=((F,F^0),(L,L^0))\in\bA\lot$ be any pair. Let
$\Ext^1(F^G,L^G)$ be the space of all extensions of GPBs
$$ 0\lra L^G\lra V^G\lra F^G\lra 0.
$$
It is a $\CC$-vector space which fits into the long exact sequence
$$0\lra \Hom(F^G,L^G)\lra \Hom(F,L)\lra \Hom(F^0,L|_{p_1+p_2}/L^0)\lra
$$
$$\quad\lra
\Ext^1(F^G,L^G)\lra \Ext^1(F,L)\lra0.
$$
Thus we have
$$(\eta^-)^{-1}((F^G,L^G))=\bP \Ext^1(F^G,L^G).
$$


Since $L^G$ and $F^G$ are both $1/3$-stable with the same slope,
$\Hom(F^G,L^G)=0$ by a standard argument.
 Hence by the Riemann-Roch theorem, we have
$$\dim\Ext^1(F^G,L^G)=-\chi(\Ext^\cdot(F,L))+6=2g.
$$
(Recall $g(X)=g-1$.) As to the base $\bA\lot$, we have
$$\dim \bG^{1/3}_{2,4,3}=\dim\Ext^1(F^G,F^G)=-\chi(\Ext^\cdot(F,F))+1+3=4g-4
$$
and
$$\dim \bG^{1/3}_{1,3,0}=\dim\Ext^1(L^G,L^G)=-\chi(\Ext^\cdot(L,L))+1=g-1.
$$
Thus $\dim\bA\lot=5g-5$ and
$$\dim \bl\bG\uz\lthreeseven-\bG\uo\lthreeseven\br=(2g-1)+(5g-5)=7g-6.
$$
Following the standard procedure, we pick a universal family
$\cF^G=(\cF,\cF^0)$ over $\bG^{1/3}_{2,4,3}\times X^+$ and a
universal family $\cL^G=(\cL,\cL^0)$ over $\bG^{1/3}_{1,3,0}\times
X^+$. We let $\pi_{ij}$ be the projection from $\bA\times
X=\bG^{1/3}_{2,4,3}\times \bG^{1/3}_{1,3,0}\times X$ to the
product of the $i$-th and the $j$-th factor. We then form the
locally free sheaves
$$ W^-=\Ext^1_{\pi_{12}}(\pi_{13}\sta\cF^G,\pi_{23}\sta\cL^G)
\and
 W^+=\Ext^1_{\pi_{12}}(\pi_{23}\sta\cL^G,\pi_{13}\sta\cF^G).
$$
Note that fibers of $ W\lot^-$ and $ W\lot^+$ over $(F^G,L^G)$ are
exactly $\Ext^1(F^G,L^G)$ and $\Ext^1(L^G,F^G)$, respectively.
Again, as in \cite{Tha1} one shows that
$$\bG\uz\lthreeseven-\bG\uo\lthreeseven=\PP  W^-\lot
\and N_{\PP W\lot^-/\bG\uz}\cong\pi_-\sta W\lot^+\otimes\cO_{\PP
W^-}(-1).
$$
Here we use $N_{A/B}$ to denote the normal bundle of $A\sub B$ and
$\pi_-\mh \PP W^-\to \bA$ is the projection. Notice that
$$\dim \bM\uz\lthreefour=\dim\bG\uz\lthreeseven
=(7g-6)+\dim \Ext^1(L^G,F^G)=9g-8,
$$
(the dimension of $\Ext^1$ is calculated below) which is exactly
the dimension of the moduli of rank three vector bundles over a
genus $g$ curve.

Similarly, the vector space $\Ext^1(L^G,F^G)$ that parameterize
all extensions
$$0\lra F^G\lra V^G\lra L^G\lra 0
$$
satisfies a similar long exact sequence and by the stability of
$F^G$ we have
$$\dim\Ext^1(L^G,F^G)=-\chi(\Ext^\cdot(L,F))=2g-2.
$$
Hence for the same reason,
$$\bG\uo\lthreeseven-\bG\uz\lthreeseven=\PP  W^+\lot
\and N_{\PP W\lot^+/\bG\uo\lthreeseven}\cong \pi_+\sta
W\lot^-\otimes\cO_{\PP W^+}(-1).
$$

Again, following the work of Thaddeus \cite{Tha1}, one checks that
the blow-up of $\bG\uz$ along $\PP W\lot^-$ is isomorphic to the
blow-up of $\bG\uo$ along $\PP W\lot^+$, extending the birational
map $\bG\uz-\to\bG\uo$. Since the details are routine, we omit it
here. This is the explicit description of the flip between
$\bG\uz$ and $\bG\uo$.

\subsection{Flip loci in $\bM^{\alpha}$--First approach}

In this subsection, we will study the flip loci of
$\bM^{0}\sim_{\text{bir}}\bM^{1/2}\lthreefour$ utilizing the fact
that both are moduli of stable vector bundles over $X\shar$.

Let $\Sigma^-\sub \bM^0$ and $\Sigma^+\sub \bM\uo$ be the
indeterminacy loci of the above birational map, namely the
smallest closed subsets so that
\begin{equation}
\lab{3.5} \bM^0-\Sigma^-\mapright{\cong}\bM\uo-\Sigma^+
\end{equation}
is an isomorphism. Our first approach to determine the set
$\Sigma^\pm$ is to characterize all members in $\Sigma^\pm$.

\begin{lemm}
\lab{3.7} Let $E\in\Sigma\upm$ be any member and $F\sub E$ its
$1/3\ump$-destabilizing subsheaf. Then
$$
(r_0(F),r\shar(F),\chi(F))=\left\{
\begin{array}{ll}
(1,0,1) & \text{ if }E\in \Sigma^-\\
(2,3,3) &\text{ if }E\in \Sigma^+
\end{array}\right.
$$
\end{lemm}

\begin{proof}
By Lemma \ref{3.16}, we must have $r_0(F)=1$ or $2$. In the case
$r_0(F)=1$, $r\shar(F)$ and $\chi(F)$ must satisfy the equation
$\chi(F)-r\shar(F)/3=(4-1)/3$, which is possible only if
$r\shar(F)=0$ and $\chi(F)=1$. To determine if $E$ is in
$\Sigma^-$ or $\Sigma^+$, we only need to compute
$$\mu(F,0)\sim 1< 4/3
\quad\text{and}\quad \mu(F,1/2)\sim 1>4/3-1/2,
$$
which implies that $E\not\in \bM\uo$. Thus $E\in \Sigma^-$.

Similarly, when $r_0(F)=2$, the restraint is
$\chi(F)/2-r\shar(F)/6=(4-1)/3$, which has solution $r\shar(F)=3$
and $\chi(F)=3$. A simple computation shows that $E\in\Sigma^+$.
\end{proof}

We now give a more detailed description of pairs $F\sub E\in
\Sigma^+$ which must have $(r_0(F),r\shar(F),\chi(F))=(2,3,3)$. We
assume $E$ is over $X_{n+m}\shar$. As before, let $q_0,
q_1,\cdots$ be the nodes of $X_{n+m}\shar$ with $q_0=p_-$ in $X$.
We define $r(F,q_i)=\dim\{ F|_{q_i}\to E|_{q_i}\}$. Then since
$r_0(F)=2$, we must have $r(F,q_0)$, $r(F,q_n)\leq 2$. On the
other hand, let $F|_D^0=\oplus_1^k \cL_i$ \footnote{$F|_D^0$ is
the torsion free part of $F|_D$} be the decomposition given by
Lemma \ref{3.15}. Then because $r\shar(F)=3$, $k\geq 3$. We claim
that $k=3$. First, because
$$4\geq r(F,q_0)+r(F,q_n)=\sum r(\cL_i,q_0)+r(\cL_i,q_n),
$$
and $r(\cL_i,q_0)+r(\cL_i,q_n)\geq 1$, $k>4$ is impossible. When
$k=4$, by Lemma \ref{3.7}
$$(F|_D)\utf\cong \cO_{[0,i_1)}\oplus \cO_{[0,i_2)}\oplus \cO_{(j_1,n]}\oplus \cO_{(j_2,n]}.
$$
Clearly, this is possible only if $\chi(E|_D)\geq 4+4\chi(\cO_D)$,
which contradicts to the regularity of $E|_D$. Hence $k=3$. In
this case we must have
$$(F|_D)\utf\cong \cO_{[0,i)}\oplus \cO_{(j,n]}\oplus \cO_{[0,n]}
\quad\text{or}\quad \cO_{[0,i)}\oplus \cO_{(j,n]}\oplus
\cO_{[0,n]}^{[k]}.
$$
Now let the marked node of $X_n\shar$ be $q_l$. Since
$r(F,q_l)>r(F,q_0)$ and $r(F,q_n)$, we must have $0<l<n$, thus
$n=2$ or $3$. When $n=2$, $F|_D^0$ must be one of the list
\begin{equation}
\cO_{[0,2)}\oplus\cO_{(0,2]}\oplus \cO_{[0,2]},\
\cO_{[0,2)}\oplus\cO_{(0,2]}\oplus \cO_{[0,2]}^{[1]},\
\cO_{[0,2)}\oplus\cO_{(0,2]}\oplus \cO_{[0,2]}^{[2]}.
\end{equation}
When $n=3$, the $\alpha^\pm$-stable condition on $E$ forces
$E|_{D_i}$ to have at least one degree $1$ factor and combined
there are at most three degree $1$ factors. Hence each $E|_{D_i}$
has exactly one degree $1$ factor. Following this, it is easy to
see that when $q\shar=q_1$, $F|_D^0$ must be one of
\begin{equation}
\cO_{[0,2)}\oplus \cO_{(0,3]}\oplus \cO_{[0,3]}^{[3]},\
\cO_{[0,3)}\oplus \cO_{(0,3]}\oplus \cO_{[0,3]}^{[2]}
\end{equation}
and when $q\shar=q_2$, $F|_D$ must be one of
$$
\cO_{[0,3)}\oplus \cO_{(1,3]}\oplus \cO_{[0,3]}^{[1]},\
\cO_{[0,3)}\oplus \cO_{(0,3]}\oplus \cO_{[0,3]}^{[2]}.
$$

To make our presentation easier to follow, we represent such
subsheaves by graphs. Here is the rule we will follow: for each
$D_k\sub D$, we will encounter invertible subsheaves $\sigma\mh
\cL\to E|_{D_k}$, where $\cL$ is either $\cO_{D_k}$ or
$\cO_{D_k}(1)$. There are two cases, depending on whether the
image sheaf $\sigma(\cL)$ lies in a factor $\cO_{D_k}$ or a factor
$\cO_{D_k}(1)$ of $E|_{D_k}$. In case $\sigma(q_{k-1})=0$ (resp.
$\ne 0$) we will attach a circle (resp. dot) to the left end point
of this line segment. We attach a circle or a dot to the right end
point of the line segment according to whether $\sigma(q_k)=0$ or
$\ne 0$. Following this rule, we will represent a sheaf of
$\cO_D$-modules whose restriction to each $D_k$ is as mentioned by
a chain of line segments, with dots or circles attached. The
following is the list of such subsheaves on $D=D_{[0,2]}$:

\begin{picture}(150,30)(0,0)


\put(12,11){\line(1,0){16}} \put(32,11){\line(1,0){16}}
\put(10,11){\circle{4}} \put(30,11){\circle*{4}}
\put(50,11){\circle*{4}} \put(18,13){\small{${}_{1}$}}

\put(82,11){\line(1,0){16}} \put(102,11){\line(1,0){16}}
\put(80,11){\circle*{4}} \put(100,11){\circle*{4}}
\put(120,11){\circle{4}} \put(108,13){\small{${}_{1}$}}

\put(152,11){\line(1,0){16}} \put(172,11){\line(1,0){16}}
\put(150,11){\circle*{4}} \put(170,11){\circle*{4}}
\put(190,11){\circle*{4}}

\put(222,11){\line(1,0){16}} \put(242,11){\line(1,0){16}}
\put(220,11){\circle*{4}} \put(240,11){\circle*{4}}
\put(260,11){\circle*{4}} \put(228,13){\small{${}_{1}$}}

\put(292,11){\line(1,0){16}} \put(312,11){\line(1,0){16}}
\put(290,11){\circle*{4}} \put(310,11){\circle*{4}}
\put(330,11){\circle*{4}} \put(318,13){\small{${}_{1}$}}
\end{picture}

{\small Figure 1: {\sl These represents subsheaves  $\cO_{(0,2]}$,
$\cO_{[0,2)}$, $\cO_{[0,2]}$, $\cO_{[0,2]}^{[1]}$ and
$\cO_{[0,2]}^{[2]}$.}} \vskip5pt

We next indicate how to represent a pair of sheaves $F\sub E$ near
$D\sub X_2\shar$. Let $U_1$ and $U_2$ be small (analytic) disks
containing $p_1$ and $p_2\in X$ and let $\hat D=U_1\cup D\cup
U_2\sub X_2$ be an analytic neighborhood of $D\sub X_2$. The
following are three examples of pairs $F\sub E$ of subsheaves in
invertible sheaves $E$:

\begin{picture}(150,30)(0,0)


\put(30,16){\line(1,0){8}}
  \put(42,16){\line(1,0){16}}
  \put(62,16){\line(1,0){16}}
  \put(82,16){\line(1,0){8}}
\put(40,16){\circle*{4}}
   \put(60,16){\circle*{4}}
   \put(80,16){\circle*{4}}
\put(48,18){\small{${}_{1}$}}
\put(40,20){\vector(0,-1){2}}

\put(130,16){\qbezier[5](0,0)(4,0)(7,0)}
   \put(142,16){\line(1,0){16}}
   \put(162,16){\line(1,0){16}}
   \put(182,16){\line(1,0){8}}
\put(140,16){\circle{4}}
   \put(160,16){\circle*{4}}
   \put(180,16){\circle*{4}}
\put(148,18){\small{${}_{1}$}}
\put(140,20){\vector(0,-1){2}}

\put(230,16){\line(1,0){8}}
   \put(242,16){\line(1,0){16}}
   \put(262,16){\line(1,0){16}}
   \put(282,16){\qbezier[5](0,0)(4,0)(7,0)}
\put(240,16){\circle*{4}}
   \put(260,16){\circle*{4}}
   \put(280,16){\circle{4}}
\put(268,18){\small{${}_{1}$}} \put(260,20){\vector(0,-1){2}}


\end{picture}

{\small Figure 2: The sheaves $E$ in all three cases are
invertible sheaves of $\cO_{\hat D}$-modules. Its degrees along
$D_1$ and $D_2$ are $1$ and $0$ (abbreviated $\cO_{\hat D}^{[1]}$)
in the first two cases and are $0$ and $1$ (abbreviated $\cO_{\hat
D}^{[2]}$) in the last case. In the first case, $F\cong E$ and
$E/F=0$; In the second case, $F=\cO_{D\cup U_2}$ and $E/F\cong
\cO_{U_1}$; In the last case, $F=\cO_{U_1\cup D}$ and $E/F\cong
\cO_{U_2}$. The three arrows indicate that the marked node of the
first two examples are $q_0$ and of the last example is $q_1$. }
\vskip5pt

Accordingly, a pair $F\sub E$ with $\rank E=3$ along $\hat D$ will
be represented by three horizontal lines, each representing a
direct summand of $E|_{\hat D}$. The following is such an example:

\begin{picture}(150,40)(0,0)


\put(130,6){\line(1,0){8}}
  \put(142,6){\line(1,0){16}}
  \put(162,6){\line(1,0){16}}
  \put(182,6){\line(1,0){8}}
\put(140,6){\circle*{4}}
   \put(160,6){\circle*{4}}
   \put(180,6){\circle*{4}}
\put(148,8){\small{${}_{1}$}}

\put(130,16){\qbezier[5](0,0)(4,0)(7,0)}
   \put(142,16){\line(1,0){16}}
   \put(162,16){\line(1,0){16}}
   \put(182,16){\line(1,0){8}}
\put(140,16){\circle{4}}
   \put(160,16){\circle*{4}}
   \put(180,16){\circle*{4}}
\put(148,18){\small{${}_{1}$}}

\put(130,26){\line(1,0){8}}
   \put(142,26){\line(1,0){16}}
   \put(162,26){\line(1,0){16}}
   \put(182,26){\qbezier[5](0,0)(4,0)(7,0)}
\put(140,26){\circle*{4}}
   \put(160,26){\circle*{4}}
   \put(180,26){\circle{4}}
\put(168,28){\small{${}_{1}$}} \put(160,30){\vector(0,-1){2}}


\end{picture}

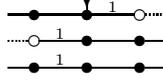
{\small Figure 3: In this example, $E|_{\hat D}$ is a direct sum
of (from top to bottom) $\cO_{\hat D}^{[2]}\oplus\cO_{\hat
D}^{[1]} \oplus \cO_{\hat D}^{[1]}$; The solid lines represent the
subsheaf $F|_{\hat D}\sub E|_{\hat D}$,
which is a direct sum of $\cO_{U_1\cup D}\sub \cO_{\hat D}^{[2]}$,
$\cO_{D\cup U_2}\sub \cO_{\hat D}^{[1]}$ and $\cO_{\hat
D}^{[1]}\sub \cO_{\hat D}^{[1]}$. } \vskip5pt

By analyzing the possible structures of $\{F\sub E\}\in \Sigma^+$
over $X_{n+m}\shar$, we arrive at the following complete lists of
such sheaves:

\begin{picture}(340,40)(0,0)

\put(0,16){{\small $I^{+0}_a$:}} \put(20,6){\line(1,0){8}}
  \put(32,6){\line(1,0){16}}
  \put(52,6){\line(1,0){16}}
  \put(72,6){\line(1,0){8}}
\put(30,6){\circle*{4}}
   \put(50,6){\circle*{4}}
   \put(70,6){\circle*{4}}

\put(20,16){\qbezier[5](0,0)(4,0)(7,0)}
   \put(32,16){\line(1,0){16}}
   \put(52,16){\line(1,0){16}}
   \put(72,16){\line(1,0){8}}
\put(30,16){\circle{4}}
   \put(50,16){\circle*{4}}
   \put(70,16){\circle*{4}}
\put(38,18){\small{${}_{1}$}}

\put(20,26){\line(1,0){8}}
   \put(32,26){\line(1,0){16}}
   \put(52,26){\line(1,0){16}}
   \put(72,26){\qbezier[5](0,0)(4,0)(7,0)}
\put(30,26){\circle*{4}}
   \put(50,26){\circle*{4}}
   \put(70,26){\circle{4}}
\put(58,28){\small{${}_{1}$}} \put(50,30){\vector(0,-1){2}}

\put(110,16){{\small $I^{+1}_a$:}}
 \put(130,6){\line(1,0){8}}
  \put(142,6){\line(1,0){16}}
  \put(162,6){\line(1,0){16}}
  \put(182,6){\line(1,0){16}}
  \put(202,6){\line(1,0){8}}
\put(140,6){\circle*{4}}
   \put(160,6){\circle*{4}}
   \put(180,6){\circle*{4}}
   \put(200,6){\circle*{4}}
\put(168,28){\small{${}_{1}$}}

\put(130,16){\qbezier[5](0,0)(4,0)(7,0)} 
  \put(142,16){\line(1,0){16}}
  \put(162,16){\line(1,0){16}}
  \put(182,16){\line(1,0){16}}
  \put(202,16){\line(1,0){8}}
\put(140,16){\circle{4}}
   \put(160,16){\circle*{4}}
   \put(180,16){\circle*{4}}
   \put(200,16){\circle*{4}}
\put(148,18){\small{${}_{1}$}}

\put(130,26){\line(1,0){8}}
  \put(142,26){\line(1,0){16}}
  \put(162,26){\line(1,0){16}}
  \put(182,26){\qbezier[10](0,0)(4,0)(15,0)}
  \put(202,26){\qbezier[5](0,0)(4,0)(7,0)}
\put(140,26){\circle*{4}}
   \put(160,26){\circle*{4}}
   \put(180,26){\circle{4}}
   \put(200,26){\circle{4}}
\put(188,8){\small{${}_{1}$}} \put(160,30){\vector(0,-1){2}}

\put(240,16){{\small $I^{+2}_a$:}}
 \put(260,6){\line(1,0){8}}
  \put(272,6){\line(1,0){16}}
  \put(292,6){\line(1,0){16}}
  \put(312,6){\line(1,0){16}}
  \put(332,6){\line(1,0){8}}
\put(270,6){\circle*{4}}
   \put(290,6){\circle*{4}}
   \put(310,6){\circle*{4}}
   \put(330,6){\circle*{4}}
\put(278,8){\small{${}_{1}$}}

\put(260,16){\qbezier[5](0,0)(4,0)(7,0)}
  \put(272,16){\qbezier[10](0,0)(4,0)(15,0)}
  \put(292,16){\line(1,0){16}}
  \put(312,16){\line(1,0){16}}
  \put(332,16){\line(1,0){8}}
\put(270,16){\circle{4}}
   \put(290,16){\circle{4}}
   \put(310,16){\circle*{4}}
   \put(330,16){\circle*{4}}
\put(298,18){\small{${}_{1}$}}

\put(260,26){\line(1,0){8}}
  \put(272,26){\line(1,0){16}}
  \put(292,26){\line(1,0){16}}
  \put(312,26){\line(1,0){16}}
  \put(332,26){\qbezier[5](0,0)(4,0)(7,0)}
\put(270,26){\circle*{4}}
   \put(290,26){\circle*{4}}
   \put(310,26){\circle*{4}}
   \put(330,26){\circle{4}}
\put(318,28){\small{${}_{1}$}} \put(310,30){\vector(0,-1){2}}
\end{picture}

\begin{picture}(350,50)(0,0)

\put(0,16){{\small $I^{+0}_b$:}} \put(20,6){\line(1,0){8}}
  \put(32,6){\line(1,0){16}}
  \put(52,6){\line(1,0){16}}
  \put(72,6){\line(1,0){8}}
\put(30,6){\circle*{4}}
   \put(50,6){\circle*{4}}
   \put(70,6){\circle*{4}}
\put(38,8){\small{${}_{1}$}}

\put(20,16){\qbezier[5](0,0)(4,0)(7,0)}
   \put(32,16){\line(1,0){16}}
   \put(52,16){\line(1,0){16}}
   \put(72,16){\line(1,0){8}}
\put(30,16){\circle{4}}
   \put(50,16){\circle*{4}}
   \put(70,16){\circle*{4}}
\put(38,18){\small{${}_{1}$}}

\put(20,26){\line(1,0){8}}
   \put(32,26){\line(1,0){16}}
   \put(52,26){\line(1,0){16}}
   \put(72,26){\qbezier[5](0,0)(4,0)(7,0)}
\put(30,26){\circle*{4}}
   \put(50,26){\circle*{4}}
   \put(70,26){\circle{4}}
\put(58,28){\small{${}_{1}$}} \put(50,30){\vector(0,-1){2}}

\put(110,16){{\small $I^{+1}_b$:}}
 \put(130,6){\line(1,0){8}}
  \put(142,6){\line(1,0){16}}
  \put(162,6){\line(1,0){16}}
  \put(182,6){\line(1,0){16}}
  \put(202,6){\line(1,0){8}}
\put(140,6){\circle*{4}}
   \put(160,6){\circle*{4}}
   \put(180,6){\circle*{4}}
   \put(200,6){\circle*{4}}
\put(168,8){\small{${}_{1}$}}

\put(130,16){\qbezier[5](0,0)(4,0)(7,0)} 
  \put(142,16){\line(1,0){16}}
  \put(162,16){\line(1,0){16}}
  \put(182,16){\line(1,0){16}}
  \put(202,16){\line(1,0){8}}
\put(140,16){\circle{4}}
   \put(160,16){\circle*{4}}
   \put(180,16){\circle*{4}}
   \put(200,16){\circle*{4}}
\put(148,18){\small{${}_{1}$}}

\put(130,26){\line(1,0){8}}
  \put(142,26){\line(1,0){16}}
  \put(162,26){\line(1,0){16}}
  \put(182,26){\line(1,0){16}}
  \put(202,26){\qbezier[5](0,0)(4,0)(7,0)}
\put(140,26){\circle*{4}}
   \put(160,26){\circle*{4}}
   \put(180,26){\circle*{4}}
   \put(200,26){\circle{4}}
\put(188,28){\small{${}_{1}$}} \put(180,30){\vector(0,-1){2}}

\put(240,16){{\small $I^{+2}_b$:}}
 \put(260,6){\line(1,0){8}}
  \put(272,6){\line(1,0){16}}
  \put(292,6){\line(1,0){16}}
  \put(312,6){\line(1,0){16}}
  \put(332,6){\line(1,0){8}}
\put(270,6){\circle*{4}}
   \put(290,6){\circle*{4}}
   \put(310,6){\circle*{4}}
   \put(330,6){\circle*{4}}
\put(278,8){\small{${}_{1}$}}

\put(260,16){\qbezier[5](0,0)(4,0)(7,0)}
  \put(272,16){\qbezier[10](0,0)(4,0)(15,0)}
  \put(292,16){\line(1,0){16}}
  \put(312,16){\line(1,0){16}}
  \put(332,16){\line(1,0){8}}
\put(270,16){\circle{4}}
   \put(290,16){\circle{4}}
   \put(310,16){\circle*{4}}
   \put(330,16){\circle*{4}}
\put(298,18){\small{${}_{1}$}}

\put(260,26){\line(1,0){8}}
  \put(272,26){\line(1,0){16}}
  \put(292,26){\line(1,0){16}}
  \put(312,26){\line(1,0){16}}
  \put(332,26){\qbezier[5](0,0)(4,0)(7,0)}
\put(270,26){\circle*{4}}
   \put(290,26){\circle*{4}}
   \put(310,26){\circle*{4}}
   \put(330,26){\circle{4}}
\put(318,28){\small{${}_{1}$}} \put(310,30){\vector(0,-1){2}}
\end{picture}

\begin{picture}(350,50)(0,0)

\put(0,16){{\small $I^{+0}_c$:}}
\put(20,6){\line(1,0){8}}
  \put(32,6){\line(1,0){16}}
  \put(52,6){\line(1,0){16}}
  \put(72,6){\line(1,0){8}}
\put(30,6){\circle*{4}}
   \put(50,6){\circle*{4}}
   \put(70,6){\circle*{4}}
\put(58,8){\small{${}_{1}$}}

\put(20,16){\qbezier[5](0,0)(4,0)(7,0)}
   \put(32,16){\line(1,0){16}}
   \put(52,16){\line(1,0){16}}
   \put(72,16){\line(1,0){8}}
\put(30,16){\circle{4}}
   \put(50,16){\circle*{4}}
   \put(70,16){\circle*{4}}
\put(38,18){\small{${}_{1}$}}

\put(20,26){\line(1,0){8}}
   \put(32,26){\line(1,0){16}}
   \put(52,26){\line(1,0){16}}
   \put(72,26){\qbezier[5](0,0)(4,0)(7,0)}
\put(30,26){\circle*{4}}
   \put(50,26){\circle*{4}}
   \put(70,26){\circle{4}}
\put(58,28){\small{${}_{1}$}} \put(50,30){\vector(0,-1){2}}

\put(110,16){{\small $I^{+1}_c$:}}
 \put(130,6){\line(1,0){8}}
  \put(142,6){\line(1,0){16}}
  \put(162,6){\line(1,0){16}}
  \put(182,6){\line(1,0){16}}
  \put(202,6){\line(1,0){8}}
\put(140,6){\circle*{4}}
   \put(160,6){\circle*{4}}
   \put(180,6){\circle*{4}}
   \put(200,6){\circle*{4}}
\put(168,8){\small{${}_{1}$}}

\put(130,16){\qbezier[5](0,0)(4,0)(7,0)} 
  \put(142,16){\line(1,0){16}}
  \put(162,16){\line(1,0){16}}
  \put(182,16){\line(1,0){16}}
  \put(202,16){\line(1,0){8}}
\put(140,16){\circle{4}}
   \put(160,16){\circle*{4}}
   \put(180,16){\circle*{4}}
   \put(200,16){\circle*{4}}
\put(148,18){\small{${}_{1}$}}

\put(130,26){\line(1,0){8}}
  \put(142,26){\line(1,0){16}}
  \put(162,26){\line(1,0){16}}
  \put(182,26){\line(1,0){16}}
  \put(202,26){\qbezier[5](0,0)(4,0)(7,0)}
\put(140,26){\circle*{4}}
   \put(160,26){\circle*{4}}
   \put(180,26){\circle*{4}}
   \put(200,26){\circle{4}}
\put(188,28){\small{${}_{1}$}} \put(160,30){\vector(0,-1){2}}

\put(240,16){{\small $I^{+2}_c$:}}
 \put(260,6){\line(1,0){8}}
  \put(272,6){\line(1,0){16}}
  \put(292,6){\line(1,0){16}}
  \put(312,6){\line(1,0){16}}
  \put(332,6){\line(1,0){8}}
\put(270,6){\circle*{4}}
   \put(290,6){\circle*{4}}
   \put(310,6){\circle*{4}}
   \put(330,6){\circle*{4}}
\put(318,8){\small{${}_{1}$}}

\put(260,16){\qbezier[5](0,0)(4,0)(7,0)}
  \put(272,16){\line(1,0){16}}
  \put(292,16){\line(1,0){16}}
  \put(312,16){\line(1,0){16}}
  \put(332,16){\line(1,0){8}}
\put(270,16){\circle{4}}
   \put(290,16){\circle*{4}}
   \put(310,16){\circle*{4}}
   \put(330,16){\circle*{4}}
\put(278,18){\small{${}_{1}$}}

\put(260,26){\line(1,0){8}}
  \put(272,26){\line(1,0){16}}
  \put(292,26){\line(1,0){16}}
  \put(312,26){\qbezier[10](0,0)(4,0)(15,0)}
  \put(332,26){\qbezier[5](0,0)(4,0)(7,0)}
\put(270,26){\circle*{4}}
   \put(290,26){\circle*{4}}
   \put(310,26){\circle{4}}
   \put(330,26){\circle{4}}
\put(298,28){\small{${}_{1}$}} \put(290,30){\vector(0,-1){2}}
\end{picture}

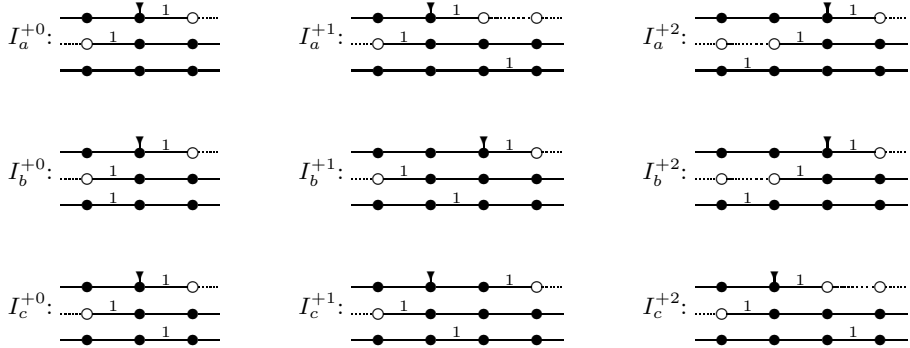
{\small Figure 4: This is the complete list of sheaves in
$\Sigma^+$. Note that $I_c^{+i}$ are the reflections of
$I_b^{+i}$.}

\begin{lemm}
The above is a complete list of sheaves in $\Sigma^+$; The sheaves
of types $I_i^{+j}$ can be (small) deformed to sheaves of type
$I_i^{+0}$; Let $\II_a^+$, $\II_b^+$ and $\II_c^+$ be the set of
sheaves in $\Sigma^+$ of types $I_a^{+\cdot}$, $I_b^{+\cdot}$ and
$I_c^{+\cdot}$, respectively. Then $\II_a^+$, $\II_b^+$ and
$\II_c^+$ are three irreducible components of $\Sigma^+$; Finally,
$\II_b^+\cap \II_c^+=\emptyset$ while $\II_a^+\cap \II_b^+$ (resp.
$\II_a^+\cap \II_c^+$) is the set consisting of sheaves of type
$I_a^{+1}=I_c^{+2}$ (resp. $I_a^{+2}=I_b^{+2}$).
\end{lemm}

\begin{proof}
We have already shown that sheaves of types in the above list can
not be $0$-stable. On the other hand, it is easy to construct
examples of sheaves of each of the type in the list that are
$1/2$-stable. Hence, all types in the list classifies some sheaves
in $\Sigma^+$. It remains to show that this is a complete list. If
$\{F\sub E\}\in\Sigma^+$ is over $X_2\shar$, then we have already
shown that it must be from the list
$\{I_a^{+0},I_b^{+0},I_c^{+0}\}$. The case where $E$ is over
$X_3\shar$ is similar, and will be omitted.

The statement that sheaves in $I_a^{+1}$ can be deformed to
sheaves in $I_a^{+0}$ is straightforward. Let $E$ be such a sheaf,
over $X_3\shar$ with $q\shar=q_1$. Clearly, by smoothing the node
$q_3\in X_3\shar$ we can deform $X_3\shar$ to $X_2\shar$. Then it
is direct to see that we can deform $E$ to $E\pri$ on $X_2\shar$
so that $E\pri$ has type $I_a^{+0}$. The proof of the remaining
statements are similar. We omit the details here since a direct
construction will be given when we study the flips of $\bM\ualp$
later.
\end{proof}

We next give the graphs of all types of sheaves in $\Sigma^-$.
Since the proof is parallel, we will omit it here.

\begin{picture}(340,40)(0,0)

\def\dotline{\qbezier[10](0,0)(8,0)(15,0)}
\def\halfdotline{\qbezier[5](0,0)(4,0)(7,0)}

\put(0,16){{\small $I^{-0}_a$:}} \put(20,26){\line(1,0){8}}
   \put(32,26){\halfdotline}
\put(30,26){\circle{4}}
\put(30,30){\vector(0,-1){2}}

\put(20,16){\halfdotline}
   \put(32,16){\line(1,0){8}}
\put(30,16){\circle{4}}

\put(20,6){\halfdotline}
  \put(32,6){\halfdotline}
\put(30,6){\circle{4}}


\put(70,16){{\small $I^{-1}_a$:}} \put(90,26){\line(1,0){8}}
  \put(102,26){\dotline}
  \put(122,26){\halfdotline}
\put(100,26){\circle{4}}
   \put(120,26){\circle{4}}
\put(100,30){\vector(0,-1){2}}

\put(90,16){\halfdotline} 
  \put(102,16){\line(1,0){16}}
  \put(122,16){\line(1,0){8}}
\put(100,16){\circle{4}}
   \put(120,16){\circle*{4}}
\put(108,18){\small{${}_{1}$}}

\put(90,6){\halfdotline} %
  \put(102,6){\dotline} %
  \put(122,6){\halfdotline} %
\put(100,6){\circle{4}}
   \put(120,6){\circle{4}}

\put(160,16){{\small $I^{-3}_a$:}}
\put(180,26){\halfdotline}
\put(192,26){\dotline}
   \put(212,26){\dotline}
   \put(232,26){\line(1,0){8}}
\put(190,26){\circle{4}}
   \put(210,26){\circle{4}}
   \put(230,26){\circle{4}}
\put(210,30){\vector(0,-1){2}}

\put(180,16){\line(1,0){8}} \put(192,16){\line(1,0){16}}
   \put(212,16){\dotline}
   \put(232,16){\halfdotline}
\put(190,16){\circle*{4}}
   \put(210,16){\circle{4}}
   \put(230,16){\circle{4}}
\put(200,18){\small{${}_{1}$}}

\put(180,6){\halfdotline}
\put(192,6){\dotline}
   \put(212,6){\dotline}
   \put(232,6){\halfdotline}
\put(190,6){\circle{4}}
   \put(210,6){\circle{4}}
   \put(230,6){\circle{4}}
\put(220,8){\small{${}_{1}$}}

\put(260,16){{\small $I^{-4}_a$:}} \put(280,26){\line(1,0){8}}
\put(292,26){\line(1,0){16}}
   \put(312,26){\dotline}
   \put(332,26){\halfdotline}
\put(290,26){\circle*{4}}
   \put(310,26){\circle{4}}
   \put(330,26){\circle{4}}
\put(298,28){\small{${}_{1}$}} \put(310,30){\vector(0,-1){2}}

\put(280,16){\halfdotline}
\put(292,16){\dotline}
   \put(312,16){\line(1,0){16}}
   \put(332,16){\line(1,0){8}}
\put(290,16){\circle{4}}
   \put(310,16){\circle{4}}
   \put(330,16){\circle*{4}}
\put(318,18){\small{${}_{1}$}}

\put(280,6){\halfdotline}
\put(292,6){\dotline}
   \put(312,6){\dotline}
   \put(332,6){\qbezier[5](0,0)(4,0)(8,0)}
\put(290,6){\circle{4}}
   \put(310,6){\circle{4}}
   \put(330,6){\circle{4}}
\end{picture}

\begin{picture}(340,40)(0,0)

\def\dotline{\qbezier[10](0,0)(8,0)(15,0)}
\def\halfdotline{\qbezier[5](0,0)(4,0)(7,0)}

\put(0,16){{\small $I^{-2}_a$:}} \put(20,6){\halfdotline}
  \put(32,6){\dotline}
  \put(52,6){\halfdotline}
\put(30,6){\circle{4}}
   \put(50,6){\circle{4}}
\put(38,8){\small{${}_{1}$}}

\put(20,16){\halfdotline}
   \put(32,16){\dotline}
   \put(52,16){\line(1,0){8}}
\put(30,16){\circle{4}}
   \put(50,16){\circle{4}}

\put(20,26){\line(1,0){8}}
   \put(52,26){\halfdotline}
   \put(32,26){\dotline}
\put(30,26){\circle{4}}
   \put(50,26){\circle{4}}
\put(30,30){\vector(0,-1){2}}


\put(80,16){{\small $I^{-0}_b$:}} \put(100,26){\line(1,0){8}}
  \put(112,26){\dotline}
  \put(132,26){\halfdotline}
\put(110,26){\circle{4}}
   \put(130,26){\circle{4}}
\put(110,30){\vector(0,-1){2}}

\put(100,16){\halfdotline} 
  \put(112,16){\line(1,0){16}}
  \put(132,16){\line(1,0){8}}
\put(110,16){\circle{4}}
   \put(130,16){\circle*{4}}
\put(118,18){\small{${}_{1}$}}

\put(100,6){\halfdotline} %
  \put(112,6){\dotline} %
  \put(132,6){\halfdotline} %
\put(110,6){\circle{4}}
   \put(130,6){\circle{4}}
\put(118,8){\small{${}_{1}$}}

\put(160,16){{\small $I^{-1}_b$:}} \put(180,26){\line(1,0){8}}
\put(192,26){\dotline}
   \put(212,26){\dotline}
   \put(232,26){\halfdotline}
\put(190,26){\circle*{4}}
   \put(210,26){\circle{4}}
   \put(230,26){\circle{4}}
\put(190,30){\vector(0,-1){2}}

\put(180,16){\halfdotline}
\put(192,16){\line(1,0){16}}
   \put(212,16){\line(1,0){16}}
   \put(232,16){\line(1,0){8}}
\put(190,16){\circle{4}}
   \put(210,16){\circle*{4}}
   \put(230,16){\circle*{4}}
\put(200,18){\small{${}_{1}$}}

\put(180,6){\halfdotline}
\put(192,6){\dotline}
   \put(212,6){\dotline}
   \put(232,6){\halfdotline}
\put(190,6){\circle{4}}
   \put(210,6){\circle{4}}
   \put(230,6){\circle{4}}
\put(220,8){\small{${}_{1}$}}

\put(260,16){{\small $I^{-2}_b$:}} \put(280,26){\line(1,0){8}}
\put(292,26){\dotline}
   \put(312,26){\dotline}
   \put(332,26){\halfdotline}
\put(290,26){\circle{4}}
   \put(310,26){\circle{4}}
   \put(330,26){\circle{4}}
\put(290,30){\vector(0,-1){2}}

\put(280,16){\halfdotline}
\put(292,16){\dotline}
   \put(312,16){\line(1,0){16}}
   \put(332,16){\line(1,0){8}}
\put(290,16){\circle{4}}
   \put(310,16){\circle{4}}
   \put(330,16){\circle*{4}}
\put(318,18){\small{${}_{1}$}}

\put(280,6){\halfdotline}
\put(292,6){\dotline}
   \put(312,6){\dotline}
   \put(332,6){\qbezier[5](0,0)(4,0)(8,0)}
\put(290,6){\circle{4}}
   \put(310,6){\circle{4}}
   \put(330,6){\circle{4}}
\put(298,8){\small{${}_{1}$}}
\end{picture}

\begin{picture}(340,40)(0,0)

\def\dotline{\qbezier[10](0,0)(8,0)(15,0)}
\def\halfdotline{\qbezier[5](0,0)(4,0)(7,0)}

\put(0,16){{\small $I^{-3}_b$:}}
\put(20,6){\halfdotline}
  \put(32,6){\dotline}
  \put(52,6){\dotline}
  \put(72,6){\halfdotline}
\put(30,6){\circle{4}}
   \put(50,6){\circle{4}}
   \put(70,6){\circle{4}}
\put(59,8){\small{${}_{1}$}}

\put(20,16){\halfdotline}
   \put(32,16){\dotline}
   \put(52,16){\line(1,0){16}}
   \put(72,16){\line(1,0){8}}
\put(30,16){\circle{4}}
   \put(50,16){\circle{4}}
   \put(70,16){\circle*{4}}
\put(59,18){\small{${}_{1}$}}

\put(20,26){\line(1,0){8}}
   \put(32,26){\line(1,0){16}}
   \put(52,26){\dotline}
   \put(72,26){\halfdotline}
\put(30,26){\circle*{4}}
   \put(50,26){\circle{4}}
   \put(70,26){\circle{4}}
\put(39,28){\small{${}_{1}$}} \put(50,30){\vector(0,-1){2}}

\put(110,16){{\small $I^{-4}_b$:}}
 \put(130,6){\halfdotline}
  \put(142,6){\dotline}
  \put(162,6){\dotline}
  \put(182,6){\dotline}
  \put(202,6){\halfdotline}
\put(140,6){\circle{4}}
   \put(160,6){\circle{4}}
   \put(180,6){\circle{4}}
   \put(200,6){\circle{4}}
\put(169,8){\small{${}_{1}$}}

\put(130,16){\halfdotline} 
  \put(142,16){\dotline}
  \put(162,16){\dotline}
  \put(182,16){\line(1,0){16}}
  \put(202,16){\line(1,0){8}}
\put(140,16){\circle{4}}
   \put(160,16){\circle{4}}
   \put(180,16){\circle{4}}
   \put(200,16){\circle*{4}}
\put(188,18){\small{${}_{1}$}}

\put(130,26){\line(1,0){8}}
  \put(142,26){\line(1,0){16}}
  \put(162,26){\dotline}
  \put(182,26){\dotline}
  \put(202,26){\halfdotline}
\put(140,26){\circle*{4}}
   \put(160,26){\circle{4}}
   \put(180,26){\circle{4}}
   \put(200,26){\circle{4}}
\put(149,28){\small{${}_{1}$}} \put(160,30){\vector(0,-1){2}}


\put(240,16){{\small $I^{-5}_b$:}}
 \put(260,6){\halfdotline}
  \put(272,6){\dotline}
  \put(292,6){\dotline}
  \put(312,6){\dotline}
  \put(332,6){\halfdotline}
\put(270,6){\circle{4}}
   \put(290,6){\circle{4}}
   \put(310,6){\circle{4}}
   \put(330,6){\circle{4}}
\put(318,8){\small{${}_{1}$}}

\put(260,16){\halfdotline}
  \put(272,16){\dotline}
  \put(292,16){\line(1,0){16}}
  \put(312,16){\line(1,0){16}}
  \put(332,16){\line(1,0){8}}
\put(270,16){\circle{4}}
   \put(290,16){\circle{4}}
   \put(310,16){\circle*{4}}
   \put(330,16){\circle*{4}}
\put(298,18){\small{${}_{1}$}}

\put(260,26){\line(1,0){8}}
  \put(272,26){\line(1,0){16}}
  \put(292,26){\dotline}
  \put(312,26){\dotline}
  \put(332,26){\halfdotline}
\put(270,26){\circle*{4}}
   \put(290,26){\circle{4}}
   \put(310,26){\circle{4}}
   \put(330,26){\circle{4}}
\put(279,28){\small{${}_{1}$}}

\put(290,30){\vector(0,-1){2}}
\end{picture}

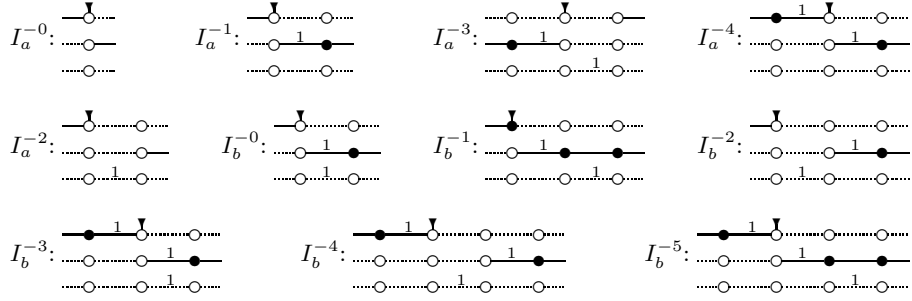
 {\small Figure 5: {\sl These are graphs of some of the sheaves in $\Sigma^-$.}} \vskip5pt

\begin{prop}
The sheaves in $\Sigma^-$ can be divided into subsets:
\newline
1. $I_a^{-i}$ and $I_a^{-i\pri}$ where $i=0,\cdots, 4$; The graph
$I_a^{-i\pri}$ is the reflection of $I_a^{-i}$ along the vertical
axis passing through the arrow;
\newline
2. $I_b^{-i}$ where $i=0,\cdots,5$;
\newline
3. $I_c^{-i}$, where $i=0,\cdots, 5$. Here the graph of $I_c^{-i}$
is the reflection of $I_b^{-i}$ along the axis passing the arrow.

The set $\Sigma^-$ 
is an irreducible subvariety of $\bM\uz$.
\end{prop}

\subsection{Flip loci in $\bM^{\alpha}$--Second approach}

In this subsection we will give an alternative description of the
flip loci $\Sigma^\pm$ of
$\bM\uz\sim_{\text{bir}}\bM^{1/2}\lthreefour$, based on the flip
loci of the moduli of GPBs described before.

First let us describe $\Sigma^+$. We know that $\bM^{1/2}$ is the
result of blowing-up $\bG^{1/2}$ along $\bY_0\cup \bZ_0$ and then
blowing-up the proper transform of $\bY_1\cup \bZ_1$. It is
obvious that $\Sigma^+\lot$ lies in the inverse image of
$\bG^{1/2}-\bG\uz=\PP W^+\lot$. It is easy to see that the
varieties $\bY_0$ and $\bZ_0$ are disjoint from $\PP W^+\lot$.
Thus they do not contribute to the flip loci $\Sigma^+$. Let
$(V,V^0)$ be a non-split extension of $(L,0)\in
\bG^{1/3}_{1,3,0}=M_{1,3}(X)$ by $(F,F^0)\in \bG^{1/3}_{2,4,3}$.
Then $(V,V^0)$ lies in $\bY_1$ (resp. $\bZ_1$) if and only if
$F|_{p_1}\subset F^0$ (resp. $F|_{p_2}\subset F^0$). Let
$$\Xi_{1}=\{\bl(F,F^0),(L,0)\br\mid F|_{p_{2}}\sub F^0\}
\and \Xi_2=\{\bl(F,F^0),(L,0)\br\mid F|_{p_{1}}\subset F^0\},
$$
both are subsets of $\bA$. Since $\bY_1$ is defined to be the loci
where $V^0\to V|_{p_1}$ have dimensions at most one, $\bY_1\cap
\PP W^+\lot$ is the preimage of $\Xi_1\sub\bA$, which has
codimension 2 in $\PP W^+\lot$ and the normal bundle $N_{\bY_1\cap
\PP W^+\lot/\PP W^+\lot}$ is the pull-back of the normal bundle
$N_{\Xi_{1}/\bA\lot}$. A similar statement holds for $\bZ_1\cap
\PP W^+\lot$ as well.
To determine the preimage of $\PP W^+$ in $\bM\uo$, we need to
know the normal bundle $N_{\bY_1/\bG\uo}$, especially its
restriction to $\bY_1\cap \PP W^+$. Now let $\xi\in \bY_1\cap \PP
W^+$ be any point lying over $(F^G,L^G)\in \bA$. It is direct to
see that
\begin{equation}
\lab{3.10} T_\xi \bG\uo/\bl T_\xi\bY_1+ T_\xi\PP W^+\lot\br\cong
\Hom(F|_{p_1},L|_{p_2})
\end{equation}
canonically. Since $\bY_1\cap\PP W^+$ has codimension two in $\PP
W^+$, $\dim T_\xi \bG\uo/T_\xi \bY_1=4$. A similar picture holds
for the intersection $\bZ_1\cap \PP W^+$.

Now let $\bB\lot$ be the blow-up of $\bA\lot$ along
$\Xi_1\cup\Xi_2$ with $\Upsilon_{1}$ and $\Upsilon_2$ the two
exceptional divisors in $\bA\lot$ over $\Xi_{1}$ and $\Xi_2$
respectively. Then the preimage of $\PP W^+$ in $\bM\uo$ is the
union of three smooth irreducible varieties: the first is the
blow-up of $\PP W^+\lot$ along $\PP W^+\lot\cap (\bY_1\cup
\bZ_1)$, which is $\PP W^+\times_\bA\bB$, a projective bundle over
$\bB$; The second is a $\bP^3$-bundles over $\PP
W^+\lot\times_\bA\Xi_1$ and the third is a $\bP^3$-bundles over
$\PP W^+\lot\times_\bA\Xi_2$. We denote these three components by
$\II_a^+$, $\II_b^+$ and $\II_c^+$, respectively. Note that
\begin{equation}
\lab{3.21} \dim \bB\lot=5g-5, \quad\dim \II_a^+=7g-8\and \dim
\II_b^+=\dim \II_c^+=7g-7.
\end{equation}
The intersections are
$$\II_a^+\cap \II_b^+=\II_a^+\times_\bB\Upsilon_1
\and \II_a^+\cap \II_c^+=\II_a^+\times_\bB\Upsilon_1.
$$

We close this subsection by showing that the subsets
$\II_\cdot^\pm$ just defined are exactly the corresponding subsets
described in the previous subsection. Indeed, a general sheaf
$[\cE]\in\II_a$ defined in this section has associated GPB $E^G$
fitting into the exact sequence
$$ 0\lra (F,F^0)\lra (E,E^0)\lra (L,0)\lra 0
$$
so that $F^0\cap F|_{p_1}$ and $F^0\cap F|_{p_2}$ are
1-dimensional. Hence the associated subsheaf $\cF\sub \cE$ has no
$\cO_{\Po}(1)$ when restricts to any rational component in the
base curve of $\cE$. This shows that the two $\II_a$ defined in
the above two subsections are identical. Now consider a general
sheaf $[\cE]$ in the component $\II_b$ defined in this section.
Then since its associated GPB $E^G$ is in $\bY_1$, it still fits
into the above exact sequence with $F|_{p_2}\sub F^0$. In
particular, there is a rational $\Po$ in the base of $\cE$ which
is to the left of the marked node so that $\cF|_{\Po}$ has a
factor $\cO_{\Po}(1)$. This shows that the two definitions of
$\II_b$ are identical.

\subsection{Flipping $\bM\uo$}

The goal of this subsection is to show that we can flip $\bM\uo$
along $\II_a^+$ and then flip the resulting variety along the
proper transform of $\II^+_a\cup\II_c^+$. We will show in the next
section that the resulting variety is isomorphic to $\bM\uz$.

We begin with determining the normal bundles of $\II_\bullet^\pm$.
In the following we adopt the convention that for $S\sub P$ we
denote by $T_SP$ the restriction of the tangent bundle $TP$ to
$S$.

\begin{lemm} \label{normallemma} Let $P,Q$ be smooth subvarieties of a nonsingular
variety $R$ such that $S\defeq P\cap Q$ is smooth. Let
$\pi\mh\tilde{R}\to R$ be the blowing-up  along $Q$ and
$\tilde{P}$ be the proper transform of $P$. Then we have an exact
sequence of vector bundles
$$
0\lra N_{\tilde{P}/\tilde{R}}\lra \pi^*N_{P/R} \lra
\pi^*[T_SR/(T_SP+T_SQ)]\lra 0.
$$
\end{lemm}
\begin{proof} It follows from Lemma 15.4 (i), (iv) in \cite{Ful}.
\end{proof}

Our first application of the lemma is a description of the normal
bundle to $\II_a^+$. We put $R=\bM\uo$, $P=\PP W^+\lot$ and
$Q=\bY_1\cup \bZ_1$. Since $S\defeq P\cap Q$ is $\PP
W^+\times_{\bA}(\Xi_1\cup \Xi_2)$, we have $\tilde P=\PP \tilde
W^+$, where $\tilde W^+$ is the pull back of $W^+$ to $\bB$, where
the latter is the blown-up of $\bA$ along $\Xi_1\cup \Xi_2$. To
determine the normal bundle $N_{\tilde P/\tilde R}$, we need to
find the other two terms in the above exact sequence. The normal
bundle $N_{P/R}\cong \rho^*W^-\lot(-1)$. Also a globalized version
of the isomorphism (\ref{3.10})
shows that the quotient bundle $T_SR/(T_SP+T_SQ)$
is the pull-back of a vector bundle on $\Xi_1\cup\Xi_2$, tensored
with $\cO_{\PP W^+}(-1)$. Thus by Lemma \ref{normallemma}, the
normal bundle $N_{\II_a^+/\bM\uo}$, which is $N_{\tilde P/\tilde
R}$ in the statement of Lemma, becomes the pull-back of a vector
bundle over $\bB\lot$ tensored with $\cO_{\PP \tilde W^+}(-1)$.
Hence by the standard theory in birational geometry we can flip
$\bM^{1/2}$ along $\II_a^+$. Let $\bM_1$ be the result of this
flip, let $\tilde\II_a$ be the flipped loci and let $\tilde\II_b$
and $\tilde\II_c$ be the proper transforms of $\II_b$ and $\II_c$.

Our next step is to show that we can flip $\bM_1$ along
$\tilde\II_b\cup\tilde\II_c$. We begin with a detailed description
of $\tilde\II_b$. As mentioned before, $\II^+_a\cap\II^+_b$ is the
projective bundle $\PP_{\Upsilon_1} W^+$ while $\II^+_b$ is a
$\bP^3$-bundle over $\PP_{\Xi_1} W^+$. To proceed, we need a more
detailed description of $\II^+_b$. Let $\cF^G=(\cF,\cF^0)$ and
$\cL^G=(\cL,0)$ be the restrictions to $\Xi_1$ of the pull-backs
of the universal families of $\bG^{1/3}_{2,4,3}$ and
$\bG^{1/3}_{1,3,0}$, respectively. Since
$\cF_\xi|_{p_2}\sub\cF_\xi^0$ for each $\xi\in\Xi_1$, there is a
tautological line subbundle $\ell\sub \cF|_{p_1}$ so that for each
$\xi$
$$\cF^0_{\xi}=\cF_\xi|_{p_2}\oplus \ell_\xi.
$$
Then the normal bundle $N_{\bY_1/\bG}$ to $\bY_1$ in $\bG\defeq
\bG^{1/2}$, restricting to $\PP_{\Xi_1} W^+\sub \bG$, is
$$N_{\bY_1/\bG}|_{\PP_{\Xi_1} W^+}=
\hom\bl\psi_1\sta \cF|_{p_2}, \psi_1\sta\cF|_{p_1}/\ell\oplus
\psi_1\sta\cL|_{p_1}(-1)\br
$$
where $\psi_1\mh \PP_{\Xi_1}W^+\to \Xi_1$ is the projection and
the sheaf $\psi_1\sta\cL|_{p_1}$ is twisted by $\cO_{\PP W^+}(-1)$
because the universal GPB over $\PP_{\Xi_1}W^+$ is given by
$$0\lra \psi_1\sta\cF^G\lra \cE^G\lra \psi_1\sta\cL^G\otimes\cO_{\PP W^+}(-1)\lra 0.
$$
Therefore,
\begin{equation}
\lab{3.102} \II^+_b=\PP \bl\psi_1\sta \cF|_{p_2}\dual\otimes(
\psi_1\sta\cF|_{p_1}/\ell \oplus \psi_1\sta\cL|_{p_1}(-1))\br,
\end{equation}
as a $\bP^3$-bundle over $\PP_{\Xi_1} W^+$, which itself is a
smooth subvariety of $\bM\uo$.

Based on this description, it is easy to see that
$\II^+_a\cap\II^+_b$ is the sub-bundle
$$\II^+_a\cap \II^+_b=\PP \bl\psi_1\sta (\cF|_{p_2}\dual\otimes (\cF|_{p_1}/\ell))\br
=\PP\bl\cF|_{p_2}\dual\otimes\cF|_{p_1}/\ell\br\times_{\Xi_1}\PP_{\Xi_1}W^+.
$$
Let
$$\begin{CD}
\PP_{\Xi_1}W^+\times_{\Xi_1}\PP(\cF|_{p_2}\dual) @>{\pi_2}>>\PP\cF|_{p_2}\dual\\
@VV{\pi_1}V @VV{\psi_2}V\\
\PP_{\Xi_1}W^+ @>{\psi_1}>> \Xi_1
\end{CD}
$$
be the projections and let $\cK$ be the tautological line
subbundle of $\cF|_{p_2}\dual$ on $\PP \cF|_{p_2}\dual$.
Then after we blow up $\bM\uo$ along $\II^+_a$, the proper
transform $\bE_b$ is the blowing up of $\II^+_b$ along
$\II^+_a\cap\II^+_b$:
\begin{equation}
\lab{3.101} \bE_b=\PP\bl\pi_1\sta\psi_1\sta \cF|_{p_2}\dual\otimes
(\cF|_{p_1}/\ell)
\oplus\pi_2\sta\cK\otimes\pi_2\sta\psi_2\sta\cL|_{p_1}\otimes\pi_1\sta\cO(-1)\br,
\end{equation}
as a $\bP^2$-bundle over
$\PP_{\Xi_1}W^+\times_{\Xi_1}\PP\cF|_{p_2}\dual$. Inside this
$\bP^2$-bundle there is a subbundle
$$\PP\bl\pi_1\sta\psi_1\sta \cF|_{p_2}\dual\otimes (\cF|_{p_1}/\ell)\br
\quad\text{over}\quad
\PP W^+\times_{\Xi_1} \PP \cF|_{p_2}\dual,
$$
which is the intersection $\bE_a\cap\bE_b$. Viewed as a bundle
over $\PP \cF|_{p_2}\dual$, we have
$$\bE_a\cap\bE_b=
\bl\PP W^+\times_{\Xi_1} \PP \cF|_{p_2}\dual\br\times_{\Xi_1}\PP
\cF|_{p_2}\dual.
$$
The proper transform $\tilde\II_b\sub\bM_1$ is then the
contraction of $\bE_a\cap\bE_b$ along all $\PP W^+$ factors, which
is a bundle over $\PP \cF|_{p_2}\dual$. We claim that
$\tilde\II_b$ is a projective bundle over $\PP \cF|_{p_2}\dual$.
Indeed, $\bE_b$, considered as a bundle over $\PP
\cF|_{p_2}\dual$, has a subbundle $\PP_{\Xi_1}W^+\times_{\Xi_1}\PP
\cF|_{p_2}\dual$. Further, using the explicit expression
(\ref{3.101}) its normal bundle in $\bE_b$ along each slice
$$\PP W_\xi^+\times\eta\sub \PP_{\Xi_1}W^+\times_{\Xi_1}\PP \cF|_{p_2}\dual,
$$
where $\eta\in \PP\cF|_{p_2}\dual$ is over $\xi\in\Xi_1$, is
isomorphic to $\cO_{\PP W_\xi^+}(1)^{\oplus 2}$. Hence
$\tilde\II_b$ is a projective bundle over $\PP \cF|_{p_2}\dual$.
By a moment of thought, one sees that
$$\tilde\II_b=\PP W_+\quad
\text{where}\quad W_+=\cK\dual\otimes
\psi_2\sta(\cL|_{p_1}\dual\otimes\cF|_{p_2}\dual\otimes\cF|_{p_1}/\ell)
\oplus\psi_2\sta W^+.
$$

Let $\pi\mh \tilde\II_b=\PP W_+\to \PP \cF|_{p_2}\dual$ be the
projection.

\begin{lemm}
\lab{3.25} There is a vector bundle $W_-$ over $\PP
\cF|_{p_2}\dual$ so that the normal bundle $N_{\tilde\II_b/\bM_1}$
is isomorphic to $\pi\sta W_-\otimes\cO_{\PP W_+}(-1)$. The same
conclusion holds for $\tilde\II_c$ as well.
\end{lemm}

As a corollary, we can flip $\bM_1$ along
$\tilde\II_b\cup\tilde\II_c$ to obtain a new smooth variety
$\bM_2$.

\begin{proof}
We continue to use the notation developed earlier. First, the
exceptional divisor $\PP N_{\bY_1/\bG}$ of the blowing up  of
$\bG=\bG^{1/2}$ is a fiber bundle over $\bY_1$. Because
$\II^+_b=\PP N_{\bY_1/\bG}\times_{\bY_1}(\bY_1\cap\PP W^+)$, we
have the exact sequence of vector bundles
$$0\lra N_{\II^+_b/\PP N_{\bY_1/\bG}}\lra
N_{\II^+_b/\bM^{1/2}}\lra N_{\PP N_{\bY_1/\bG}/\bM^{1/2}}\lra 0
$$
and the identity
$$N_{\II^+_b/\PP N_{\bY_1/\bG}}=\pi\sta N_{\bY\cap\PP W^+/\bY_1},
$$
where
$$\pi: \PP N_{\bY_1/\bG}\times_{\bY_1}(\bY_1\cap\PP W^+)\to \bY_1\cap\PP W^+
$$
is the projection. Let $\eta\in\PP\cF|_{p_2}\dual$ be any point
over $\xi\in\Xi_1$. Then based on the description (\ref{3.102})
$\eta$ defines naturally a subvariety
$$\bP_\eta\defeq \PP\bl 0\oplus \eta\otimes \cL_\xi|_{p_1}(-1)\br
\sub \II^+_b
$$
which is a section of
$$\II^+_b\times_{\Xi_1}\xi\lra \PP_{\Xi_1}W^+\times_{\Xi_1}\xi\defeq \PP W^+_\xi.
$$
Hence $\bP_\eta$ is isomorphic to $\PP W^+_\xi$.

We claim that the normal bundle $N_{\II^+_b/\bM^{1/2}}$
restricting to $\bP_\eta$ is isomorphic to a vector space $V$
tensored by $\cO_{\bP_\eta}(-1)$. Indeed, since the normal bundle
$N_{\PP N_{\bY_1/\bG}/\bM^{1/2}}$ is the tautological sub-line
bundle of the pull back of $N_{\bY_1/\bG}$ over $\PP
N_{\bY_1/\bG}$, its restriction to $\bP_\eta=\PP\bl 0\oplus
\eta\otimes \cL_\xi|_{p_1}(-1)\br$ is isomorphic to
$\cO_{\bP_\eta}(-1)$. As to the term $N_{\bY_1\cap\PP W^+/\bY_1}$,
it is clear that its restriction to $\PP W_\xi^+$ is isomorphic to
$V\otimes \cO(-1)$ for a linear subspace $V\sub
\Ext^1(\cF^G_\xi,\cL^G_\xi)$. Hence $N_{\II^+_b/\PP
N_{\bY_1/\bG}}|_{\bP_\eta}$, and therefore
$N_{\II^+_b/\bM^{1/2}}|_{\bP_\eta}$, are of the forms
$V\pri\otimes \cO_{\bP_\eta}(-1)$ for some vector spaces $V\pri$.

Finally, since the flip loci of $\bM^{1/2}\sim \bM_1$ are away
from $\iota(\PP W_\xi^+)$,
$$N_{\tilde\II_b/\bM_1}|_{\bP_\eta}\cong
N_{\II^+_b/\bM^{1/2}}|_{\bP_\eta}.
$$
By a theorem in \cite{Okoneck}, the restriction of
$N_{\tilde\II_b/\bM_1}$ to the fiber of $\PP W_+$ over
$\eta\in\PP\cF|_{p_2}\dual$ is of the form $W\pri\otimes\cO(-1)$.
Because this is true for all $\eta\in\PP\cF|_{p_2}\dual$, there
must be a vector bundle $W_-$ satisfying the requirement of the
lemma. The case for $\tilde\II_c$ is exactly the same and will be
omitted.
\end{proof}

\section{The isomorphism of the two flips}

The goal of this section is to prove the following.

\begin{prop}
\lab{3.17}
 The birational $\bM\uz\sim \bM\uo$ induces an
isomorphism $\bM^0\cong\bM_2$.
\end{prop}

We first briefly outline the strategy. As argued in section 3, we
can flip $\bM\uhalf$ along $\II_a^+$ to obtain a new variety
$\bM_1$. Let $\tilde\II_a$, etc., be the flipped loci of
$\II^+_a$, etc. Then Lemma \ref{3.25} tells us that we can flip
$\bM_1$ again along $\tilde\II_b$ and $\tilde\II_c$ to obtain a
new variety $\bM_2$. The key is to show that the birational map
$\bM_2\sim\bM\uz$ extends to a morphism $\bM_2\to \bM^0$, because
it is an isomorphism away from a subset of codimension at least
two, it is an isomorphism.

We now list the main steps in constructing the morphism
$\bM_2\to\bM\uz$. We first blow up $\bM\uhalf$ along $\II_a^+$ to
obtain the variety $\tilde\bM_1$, and then contract its
exceptional divisor to get the flip $\bM_1$. Let $\tilde\II_b$ and
$\tilde\II_c$ be the flipped loci of $\II_b^+$ and $\II_c^+$. We
then blow up $\bM_1$ along $\tilde\II_a\cup\tilde\II_c$ and then
contract the exceptional divisor to get the second flip $\bM_2$,
as shown below. It is easy to see that if we blow up $\tilde\bM_1$
along $\pi_2\upmo(\tilde\II_b \cup\tilde\II_c)$, the resulting
variety $\tilde\bM$ fits into the diagram below. The first main
technical part of the proof is to show that the birational maps
extend to a morphism $\phi\mh\tilde\bM\to\bM^0$. This is achieved
by first picking a (local) universal family of $\bM\uhalf$ and
then performing an elementary modification to it to obtain a
family on $\tilde\bM_1$, and then performing another elementary
modification to the new family to get a family over $\tilde\bM$.
We will show that the latter is a family of $0^+$-stable vector
bundles, and thus induces a morphism $\Psi$, extending the
birational map.

\begin{picture}(350,90)(20,0)

\put(50,10){{$\bM^{1/2}$}} \put(170,10){{$\bM_1$}}
\put(290,10){{$\bM_2$}}

\put(95,33){\vector(-2,-1){20}} \put(80,33){$\Psi_1$}
\put(135,33){\vector(2,-1){20}} \put(140,33){$\Psi_2$}
\put(215,33){\vector(-2,-1){20}}\put(200,33){$\Psi_3$}
\put(255,33){\vector(2,-1){20}} \put(260,33){$\Psi_4$}

\put(110,40){$\tilde\bM_1$} \put(230,40){$\tilde\bM_2$}

\put(155,63){\vector(-2,-1){20}} \put(140,63){$\tilde\Psi_1$}
\put(195,63){\vector(2,-1){20}}  \put(200,63){$\tilde\Psi_2$}

\put(170,70){$\tilde\bM$} \put(290,70){$\bM^0$}

\put(220,73){\vector(1,0){35}} \put(233,76){$\Psi$}

\end{picture}

\noindent In the end, we will show that $\Psi$ descends to a
morphism $\bM_2\to\bM^0$, as desired.

\subsection{The family over $\tilde\bM_1$}

Our first step is to construct the (local) tautological family
over $\tilde\bM_1$. Let $\xi\in\II_a^+\cap \II_b^+$ be any point
and let $\bU\sub\bM\uo-\II_c$ be an open subset containing $\xi$.
Without loss of generality, we can assume that the moduli space
$\bM\uhalf$ admits a universal family $\EE$ over $\bU$ that is a
sheaf over a family of nodal curves $\cW$ over $\bU$. The desired
family over $\tilde\bM$ will be the result of an elementary
modification to the pull back of $\EE$ to $\tilde\bM$.

To this end, we first need to construct the associated family
$\tilde\cW$ and $\cW\st$ over $\tilde\bM_1$. Let $\tilde\bM_1$ be
the blow up of $\bM^{1/2}$ along $\II^+_a$, let
$\tilde\bU\sub\tilde\bM_1$ be the pre-image of $\bU\sub\bM^{1/2}$,
and let $\bar\cW=\cW\times_{\bU}\tilde\bU$ be the pull back
family. The family $\bar\cW$ is singular, and hence needs to be
smoothed first. We now set up the notation for the singular loci
of the fibers of $\cW/\bU$. Let $\xi\in \bU$ and $\cW/\bU$ be as
before, and let $\cN\sub\cW$ be the singular loci of the fibers of
$\cW/\bU$. Then $\cN\cap\cW_\xi$ consists of four nodes: $q_0$,
$q_1$, $q\shar=q_2$ and $q_3$  of $\cW_\xi$. By shrinking $\bU$ if
necessary, we can assume that $\cN$ is a disjoint union of four
varieties $\cN_0,\cdots,\cN_3$ indexed so that $\cN_i\cap
\cW_\xi=q_i$. Clearly, $\cN_2$ is a section over $\bU$ while all
others are codimension two smooth subvarieties of $\cW$. Let
$\bD_i\sub \bU$ be the image of $\cN_i$ under the projection
$\cW\to \bU$. Then $\II^+_a\cap\bU\sub \bD_1\cap \bD_3$ while
$(\II^+_a\cap\II^+_b)\cap\bU=\II^+_a\cap\bD_0$.

Next, let $\bE_a\sub\tilde\bM_1$ be the exceptional divisor of
$\Psi_1$, let $\tilde\bD_i\sub\tilde\bU$ be the proper transform
of $\bD_i$ and let
$\bar\cN_i=\cN_i\times_{\bU}\tilde\bU\sub\bar\cW$ be the
associated subscheme. Since $\bD_0$ intersects transversally with
$\II^+_a$, the total family $\bar \cW$ is smooth along
$\bar\cN_0$. Obviously, $\bar\cN_2$ remains a section over
$\tilde\bU$. As to $\bar\cN_1$ and $\bar\cN_3$, first of all,
$\bar\cW$ remains smooth along $\bar\cN_1$ (resp. $\bar\cN_3$)
away from $\bar\cP_1\defeq \bar\cN_1|_{\bE_a\cap\tilde\bD_1}$
(resp. $\bar\cP_3\defeq \bar\cN_3|_{\bE_a\cap\tilde\bD_3}$).
Secondly, the normal slice to $\bar\cP_1$ and $\bar\cP_3$ in
$\bar\cW$ are isomorphic to the singularity of $z_1z_2=z_3z_4$.
Hence we can find a small resolution of the singularities of $\bar
\cW$ to obtain a new family $\tilde\cW$. It is known that the
small resolution is obtained by first blowing up the singular loci
of $\bar\cW$ and then contract one $\Po$ factor of the exceptional
divisors\footnote{See \cite[\S1]{Li} for details.}. Since there
are two $\Po$ factors, to proceed we need to specify our cohice of
contraction. Let $\tilde\cP_1$ and $\tilde\cP_3$ be the
exceptional loci of $\tilde\cW\to\bar\cW$, which are $\Po$-bundles
over $\bar\cP_1$ and $\bar\cP_3$ respectively. We next pick a lift
$\eta\in \bE_a\cap\tilde\bD_1\cap\tilde\bD_3$ of
$\xi\in\II_a^+\cap\II_b^+$ and consider the fiber
$\tilde\cW_{\eta}$ of $\tilde\cW$ over $\eta$. As it stands, it
contains five rational curves, indexed by $R_1, R_{1+}, R_2, R_3$
and $R_{3+}$ so that the first intersects with $X$ at $p_1$ while
any two consecutive $R_{\bullet}$'s insect at one point. The small
resolution is the one so that $R_{1+}=\tilde\cW_{\eta} \cap
\tilde\cP_1$ and $R_{3+}=\tilde\cW_{\eta}\cap\tilde \cP_3$, and
that if we let $S\sub \tilde \bD_0\cap\tilde\bD_1\cap\tilde\bD_3$
be a smooth curve that contains $\eta$ and is transversal to
$\bE_a$, then the family $\tilde\cW_S=\tilde\cW
\times_{\tilde\bU}S$
smooth the nodes $q_{1+}=R_{1+}\cap R_2$ and $q_{3-}=R_3\cap
R_{3+}$ of the central fiber $\tilde\cW_\eta$. (See the figure
below.)

\begin{picture}(350,110)(-30,0)

\put(50,87){\line(0,1){15}}

\put(46,75){\line(1,3){4}}
\put(46,75){\qbezier[10](0,0)(2,-6)(4,-12)}
\put(50,49){\line(0,1){14}} \put(46,37){\line(1,3){4}}
\put(46,37){\qbezier[10](0,0)(2,-6)(4,-12)}
\put(50,10){\line(0,1){15}}

\put(50,87){\line(1,0){40}} \put(46,75){\line(1,0){40}}
\put(50,49){\line(1,0){40}} \put(50,25){\line(1,0){40}}

\put(40,87){{\small ${}_{q_0}$}} \put(36,75){{\small ${}_{q_1}$}}
\put(38,63){{\small ${}_{q_{1\!+}}$}} \put(40,49){{\small
${}_{q_2}$}} \put(34,37){{\small ${}_{q_{3\!-}}$}}
\put(40,25){{\small ${}_{q_3}$}}


\put(200,87){\line(0,1){15}}

\put(196,75){\line(1,3){4}}
\put(196,75){\qbezier[10](0,0)(2,-6)(4,-12)}
\put(200,49){\line(0,1){14}} \put(196,37){\line(1,3){4}}
\put(196,37){\qbezier[10](0,0)(2,-6)(4,-12)}
\put(200,10){\line(0,1){15}}

\put(200,63){\line(1,0){40}} \put(200,49){\line(1,0){40}}
\put(196,37){\line(1,0){40}}

\put(190,87){{\small ${}_{q_0}$}} \put(186,75){{\small
${}_{q_1}$}} \put(188,63){{\small ${}_{q_{1\!+}}$}}
\put(190,49){{\small ${}_{q_2}$}} \put(184,37){{\small
${}_{q_{3\!-}}$}} \put(190,25){{\small ${}_{q_3}$}}


\end{picture}

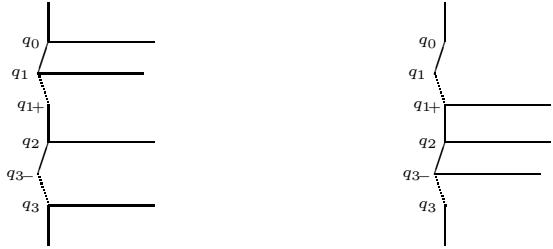
{\small Figure 6: The left one represents the total family over
$S$: the vertical chain of lines is the central fiber
$\tilde\cW_\eta$ with each corner represents a node, as labelled,
the top and the bottom lines represent the main component $X$
while others are rational curves; The two dotted lines represent
the two $\Po$ that are contracted under
$\tilde\cW_\eta\to\cW_\eta$, and the horizontal lines show that
the associated nodes are not smoothed in the family $S$.  The
right figure represents the total space over a curve $\eta\in
S\pri\sub \tilde\bD_1\cap\tilde\bD_3$, $S\pri\sub\bE_a$ and is
transversal to $\tilde\bD_0$. } \vskip5pt

The family $\cW\st$ is a contraction of $\tilde\cW$. Let
$\tilde\cR_2$ and $\tilde\cR_3$ be the two irreducible components
of $\tilde\cW\times_{\tilde\cU}{\bE_a}$ that contain $R_2$ and
$R_3\sub\tilde\cW_{\eta}$ respectively. It is easy to see that
each is isomorphic to $(\bE_a\cap\tilde\bU)\times\Po$ and its
normal bundle in $\tilde\cW$ has degree $-1$ along its fibers.
Therefore, we can contract $\tilde\cW$ along $\tilde\cR_2$ and
$\tilde\cR_3$ to obtain a new family of nodal curves. We denote
this new family by $\cW\st$ with projection
$$\pr: \tilde\cW\lra \cW.
$$

Next we investigate the associated families of sheaves over
$\tilde\cW$ and $\cW\st$. Let $\EE$ be the universal sheaf over
$\cW$. By our description of the sheaves in $\II^+_a$, for each
$\zeta\in \II^+_a\cap\bU$ the sheaf
$\EE_\zeta=\EE\otimes_{\cO_{\cW}}\cO_{\cW_\zeta}$ has a canonical
subsheaf $\FF_\zeta\sub\EE_\zeta$ and the associated quotient
sheaf $\LL_\zeta=\EE_\zeta/\FF_\zeta$. Let $\cZ_\zeta\sub
\cW_\zeta$ be the support of $\LL_\zeta$. Then $\LL_\zeta$ is a
rank one locally free sheaf of $\cO_{\cZ_\zeta}$-modules. Further,
it is direct to check that the union $\cup_{\zeta\in\bU\cap
\II^+_a}\cZ_\zeta$ forms a smooth subvariety of $\cW$ and is the
irreducible component of $\cW\times_{\bU}(\II^+_a\cap \bU)$ that
contains $X\times (\II_a^+\cap\bU)$. We denote this by $\cZ$ with
inclusion $\iota\mh \cZ\sub \cW$. Further, there is a locally free
sheaf $\LL$ of $\cO_{\cZ}$-modules and a quotient sheaf
homomorphism $\EE\to \iota\lsta \LL$ so that its restriction to
each fiber $\cW_\zeta$ is exactly the pair $\EE_\zeta\to\LL_\zeta$
mentioned before.

Now we are ready to perform an elementary modification on the pull
back sheaf over $\tilde\cW$. Let $\tilde\cZ\sub\tilde\cW$ be the
pre-image of $\cZ\sub\cW$ under the projection $\tilde\cW\to\cW$.
By our choice of the small resolution, $\tilde\cZ$ is a smooth
divisor of $\tilde\cW$ and the total space
$\tilde\cW\times_{\tilde\bU}(\bE_a\cap\tilde\bU)$ is a union of
three irreducible components: $\tilde\cR_2$, $\tilde\cR_3$ and
$\tilde\cZ$. We consider the pull-back family
$\tilde\EE\defeq\tilde\pr\sta\EE$ and the associated surjective
homomorphism $\tilde\EE\to \tilde\LL\defeq\tilde\pr\sta\iota\lsta
\LL$. Let $\tilde\EE\pri$ be the kernel of this homomorphism.

\begin{lemm}\lab{4.2}
The sheaf $\tilde\EE\pri$ is a locally free sheaf of
$\cO_{\tilde\cW}$-modules. Further, for any
$\eta\in\bE_a\cap\tilde\bU$, the restriction of $\tilde\EE\pri$ to
$\tilde\cR_2\times_{\tilde\bU}\eta$ and
$\tilde\cR_3\times_{\tilde\bU}\eta \cong\Po$ is isomorphic to
$\cO_{\Po}^{\oplus 3}$.
\end{lemm}

\begin{proof}
Since $\tilde\cZ\sub\tilde\cW$ is a smooth divisor and $\tilde\LL$
is an invertible sheaf of $\cO_{\tilde\cZ}$-modules, the kernel of
$\tilde\pr\sta\EE\to \tilde \LL$ is locally free.

We now prove the second part. We first consider the case
$\eta\in\bE_a\cap\tilde\bD_1\cap\tilde\bD_3$ in details. Let
$S\sub \tilde\bD_1\cap\tilde\bD_3$ be a smooth curve that contains
$\eta$ and is transversal to $\bE_a$. Since $\tilde\bD_0$ is
transversal to $\bE_a$, we can assume $S\sub\tilde\bD_0$. (See
Figure 6.) Then the irreducible component $V_1$ of $\tilde\cW_S$
that contains $R_{1+}$  and $R_2$ as $(-1)$-curves. Similarly,
$R_3$ and $R_{3+}$ are $(-1)$-curves in the irreducible component
$V_2$ as shown in figure 6. Now let $\tilde\EE_S\defeq
\tilde\EE\otimes_{\cO_{\tilde\cW}}\cO_{\tilde\cW_S}$ be the
pull-back family. As before, we denote by
$\tilde\EE_\eta\to\tilde\LL_\eta$ the restriction of
$\tilde\EE\to\tilde\LL$ to $\tilde\cW_\eta$. Since $\bE_a$ is a
smooth divisor, the sheaf
$\tilde\EE\pri_\eta\defeq\tilde\EE\pri|_{\tilde\cW_\eta}$ is
canonically isomorphic to $\ker\{\tilde\EE_S\to
\tilde\LL_\eta\}|_{\tilde\cW_\eta}$. Following our convention, the
pair $\tilde\EE_\eta\to\tilde\LL_\eta$ can be represented by the
left graph in Figure 7 below.

\begin{picture}(350,50)(0,0)

\def\dotline{\qbezier[10](0,0)(8,0)(15,0)}
\def\halfdotline{\qbezier[5](0,0)(4,0)(7,0)}

\put(-10,6){\line(1,0){8}}
  \put(2,6){\line(1,0){16}}
  \put(22,6){\line(1,0){16}}
  \put(42,6){\line(1,0){16}}
  \put(62,6){\line(1,0){16}}
  \put(82,6){\line(1,0){16}}
  \put(102,6){\line(1,0){8}}
\put(0,6){\circle*{4}} \put(20,6){\circle*{4}}
   \put(40,6){\circle*{4}}
   \put(60,6){\circle*{4}}
   \put(80,6){\circle*{4}}
   \put(100,6){\circle*{4}}
\put(8,8){\small{${}_{1}$}}

\put(-10,16){\halfdotline}
  \put(2,16){\dotline}
  \put(22,16){\dotline}
  \put(42,16){\line(1,0){16}}
  \put(62,16){\line(1,0){16}}
  \put(82,16){\line(1,0){16}}
  \put(102,16){\line(1,0){8}}
\put(0,16){\circle{4}} \put(20,16){\circle{4}}
   \put(40,16){\circle{4}}
   \put(60,16){\circle*{4}}
   \put(80,16){\circle*{4}}
   \put(100,16){\circle*{4}}
\put(48,18){\small{${}_{1}$}}

\put(-10,26){\line(1,0){8}}
  \put(2,26){\line(1,0){16}}
  \put(22,26){\line(1,0){16}}
  \put(42,26){\line(1,0){16}}
  \put(62,26){\line(1,0){16}}
  \put(82,26){\dotline}
  \put(102,26){\halfdotline}
\put(0,26){\circle*{4}} \put(20,26){\circle*{4}}
   \put(40,26){\circle*{4}}
   \put(60,26){\circle*{4}}
   \put(80,26){\circle{4}}
   \put(100,26){\circle{4}}
\put(68,28){\small{${}_{1}$}}

\put(60,30){\vector(0,-1){2}}

 \put(125,6){\qbezier[5](0,0)(4,0)(7,0)}
  \put(137,6){\qbezier[10](0,0)(8,0)(15,0)}
  \put(157,6){\qbezier[10](0,0)(8,0)(15,0)}
  \put(177,6){\qbezier[10](0,0)(8,0)(15,0)}
  \put(197,6){\qbezier[10](0,0)(8,0)(15,0)}
  \put(217,6){\qbezier[10](0,0)(8,0)(15,0)}
  \put(237,6){\qbezier[5](0,0)(4,0)(7,0)}
\put(135,6){\circle{4}}
   \put(155,6){\circle{4}}
   \put(175,6){\circle{4}}
   \put(195,6){\circle{4}}
   \put(215,6){\circle{4}}
   \put(235,6){\circle{4}}
\put(143,8){\small{${}_{1}$}}

\put(125,16){\line(1,0){8}}
  \put(137,16){\line(1,0){16}}
  \put(157,16){\line(1,0){16}}
  \put(177,16){\qbezier[10](0,0)(8,0)(15,0)}
  \put(197,16){\qbezier[10](0,0)(8,0)(15,0)}
  \put(217,16){\qbezier[10](0,0)(8,0)(15,0)}
  \put(237,16){\qbezier[5](0,0)(4,0)(7,0)}
\put(135,16){\circle*{4}}
   \put(155,16){\circle*{4}}
   \put(175,16){\circle{4}}
   \put(195,16){\circle{4}}
   \put(215,16){\circle{4}}
   \put(235,16){\circle{4}}
\put(163,18){\small{${}_{1}$}}

\put(125,26){\qbezier[5](0,0)(4,0)(7,0)}
  \put(137,26){\qbezier[10](0,0)(8,0)(15,0)}
  \put(157,26){\qbezier[10](0,0)(8,0)(15,0)}
  \put(177,26){\qbezier[10](0,0)(8,0)(15,0)}
  \put(197,26){\qbezier[10](0,0)(8,0)(15,0)}
  \put(217,26){\line(1,0){16}}
  \put(237,26){\line(1,0){8}}
\put(135,26){\circle{4}}
   \put(155,26){\circle{4}}
   \put(175,26){\circle{4}}
   \put(195,26){\circle{4}}
   \put(215,26){\circle*{4}}
   \put(235,26){\circle*{4}}
\put(223,28){\small{${}_{1}$}}

\put(195,30){\vector(0,-1){2}}

\put(260,6){\halfdotline}
  \put(272,6){\dotline}
  \put(292,6){\dotline}
  \put(312,6){\dotline}
  \put(332,6){\halfdotline}
\put(270,6){\circle{4}}
   \put(290,6){\circle{4}}
   \put(310,6){\circle{4}}
   \put(330,6){\circle{4}}
\put(278,8){\small{${}_{1}$}}

\put(260,16){\line(1,0){8}}
  \put(272,16){\line(1,0){16}}
  \put(292,16){\line(1,0){16}}
  \put(312,16){\dotline}
  \put(332,16){\halfdotline}
\put(270,16){\circle*{4}}
   \put(290,16){\circle{4}}
   \put(310,16){\circle{4}}
   \put(330,16){\circle{4}}
\put(298,18){\small{${}_{1}$}}

\put(260,26){\halfdotline}
  \put(272,26){\dotline}
  \put(292,26){\dotline}
  \put(312,26){\line(1,0){16}}
  \put(332,26){\line(1,0){8}}
\put(270,26){\circle{4}}
   \put(290,26){\circle{4}}
   \put(310,26){\circle{4}}
   \put(330,26){\circle*{4}}
\put(318,28){\small{${}_{1}$}}

\put(310,30){\vector(0,-1){2}}

\end{picture}

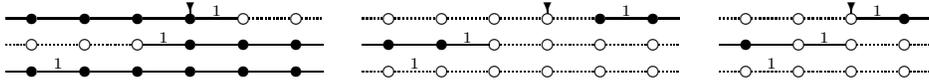
{\small Figure 7. The last graph represents the type $I_b^{-5}$ in
Figure 5. The quotient sheaf $\tilde\LL_\eta$ is represented by
the dotted lines.} \vskip5pt

We now show that the middle one represents the sheaf
$\tilde\EE\pri_\eta$. First, since $\tilde\EE_S$ fits into the
exact sequence
$$0\lra \tilde\EE\pri_S\lra \tilde\EE_S\lra \tilde\LL_\eta\lra 0,
$$
after tensoring with $\cO_{\tilde\cW_\eta}$, we obtain
$$0\lra \tilde\LL_\eta\lra\tilde\EE\pri_\eta\lra \tilde\EE_\eta\lra \tilde\LL_\eta\lra0.
$$
Because $\tilde\EE\pri_\eta$ is locally free, $\tilde\EE\pri_\eta$
must be of the type shown in the middle figure above. Here we used
the fact that the total space $\tilde\cW_S$ is smooth at the
non-locally free loci of $\tilde\LL_\eta$ (as sheaves of
$\cO_{\tilde\cW_\eta}$-modules) and the curves $R_{2+}$ and
$R_{3+}$ are $(-1)$-curves. Consequently, the restriction of
$\tilde\EE\pri_\eta$ to the two rational curves $R_2$ and $R_3$
are of the firm $\cO_{\Po}^{\oplus 3}$.

The study of the sheaves $\EE\pri_\eta$ for $\eta\in\bE_a$
belonging to $\tilde\bD_1-\tilde\bD_3$, $\tilde\bD_3-\tilde\bD_1$
and in the complement of $\tilde\bD_1\cup\tilde\bD_3$ are similar
and will be omitted.
\end{proof}

For completeness, we list their stable modifications as follows.

\begin{picture}(350,50)(0,0)

\def\dotline{\qbezier[10](0,0)(8,0)(15,0)}
\def\halfdotline{\qbezier[5](0,0)(4,0)(7,0)}


\put(0,6){\line(1,0){8}} \put(12,6){\line(1,0){16}}
  \put(32,6){\line(1,0){16}}
  \put(52,6){\line(1,0){16}}
  \put(72,6){\line(1,0){8}}
\put(10,6){\circle*{4}} \put(30,6){\circle*{4}}
   \put(50,6){\circle*{4}}
   \put(70,6){\circle*{4}}
\put(18,8){\small{${}_{1}$}}

\put(0,16){\qbezier[5](0,0)(4,0)(7,0)}
\put(12,16){\dotline}
   \put(32,16){\line(1,0){16}}
   \put(52,16){\line(1,0){16}}
   \put(72,16){\line(1,0){8}}
\put(10,16){\circle{4}} \put(30,16){\circle*{4}}
   \put(50,16){\circle*{4}}
   \put(70,16){\circle*{4}}
\put(38,18){\small{${}_{1}$}}

\put(0,26){\line(1,0){8}}
   \put(12,26){\line(1,0){16}}
   \put(32,26){\line(1,0){16}}
   \put(52,26){\line(1,0){16}}
   \put(72,26){\qbezier[5](0,0)(4,0)(7,0)}
\put(10,26){\circle*{4}}
   \put(30,26){\circle*{4}}
   \put(50,26){\circle{4}}
   \put(70,26){\circle{4}}
\put(58,28){\small{${}_{1}$}} \put(50,30){\vector(0,-1){2}}

 \put(130,6){\halfdotline}
  \put(142,6){\dotline}
  \put(162,6){\dotline}
  \put(182,6){\dotline}
  \put(202,6){\halfdotline}
\put(140,6){\circle{4}}
   \put(160,6){\circle{4}}
   \put(180,6){\circle{4}}
   \put(200,6){\circle{4}}
\put(148,8){\small{${}_{1}$}}

\put(130,16){\line(1,0){8}}
  \put(142,16){\line(1,0){16}}
  \put(162,16){\dotline}
  \put(182,16){\dotline}
  \put(202,16){\halfdotline}
\put(140,16){\circle*{4}}
   \put(160,16){\circle{4}}
   \put(180,16){\circle{4}}
   \put(200,16){\circle{4}}
\put(148,18){\small{${}_{1}$}}

\put(130,26){\halfdotline}
  \put(142,26){\dotline}
  \put(162,26){\dotline}
  \put(182,26){\dotline}
  \put(202,26){\line(1,0){8}}
\put(140,26){\circle{4}}
   \put(160,26){\circle{4}}
   \put(180,26){\circle{4}}
   \put(200,26){\circle*{4}}
\put(180,30){\vector(0,-1){2}}


\put(255,16){{\small $I^{-0}_b$:}}
 \put(280,6){\halfdotline}
  \put(292,6){\dotline}
  \put(312,6){\halfdotline}
   \put(290,6){\circle{4}}
   \put(310,6){\circle{4}}
\put(298,8){\small{${}_{1}$}}

\put(280,16){\line(1,0){8}}
  \put(292,16){\line(1,0){16}}
  \put(312,16){\halfdotline}
   \put(290,16){\circle*{4}}
   \put(310,16){\circle{4}}
\put(298,18){\small{${}_{1}$}}

\put(280,26){\halfdotline}
  \put(292,26){\dotline}
  \put(312,26){\line(1,0){8}}
   \put(290,26){\circle{4}}
   \put(310,26){\circle{4}}
\put(310,30){\vector(0,-1){2}}

\end{picture}

\begin{picture}(350,50)(0,0)

\def\dotline{\qbezier[10](0,0)(8,0)(15,0)}
\def\halfdotline{\qbezier[5](0,0)(4,0)(7,0)}

\put(0,6){\line(1,0){8}}
  \put(12,6){\line(1,0){16}}
  \put(32,6){\line(1,0){16}}
  \put(52,6){\line(1,0){16}}
  \put(72,6){\line(1,0){16}}
  \put(92,6){\line(1,0){8}}
\put(10,6){\circle*{4}}
   \put(30,6){\circle*{4}}
   \put(50,6){\circle*{4}}
   \put(70,6){\circle*{4}}
   \put(90,6){\circle*{4}}
\put(18,8){\small{${}_{1}$}}

\put(0,16){\qbezier[5](0,0)(4,0)(7,0)}
  \put(12,16){\dotline}
   \put(32,16){\dotline}
   \put(52,16){\line(1,0){16}}
   \put(72,16){\line(1,0){16}}
   \put(92,16){\line(1,0){8}}
\put(10,16){\circle{4}}
   \put(30,16){\circle{4}}
   \put(50,16){\circle{4}}
   \put(70,16){\circle*{4}}
   \put(90,16){\circle*{4}}
\put(58,18){\small{${}_{1}$}}

\put(0,26){\line(1,0){8}}
   \put(12,26){\line(1,0){16}}
   \put(32,26){\line(1,0){16}}
   \put(52,26){\line(1,0){16}}
   \put(72,26){\line(1,0){16}}
   \put(92,26){\qbezier[5](0,0)(4,0)(7,0)}
\put(10,26){\circle*{4}}
   \put(30,26){\circle*{4}}
   \put(50,26){\circle*{4}}
   \put(70,26){\circle*{4}}
   \put(90,26){\circle{4}}
\put(78,28){\small{${}_{1}$}} \put(70,30){\vector(0,-1){2}}

 \put(130,6){\qbezier[5](0,0)(4,0)(7,0)}
  \put(142,6){\qbezier[10](0,0)(8,0)(15,0)}
  \put(162,6){\qbezier[10](0,0)(8,0)(15,0)}
  \put(182,6){\qbezier[10](0,0)(8,0)(15,0)}
  \put(202,6){\qbezier[10](0,0)(8,0)(15,0)}
  \put(222,6){\qbezier[5](0,0)(4,0)(7,0)}
\put(140,6){\circle{4}}
   \put(160,6){\circle{4}}
   \put(180,6){\circle{4}}
   \put(200,6){\circle{4}}
   \put(220,6){\circle{4}}
\put(148,8){\small{${}_{1}$}}

\put(130,16){\line(1,0){8}} 
  \put(142,16){\line(1,0){16}}
  \put(162,16){\line(1,0){16}}
  \put(182,16){\qbezier[10](0,0)(8,0)(15,0)}
  \put(202,16){\qbezier[10](0,0)(8,0)(15,0)}
  \put(222,16){\qbezier[5](0,0)(4,0)(7,0)}
\put(140,16){\circle*{4}}
   \put(160,16){\circle*{4}}
   \put(180,16){\circle{4}}
   \put(200,16){\circle{4}}
   \put(220,16){\circle{4}}
\put(168,18){\small{${}_{1}$}}

\put(130,26){\qbezier[5](0,0)(4,0)(7,0)}
  \put(142,26){\dotline}
  \put(162,26){\dotline}
  \put(182,26){\dotline}
  \put(202,26){\dotline}
  \put(222,26){\line(1,0){8}}
\put(140,26){\circle{4}}
   \put(160,26){\circle{4}}
   \put(180,26){\circle{4}}
   \put(200,26){\circle{4}}
   \put(220,26){\circle{4}}
\put(200,30){\vector(0,-1){2}}

\put(256,16){{\small $I^{-1}_b$:}} \put(280,6){\halfdotline}
  \put(292,6){\dotline}
  \put(312,6){\dotline}
  \put(332,6){\halfdotline}
   \put(290,6){\circle{4}}
   \put(310,6){\circle{4}}
   \put(330,6){\circle{4}}
\put(298,8){\small{${}_{1}$}}

\put(280,16){\line(1,0){8}}
  \put(292,16){\line(1,0){16}}
  \put(312,16){\line(1,0){16}}
  \put(332,16){\halfdotline}
   \put(290,16){\circle*{4}}
   \put(310,16){\circle*{4}}
   \put(330,16){\circle{4}}
\put(318,18){\small{${}_{1}$}}

\put(280,26){\qbezier[5](0,0)(4,0)(7,0)}
  \put(292,26){\dotline}
  \put(312,26){\dotline}
  \put(332,26){\line(1,0){8}}
   \put(290,26){\circle{4}}
   \put(310,26){\circle{4}}
   \put(330,26){\circle{4}}
\put(330,30){\vector(0,-1){2}}

\end{picture}

\begin{picture}(350,50)(0,0)

\def\dotline{\qbezier[10](0,0)(8,0)(15,0)}
\def\halfdotline{\qbezier[5](0,0)(4,0)(7,0)}

\put(0,6){\line(1,0){8}}
  \put(12,6){\line(1,0){16}}
  \put(32,6){\line(1,0){16}}
  \put(52,6){\line(1,0){16}}
  \put(72,6){\line(1,0){16}}
  \put(92,6){\line(1,0){8}}
\put(10,6){\circle*{4}}
   \put(30,6){\circle*{4}}
   \put(50,6){\circle*{4}}
   \put(70,6){\circle*{4}}
   \put(90,6){\circle*{4}}
\put(18,8){\small{${}_{1}$}}

\put(0,16){\qbezier[5](0,0)(4,0)(7,0)}
  \put(12,16){\dotline}
   \put(32,16){\line(1,0){16}}
   \put(52,16){\line(1,0){16}}
   \put(72,16){\line(1,0){16}}
   \put(92,16){\line(1,0){8}}
\put(10,16){\circle{4}} \put(30,16){\circle{4}}
   \put(50,16){\circle*{4}}
   \put(70,16){\circle*{4}}
   \put(90,16){\circle*{4}}
\put(38,18){\small{${}_{1}$}}

\put(0,26){\line(1,0){8}}
   \put(12,26){\line(1,0){16}}
   \put(32,26){\line(1,0){16}}
   \put(52,26){\line(1,0){16}}
   \put(72,26){\qbezier[10](0,0)(8,0)(15,0)}
   \put(92,26){\qbezier[5](0,0)(4,0)(7,0)}
\put(10,26){\circle*{4}}
   \put(30,26){\circle*{4}}
   \put(50,26){\circle*{4}}
   \put(70,26){\circle{4}}
   \put(90,26){\circle{4}}
\put(58,28){\small{${}_{1}$}} \put(50,30){\vector(0,-1){2}}

 \put(130,6){\qbezier[5](0,0)(4,0)(7,0)}
  \put(142,6){\qbezier[10](0,0)(8,0)(15,0)}
  \put(162,6){\qbezier[10](0,0)(8,0)(15,0)}
  \put(182,6){\qbezier[10](0,0)(8,0)(15,0)}
  \put(202,6){\qbezier[10](0,0)(8,0)(15,0)}
  \put(222,6){\qbezier[5](0,0)(4,0)(7,0)}
\put(140,6){\circle{4}}
   \put(160,6){\circle{4}}
   \put(180,6){\circle{4}}
   \put(200,6){\circle{4}}
   \put(220,6){\circle{4}}
\put(148,8){\small{${}_{1}$}}

\put(130,16){\line(1,0){8}} 
  \put(142,16){\line(1,0){16}}
  \put(162,16){\qbezier[10](0,0)(8,0)(15,0)}
  \put(182,16){\qbezier[10](0,0)(8,0)(15,0)}
  \put(202,16){\qbezier[10](0,0)(8,0)(15,0)}
  \put(222,16){\qbezier[5](0,0)(4,0)(7,0)}
\put(140,16){\circle*{4}}
   \put(160,16){\circle{4}}
   \put(180,16){\circle{4}}
   \put(200,16){\circle{4}}
   \put(220,16){\circle{4}}
\put(148,18){\small{${}_{1}$}}

\put(130,26){\qbezier[5](0,0)(4,0)(7,0)}
  \put(142,26){\qbezier[10](0,0)(8,0)(15,0)}
  \put(162,26){\qbezier[10](0,0)(8,0)(15,0)}
  \put(182,26){\dotline}
  \put(202,26){\line(1,0){16}}
  \put(222,26){\line(1,0){8}}
\put(140,26){\circle{4}}
   \put(160,26){\circle{4}}
   \put(180,26){\circle{4}}
   \put(200,26){\circle{4}}
   \put(220,26){\circle*{4}}
\put(208,28){\small{${}_{1}$}} \put(180,30){\vector(0,-1){2}}

\put(256,16){{\small $I^{-3}_b$:}} \put(280,6){\halfdotline}
  \put(292,6){\dotline}
  \put(312,6){\dotline}
  \put(332,6){\halfdotline}
   \put(290,6){\circle{4}}
   \put(310,6){\circle{4}}
   \put(330,6){\circle{4}}
\put(298,8){\small{${}_{1}$}}

\put(280,16){\line(1,0){8}}
  \put(292,16){\line(1,0){16}}
  \put(312,16){\dotline}
  \put(332,16){\halfdotline}
   \put(290,16){\circle*{4}}
   \put(310,16){\circle{4}}
   \put(330,16){\circle{4}}
\put(298,18){\small{${}_{1}$}}

\put(280,26){\qbezier[5](0,0)(4,0)(7,0)}
  \put(292,26){\dotline}
  \put(312,26){\line(1,0){16}}
  \put(332,26){\line(1,0){8}}
   \put(290,26){\circle{4}}
   \put(310,26){\circle{4}}
   \put(330,26){\circle*{4}}
\put(318,28){\small{${}_{1}$}} \put(310,30){\vector(0,-1){2}}
\end{picture}

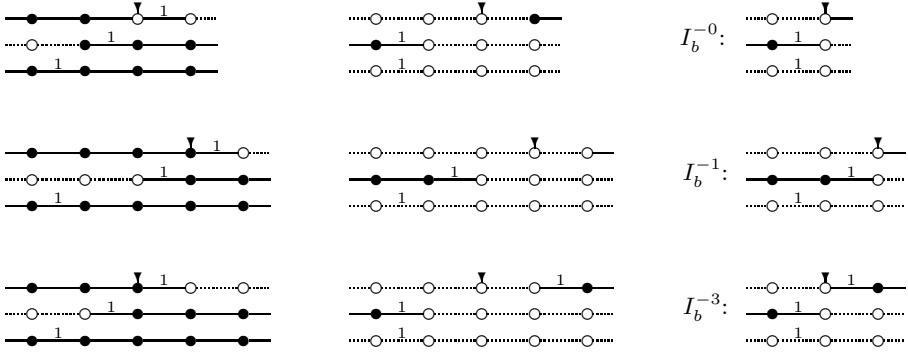
{\small Figure 8. The top, middle and the bottom figures represent
the process of stable modifications of $\EE_\eta$ for $\eta$ in
the complement of $\tilde\bD_1\cup\tilde\bD_3$, in
$\tilde\bD_3-\tilde\bD_1$ and in $\tilde\bD_1-\tilde\bD_3$
respectively.} \vskip5pt

We now construct the stable modification of $\tilde\EE$. We first
contract $\tilde\cW$ along $\tilde\cR_2$ and $\tilde\cR_3$. Since
the restriction of $\tilde\EE\pri$ to fibers of $\tilde\cR_2$ and
$\tilde\cR_3$ are isomorphic to $\cO_{\Po}^{\oplus 3}$, there is a
unique sheaf $\EE\st$ on $\cW\st$ whose pull back to $\tilde\cW$
is $\tilde\EE\pri$. The sheaf $\EE\st$ is called the stable
modification of $\tilde\EE$. The restriction of $\tilde\EE\st$ to
fibers over $\bE_a$ are represented by the right figures in Figure
7.

In the following, for any $\eta\in\bE_a$ we denote by
$\EE_\eta\st$ the restriction of $\EE\st$ to $\cW\st_\eta$. Note
that by applying the same construction to different open subsets
$\bU$, we can construct sheaves $\EE_\eta\st$ for all
$\eta\in\bE_a$, and it is independent of the choices of $\bU$.

Let $\bE_b, \bE_c\sub\tilde\bM_1$ be the proper transforms of
$\II_b^+$ and $\II^+_c$.

\begin{lemm}\lab{4.3}
The sheaves $\EE_\eta\st$ are $0^+$-stable for all $\eta\in\bE_a-
\bE_b\cup\bE_c$.
\end{lemm}

\begin{proof}
Let $\eta\in\bE_a$ be any closed point with $\EE_\eta\st$ the
associated sheaf on $X\lnm\shar$ for some appropriate integers $n$
and $m$. Let $\pi\shar\mh X\lnm\to X\shar\lnm$ be the
desingularization of the marked node and let $\pi\mh X\lnm\to X$
be the contraction of all rational curves. Then
$E_\eta=\pi\lsta\pi^{\dagger\ast}\EE_\eta\st$ is a rank three
locally free sheaves of $\cO_X$-modules with a GPB structure
$E_\eta^0\sub E_\eta|_{p_1+p_2}$ as described in Section 2.
Furthermore, the subsheaf $\LL_\eta\sub \EE_\eta\st$ and the
quotient sheaf $\EE_\eta\st\to\FF_\eta$ defines a sub and a
quotient GPB bundle, $L^G_\eta\sub E_\eta^G$ and $E_\eta^G\to
F^G_\eta$. According to Proposition \ref{2.16}, $\EE_\eta\st$ is
$0$-stable if and only if $E_\eta^G$ is $0$-stable, which is the
case when the extension
\begin{equation}
\lab{4.15} 0\lra L^G\lra E_\eta^G\lra F^G\lra 0
\end{equation}
is non-trivial.
Since $E_\eta^G$ is never $F^G_\eta\oplus L^G_\eta$ when
$\EE_\eta\st$ is of types $I_b^{-0}$ and $I_b^{-1}$, $\EE_\eta\st$
could be $0^+$-unstable only when it was of type $I_b^{-3}$ or of
$I_b^{-5}$. (Here recall that there are no strictly $0$-semistable
vector bundles.)

We now demonstrate that $\eta$ must belong to $\bE_a\cap\bE_b$
when $\EE_\eta\st$ is of type $I_b^{-3}$ and its associated
$E_\eta^G$ is a split extension. Let $S\sub
\tilde\bD_0\cap\tilde\bD_3-\tilde\bD_1$ be a smooth curve
containing $\eta$ and is transversal to $\bE_a$, with $\tilde
\cW_S$ the restriction of $\tilde\cW$ over $S$.
Let $X\times S\sub \tilde\cW_S$ be the main irreducible component
with $\iota\mh X\times\eta\sub X\times S$ the central fiber. Let
$\tilde\EE_S$ and $\tilde\EE\pri_S$ be the associated sheaf on
$\tilde\cW_S$ constructed before (Figure 7). Then that
$\EE_\eta\st$ is $0$-unstable, which is the case when (\ref{4.15})
splits, implies that there must be a surjective sheaf homomorphism
$\tilde\EE\pri_S\to \iota\lsta L(-p_1-p_2)$ so that the composite
\begin{equation}
\lab{4.20} \tilde\LL_\eta\lra
\tilde\EE_S\pri|_{\tilde\cW_\eta}\lra \iota\lsta L(-p_1-p_2)
\end{equation}
is surjective. Let $\tilde\cL_2$ be the sheaf of
$\cO_{\tilde\cW_S}$-modules that fits into the commutative diagram
with the lower sequence exact:
\begin{equation}\lab{4.19}
\begin{CD}
0 @>>> \tilde\EE_S\pri @>>> \tilde\EE_S @>>> \tilde\cL_\eta @>>> 0\\
@. @VVV @VVV @|\\
0 @>>> \iota\lsta L(-p_1-p_2) @>>> \tilde\cL_2 @>>> \tilde\cL_\eta
@>>>0
\end{CD}
\end{equation}
Now let $X_2=X\times\spec \kk[t]/(t^2)$ and let $\iota_2\mh X_2\to
\tilde\cW_S$ be the immersion extending $\iota\mh X\to
\tilde\cW_S$. Since (\ref{4.19}) is exact while the composition of
(\ref{4.20}) is surjective, the pull-back sheaf
$\iota_2\sta\tilde\cL_2$ is an invertible sheaf of
$\cO_{X_2}$-modules and is an extension of $\iota\lsta
L(-p_1-p_2)$ by $\iota\lsta L(-p_1-p_2)$. Now let $S\pri\sub
\bM\uhalf$ be the image of $S\sub\tilde\bM_1$ under the projection
$\tilde\bM_1\to \bM\uhalf$ with $\eta\pri\in S\pri$ the image of
$\eta\in S$. Since $\eta\in\bE_a\cap\bE_b$, $\eta\pri$ belongs to
$\II_a^+\cap\II_b^+$. A direct check shows that the existence of
$\tilde\EE_S\to\tilde\cL_2$ implies that the tangent
$T_{\eta\pri}S\pri$ must be contained in the span of the tangent
spaces
\begin{equation}\lab{4.230}
T_{\eta\pri}\II_a^++T_{\eta\pri}\II_b^+\sub T_{\eta\pri}\bM\uhalf.
\end{equation}
Since $\bE_a$ is the blown-up locus of $\II_a^+$ while $\bE_b$ is
the proper transform of $\II_b^+$, the curve $S$ must specialize
to a point in the intersection $\bE_a\cap\bE_b$, and hence
$\eta\in\bE_a\cap\bE_b$. This proves the claim.

In case $\EE\st_\eta$ is of type $I_b^{-5}$, a similar argument
shows that (\ref{4.230}) also hold. But then since those
$\EE\st_\eta$ of type $I_b^{-5}$ are elementary modification along
direction inside $\tilde\bD_1\cap\tilde\bD_3$, it must be inside
$T_{\eta\pri}\II_a^+$. But this is impossible, which proves that
if $\EE\st_\eta$ is of type $I_b^{-5}$ then it must be $0$-stable.
\end{proof}

\subsection{The family over $\tilde\bM$}

In this subsection, we will construct a family of sheaves over
$\tilde\bM$ after performing an elementary modification to the
family constructed in the previous subsection, and we will then
show that all members of this family are $0$-stable. This way we
can prove the proposition by applying the universal property of
the moduli space $\bM^0$:

\begin{prop}
The birational map $\tilde\bM\,-\!\to \bM^0$ extends to a morphism
$\tilde\bM\to\bM^0$.
\end{prop}

We now prove this proposition. We will sketch the steps that are
parallel to the proof in the previous subsection and provide
details when called for. As mentioned before, the main strategy is
to pull back the (local) tautological families over $\tilde\bM_1$
to $\tilde\bM$, find a small resolution of the base variety of
these families over $\tilde\bM$ and then perform an elementary
modification to these new families. The members of the resulting
families are $0^+$-stable, and hence induce a morphism
$\tilde\bM\to\bM^0$.

We begin with any analytic neighborhood
$\tilde\bU\sub\tilde\bM_1-\bE_c$ and the tautological family
$\EE\st$ on $\tilde\cW$ over $\tilde\bU$ that was constructed in
the previous subsection. When $\xi\in \bE_b$ is away from $\bE_a$,
then $\EE\st_\xi$ must be of the type $I^{+0}_b$ or $I_b^{+1}$, as
shown in figure 4, and accordingly it has an associated quotient
sheaf $\LL_\xi$ and subsheaf $\FF_\xi$ that fits into the exact
sequence
\begin{equation}
\lab{4.31} 0\lra \FF_\xi\lra\EE_\xi\st\lra\LL_\xi\lra 0.
\end{equation}
Note that sheaves of types $I_b^{+2}$ are in $\II_a\cap\II_b$, and
thus won't appear in $\bE_b-\bE_b$. If $\xi\in\bE_b\cap\bE_a$,
then it is of type $\II_b^{-3}$ with the additional property that
the associated GPB short exact sequence of $\EE\st_\xi$ splits, as
proved in Lemma \ref{4.3}. We now show that we can pick a new
associated quotient sheaf of $\EE_\xi\st$ so as to make it of type
$I_b^{+0}$ as well: since $\EE\st_\xi$ is of type $\II_b^{-3}$, it
is a sheaf on $X_{1,1}\shar$ and fits into the exact sequence
$$0\lra \cL_\xi\lra \EE\st_\xi\lra \cF_\xi\lra 0,
$$
according to the proof of Lemma \ref{4.2}. Because the associated
GPB splits, if we let $\LL_\xi$ be the cokernel of
$\cO_{D_1}(-1)\oplus\cO_{D_2}(-1)\to\cL_\xi$ we obtain a unique
surjective
\begin{equation}
\EE\st_\xi\lra \LL_\xi
\end{equation}
so that the composite $\cL_\xi\to\EE\st_\xi\to \LL_\xi$ is exactly
the defining quotient homomorphism $\cL_\xi\to\LL_\xi$. Let
$\FF_\xi$ be the kernel of $\EE\st_\xi\to\LL_\xi$. Then
$\EE\st_\xi$ fits into the exact sequence (\ref{4.31}) as well,
and the latter will be called the associated exact sequence of
$\EE\st_\xi$ and the sheaves $\FF_\xi$ and $\LL_\xi$ be called the
associated sub and quotient sheaves of $\EE\st_\xi$. Clearly,
(\ref{4.31}) makes such $\EE\st_\xi$ a sheaf of type $I_b^{+0}$.

As before, we can make the quotients $\EE\st_\xi\to\LL_\xi$ into a
family of quotients. Let $\tilde\cZ_\xi\sub\tilde\cW_\xi$ be the
support of $\LL_\xi$; $\tilde\cZ_\xi$ is $X\sub X\shar_{1,1}$ when
$\EE\st_\xi$ is of type $I_b^{+0}$, and is $X_{1,0}\sub
X\shar_{2,1}$ when $\EE\st_\xi$ is of type $I_b^{+1}$. The union
$\tilde\cZ=\supp_{\xi\in\bE_b}\tilde\cZ_\xi$ is an irreducible
variety and there is an invertible sheaf $\LL$ of
$\cO_{\tilde\cZ}$-modules and a surjective homomorphism
\begin{equation}\lab{4.332}
\EE\st\lra \iota\lsta\LL,
\end{equation}
where $\iota:\tilde\cZ\hookrightarrow\tilde\cW$ is the inclusion,
so that its restriction to each $\tilde\cW_\xi$ is exactly the
associated homomorphism in (\ref{4.31}).

Our next step is to pull back the family $\EE\st$ to a family over
$\tilde\bM$ and perform elementary modification to it. Let
$\tilde\bV$ be $\tilde\Psi_1\upmo(\tilde\bU)$, which is the
blowing up of $\tilde\bU$ along $\bE_b$; let
$\pi\mh\tilde\bV\to\tilde\bU$ be the projection; let $\tilde\bE_b$
be the exceptional divisor and let $\tilde\bE_a$ be the proper
transform of $\bE_a$. Let
$\cX=\tilde\cW\times_{\tilde\bU}\tilde\bV$ be the pull-back family
over $\tilde\bV$ and let
$\cY=\tilde\cZ\times_{\tilde\bU}\tilde\bV$ the associated
subvariety of $\cX$. As before, the total space of $\cX$ is not
smooth, and we need to small resolve its singularity. For the
moment, we consider the case where all sheaves over $\tilde\bU$
are of types $I_b^{+0}$. Hence, by shrinking $\tilde\cW$ if
necessary we can assume that the singular locus of the fibers of
$\tilde\cW/\tilde\bU$ consists of three smooth connected
codimension two subvarieties: $\cN_0$, $\cN_1$ and $\cN_2$ that
are ordered so that for $\xi\in \bE_b$ the intersection
$\tilde\cW_\xi\cap\cN_i$ is the $i$-th nodal point of
$\tilde\cW_\xi$. Now let
$\tilde\cN_i=\cN_i\times_{\tilde\bU}\tilde\bV$, let
$\bD_i\sub\tilde\bU$ be the image divisor of $\cN_i$ under
$\cX\to\tilde\bU$ and let $\tilde\bD_i\sub\tilde\bV$ be the proper
transform of $\tilde\bD_i$. Then $\tilde\cN_1$ is the marked nodes
of the whole family $\cX$, that $\cX$ is smooth along
$\tilde\cN_0$ except over those $\xi\in\tilde\bD_0\cap\tilde\bE_b$
and smooth along $\tilde\cN_2$ except over those
$\xi\in\tilde\bD_2\cap\tilde\bE_b$. Again, we blow up $\cX$ along
$\tilde\cN_0\times_{\tilde\bV}(\tilde\bD_1\cap\tilde\bE_b)$ and
$\tilde\cN_2\times_{\tilde\bV}(\tilde\bD_3\cap\tilde\bE_b)$; we
then contract one $\bP^1$-factor from each of the two exceptional
divisors to obtain a family of nodal curves $\tilde\cX$. As
before, we choose the contraction so that if we let
$\tilde\cY\sub\tilde\cX$ be the proper transform of $\cY\sub \cX$,
and let $\tilde\cD_0$ and $\tilde\cD_2$ be the two exceptional
loci of $\tilde\cX\to\cX$, then the intersection
$\tilde\cY\cap\tilde\cD_0$ and $\tilde\cY\cap\tilde\cD_2$ are
finite over $\tilde\bV$; Namely it contains no curves that lie
inside a single fiber $\tilde\cX_\xi$ over some $\xi\in\tilde\bV$.

We now pull back the family $\EE\st$ to $\tilde\cX$ and perform an
elementary modification. Let $p\mh \tilde\cX\to \tilde\cW$ be the
projection, and let
\begin{equation}\lab{4.335}
p\sta\EE\st\to p\sta\iota\lsta\LL
\end{equation}
be the pull-back of the pair (\ref{4.332}). Then the kernel
$\tilde\EE$ of the homomorphism (\ref{4.335}) is the modification
we seek for.

It remains to prove that all members in $\tilde\EE$ are
$0^+$-stable. Let $\xi\in\bE_b$ be any point and let
$\eta\in\tilde\bE_b$ be any of its lifts. Since $\xi\in
\bD_0\cap\bD_2$, $\eta$ can possibly be in
$\tilde\bD_0\cap\tilde\bD_2$, in $\tilde\bD_0-\tilde\bD_2$, in
$\tilde\bD_2-\tilde\bD_0$ or is away from
$\tilde\bD_0\cup\tilde\bD_2$. Now we investigate in detail the
sheaf $\tilde\EE_\eta$ when $\eta\in\tilde\bD_0\cap\tilde\bD_2$.
First, the curve $\tilde\cX_\eta$ is $X_{2,2}\shar$ with
$\tilde\cX_\eta\to\cX_\eta$ the contraction of the first and the
last rational curves; When $\eta_t$ is curve in $\tilde\bV$ with
$\eta_0=\eta$ and is normal to $\tilde\bE_b$, the family
$\tilde\cX_{\eta_t}$ smooths the first\footnote{Recall that our
convention is to label the nodes from the $0$-th to the $3$-rd if
there are four nodes.} and the third node of $\tilde\cX_\eta$. The
pull back
$$(p\sta\tilde\EE)_\eta\lra (p\sta\iota\lsta\LL)_\eta
$$
is represented by the left graph below with
$(p\sta\iota\lsta\LL)_\eta$ represented by the dotted lines. The
modified sheaf $\tilde\EE_\eta$ is represented by the middle
graph. After contracting the third rational curve in
$\tilde\cX_\eta$ we obtain a sheaf shown in the right graph that
is of type $I_b^{-4}$:

\begin{picture}(120,40)(0,0)

\put(0,6){\line(1,0){8}}
  \put(12,6){\line(1,0){16}}
  \put(32,6){\line(1,0){16}}
  \put(52,6){\line(1,0){16}}
  \put(72,6){\line(1,0){16}}
  \put(92,6){\line(1,0){8}}
\put(10,6){\circle*{4}}
   \put(30,6){\circle*{4}}
   \put(50,6){\circle*{4}}
   \put(70,6){\circle*{4}}
   \put(90,6){\circle*{4}}
\put(38,8){\small{${}_{1}$}}

\put(0,16){\qbezier[5](0,0)(4,0)(7,0)}
  \put(12,16){\dotline}
   \put(32,16)\solid
   \put(52,16){\line(1,0){16}}
   \put(72,16){\line(1,0){16}}
   \put(92,16){\line(1,0){8}}
\put(10,16){\circle{4}}
   \put(30,16){\circle{4}}
   \put(50,16){\circle*{4}}
   \put(70,16){\circle*{4}}
   \put(90,16){\circle*{4}}
\put(38,18){\small{${}_{1}$}}

\put(0,26){\line(1,0){8}}
   \put(12,26){\line(1,0){16}}
   \put(32,26){\line(1,0){16}}
   \put(52,26){\line(1,0){16}}
   \put(72,26)\dotline
   \put(92,26)\halfdotline
\put(10,26){\circle*{4}}
   \put(30,26){\circle*{4}}
   \put(50,26){\circle*{4}}
   \put(70,26){\circle{4}}
   \put(90,26){\circle{4}}
\put(58,28){\small{${}_{1}$}} \put(50,30){\vector(0,-1){2}}

\end{picture}
\begin{picture}(120,40)(0,0)

\put(50,30){\vector(0,-1){2}}

\put(0,26)\halfdotline
   \put(12,26)\dotline
   \put(32,26)\dotline
   \put(52,26)\dotline
   \put(72,26)\solid
   \put(92,26)\halfsolid
\put(10,26){\circle{4}}
   \put(30,26){\circle{4}}
   \put(50,26){\circle{4}}
   \put(70,26){\circle{4}}
   \put(90,26){\circle*{4}}
\put(78,28){\small{${}_{1}$}}

\put(0,16)\halfsolid
  \put(12,16)\solid
   \put(32,16)\dotline
   \put(52,16)\dotline
   \put(72,16)\dotline
   \put(92,16)\halfdotline
\put(10,16){\circle*{4}}
   \put(30,16){\circle{4}}
   \put(50,16){\circle{4}}
   \put(70,16){\circle{4}}
   \put(90,16){\circle{4}}
\put(18,18){\small{${}_{1}$}}

\put(0,6)\halfdotline
  \put(12,6)\dotline
  \put(32,6)\dotline
  \put(52,6)\dotline
  \put(72,6)\dotline
  \put(92,6)\halfdotline
\put(10,6){\circle{4}}
   \put(30,6){\circle{4}}
   \put(50,6){\circle{4}}
   \put(70,6){\circle{4}}
   \put(90,6){\circle{4}}
\put(38,8){\small{${}_{1}$}}

\end{picture}
\begin{picture}(120,40)(0,0)

\put(0,26)\halfdotline
   \put(12,26)\dotline
   \put(32,26)\dotline
   \put(52,26)\solid
   \put(72,26)\halfsolid
\put(10,26){\circle{4}}
   \put(30,26){\circle{4}}
   \put(50,26){\circle{4}}
   \put(70,26){\circle*{4}}
\put(58,28){\small{${}_{1}$}} \put(50,30){\vector(0,-1){2}}

\put(0,16)\halfsolid
  \put(12,16)\solid
   \put(32,16)\dotline
   \put(52,16)\dotline
   \put(72,16)\halfdotline
\put(10,16){\circle*{4}}
   \put(30,16){\circle{4}}
   \put(50,16){\circle{4}}
   \put(70,16){\circle{4}}
\put(18,18){\small{${}_{1}$}}

\put(0,6)\halfdotline
  \put(12,6)\dotline
  \put(32,6)\dotline
  \put(52,6)\dotline
  \put(72,6)\halfdotline
\put(10,6){\circle{4}}
   \put(30,6){\circle{4}}
   \put(50,6){\circle{4}}
   \put(70,6){\circle{4}}
\put(38,8){\small{${}_{1}$}}

\end{picture}

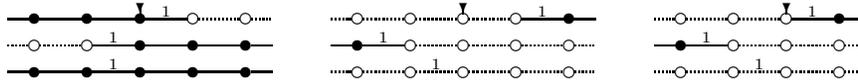
{\small Figure 9. The graphs represent the sheaf over
$\eta\in\tilde\bD_0\cap\tilde\bD_2$ before the modification, after
the modification and after the stabilization. The resulting sheaf
is of type $I_b^{-4}$.}

\vsp

\noindent We will call the sheaf obtained after contracting the
third rational curve the stable modification of
$(p\sta\EE\st)_\eta$, and will denote it by $\tilde\EE\st_\eta$.

We can derive the other stable modifications $\tilde\EE\st_\eta$
similarly and will give the graphs sketching their respective
process as follows:

\begin{picture}(120,30)(0,0)

\put(0,26){\line(1,0){8}}
   \put(12,26){\line(1,0){16}}
   \put(32,26){\line(1,0){16}}
   \put(52,26){\line(1,0){16}}
   \put(72,26)\halfdotline
\put(10,26){\circle*{4}}
   \put(30,26){\circle*{4}}
   \put(50,26){\circle*{4}}
   \put(70,26){\circle{4}}
\put(58,28){\small{${}_{1}$}} \put(50,30){\vector(0,-1){2}}

\put(0,16){\qbezier[5](0,0)(4,0)(7,0)}
  \put(12,16){\dotline}
   \put(32,16)\solid
   \put(52,16){\line(1,0){16}}
   \put(72,16){\line(1,0){8}}
\put(10,16){\circle{4}}
   \put(30,16){\circle{4}}
   \put(50,16){\circle*{4}}
   \put(70,16){\circle*{4}}
\put(38,18){\small{${}_{1}$}}

 \put(0,6){\line(1,0){8}}
  \put(12,6){\line(1,0){16}}
  \put(32,6){\line(1,0){16}}
  \put(52,6){\line(1,0){16}}
  \put(72,6){\line(1,0){8}}
\put(10,6){\circle*{4}}
   \put(30,6){\circle*{4}}
   \put(50,6){\circle*{4}}
   \put(70,6){\circle*{4}}
\put(38,8){\small{${}_{1}$}}

\end{picture}
\begin{picture}(120,30)(0,0)

\put(0,26)\halfdotline
   \put(12,26)\dotline
   \put(32,26)\dotline
   \put(52,26)\dotline
   \put(72,26)\halfsolid
\put(10,26){\circle{4}}
   \put(30,26){\circle{4}}
   \put(50,26){\circle{4}}
   \put(70,26){\circle{4}}
\put(50,30){\vector(0,-1){2}}

\put(0,16)\halfsolid
  \put(12,16)\solid
   \put(32,16)\dotline
   \put(52,16)\dotline
   \put(72,16)\halfdotline
\put(10,16){\circle*{4}}
   \put(30,16){\circle{4}}
   \put(50,16){\circle{4}}
   \put(70,16){\circle{4}}
\put(18,18){\small{${}_{1}$}}

\put(0,6)\halfdotline
  \put(12,6)\dotline
  \put(32,6)\dotline
  \put(52,6)\dotline
  \put(72,6)\halfdotline
\put(10,6){\circle{4}}
   \put(30,6){\circle{4}}
   \put(50,6){\circle{4}}
   \put(70,6){\circle{4}}
\put(38,8){\small{${}_{1}$}}

\end{picture}
\begin{picture}(120,30)(0,0)

\put(0,26)\halfdotline
   \put(12,26)\dotline
   \put(32,26)\dotline
   \put(52,26)\halfsolid
\put(10,26){\circle{4}}
   \put(30,26){\circle{4}}
   \put(50,26){\circle{4}}
\put(50,30){\vector(0,-1){2}}

\put(0,16)\halfsolid
  \put(12,16)\solid
   \put(32,16)\dotline
   \put(52,16)\halfdotline
\put(10,16){\circle*{4}}
   \put(30,16){\circle{4}}
   \put(50,16){\circle{4}}
\put(18,18){\small{${}_{1}$}}

\put(0,6)\halfdotline
  \put(12,6)\dotline
  \put(32,6)\dotline
  \put(52,6)\halfdotline
\put(10,6){\circle{4}}
   \put(30,6){\circle{4}}
   \put(50,6){\circle{4}}
\put(38,8){\small{${}_{1}$}}

\end{picture}
\begin{picture}(120,40)(0,0)

\put(0,26){\line(1,0){8}}
   \put(12,26){\line(1,0){16}}
   \put(32,26){\line(1,0){16}}
   \put(52,26)\dotline
   \put(72,26)\halfdotline
\put(10,26){\circle*{4}}
   \put(30,26){\circle*{4}}
   \put(50,26){\circle{4}}
   \put(70,26){\circle{4}}
\put(38,28){\small{${}_{1}$}} \put(30,30){\vector(0,-1){2}}

\put(0,16){\qbezier[5](0,0)(4,0)(7,0)}
  \put(12,16)\solid
   \put(32,16)\solid
   \put(52,16){\line(1,0){16}}
   \put(72,16){\line(1,0){8}}
\put(10,16){\circle{4}}
   \put(30,16){\circle*{4}}
   \put(50,16){\circle*{4}}
   \put(70,16){\circle*{4}}
\put(18,18){\small{${}_{1}$}}

 \put(0,6){\line(1,0){8}}
  \put(12,6){\line(1,0){16}}
  \put(32,6){\line(1,0){16}}
  \put(52,6){\line(1,0){16}}
  \put(72,6){\line(1,0){8}}
\put(10,6){\circle*{4}}
   \put(30,6){\circle*{4}}
   \put(50,6){\circle*{4}}
   \put(70,6){\circle*{4}}
\put(18,8){\small{${}_{1}$}}

\end{picture}
\begin{picture}(120,40)(0,0)

\put(30,30){\vector(0,-1){2}}

\put(0,26)\halfdotline
   \put(12,26)\dotline
   \put(32,26)\dotline
   \put(52,26)\solid
   \put(72,26)\halfsolid
\put(10,26){\circle{4}}
   \put(30,26){\circle{4}}
   \put(50,26){\circle{4}}
   \put(70,26){\circle*{4}}
\put(58,28){\small{${}_{1}$}}

\put(0,16)\halfsolid
  \put(12,16)\dotline
   \put(32,16)\dotline
   \put(52,16)\dotline
   \put(72,16)\halfdotline
\put(10,16){\circle{4}}
   \put(30,16){\circle{4}}
   \put(50,16){\circle{4}}
   \put(70,16){\circle{4}}

\put(0,6)\halfdotline
  \put(12,6)\dotline
  \put(32,6)\dotline
  \put(52,6)\dotline
  \put(72,6)\halfdotline
\put(10,6){\circle{4}}
   \put(30,6){\circle{4}}
   \put(50,6){\circle{4}}
   \put(70,6){\circle{4}}
\put(18,8){\small{${}_{1}$}}

\end{picture}
\begin{picture}(120,40)(0,0)

\put(30,30){\vector(0,-1){2}}

\put(0,26)\halfdotline
   \put(12,26)\dotline
   \put(32,26)\solid
   \put(52,26)\halfsolid
\put(10,26){\circle{4}}
   \put(30,26){\circle{4}}
   \put(50,26){\circle*{4}}
\put(38,28){\small{${}_{1}$}}

\put(0,16)\halfsolid
  \put(12,16)\dotline
   \put(32,16)\dotline
   \put(52,16)\halfdotline
\put(10,16){\circle{4}}
   \put(30,16){\circle{4}}
   \put(50,16){\circle{4}}

\put(0,6)\halfdotline
  \put(12,6)\dotline
  \put(32,6)\dotline
  \put(52,6)\halfdotline
\put(10,6){\circle{4}}
   \put(30,6){\circle{4}}
   \put(50,6){\circle{4}}
\put(18,8){\small{${}_{1}$}}

\end{picture}
\begin{picture}(120,40)(-10,0)

\put(30,30){\vector(0,-1){2}}

\put(0,26){\line(1,0){8}}
   \put(12,26){\line(1,0){16}}
   \put(32,26){\line(1,0){16}}
   \put(52,26)\halfdotline
\put(10,26){\circle*{4}}
   \put(30,26){\circle*{4}}
   \put(50,26){\circle{4}}
\put(38,28){\small{${}_{1}$}}

\put(0,16){\qbezier[5](0,0)(4,0)(7,0)}
  \put(12,16)\solid
   \put(32,16)\solid
   \put(52,16){\line(1,0){8}}
\put(10,16){\circle{4}}
   \put(30,16){\circle*{4}}
   \put(50,16){\circle*{4}}
\put(18,18){\small{${}_{1}$}}

 \put(0,6){\line(1,0){8}}
  \put(12,6){\line(1,0){16}}
  \put(32,6){\line(1,0){16}}
  \put(52,6){\line(1,0){8}}
\put(10,6){\circle*{4}}
   \put(30,6){\circle*{4}}
   \put(50,6){\circle*{4}}
\put(18,8){\small{${}_{1}$}}

\end{picture}
\begin{picture}(120,40)(-10,0)

\put(30,30){\vector(0,-1){2}}

\put(0,26)\halfdotline
   \put(12,26)\dotline
   \put(32,26)\dotline
   \put(52,26)\halfsolid
\put(10,26){\circle{4}}
   \put(30,26){\circle{4}}
   \put(50,26){\circle{4}}

\put(0,16)\halfsolid
  \put(12,16)\dotline
   \put(32,16)\dotline
   \put(52,16)\halfdotline
\put(10,16){\circle{4}}
   \put(30,16){\circle{4}}
   \put(50,16){\circle{4}}

\put(0,6)\halfdotline
  \put(12,6)\dotline
  \put(32,6)\dotline
  \put(52,6)\halfdotline
\put(10,6){\circle{4}}
   \put(30,6){\circle{4}}
   \put(50,6){\circle{4}}
\put(18,8){\small{${}_{1}$}}

\end{picture}
\begin{picture}(120,40)(-10,0)

\put(30,30){\vector(0,-1){2}}

\put(0,26)\halfdotline
   \put(12,26)\dotline
   \put(32,26)\halfsolid
\put(10,26){\circle{4}}
   \put(30,26){\circle{4}}

\put(0,16)\halfsolid
  \put(12,16)\dotline
   \put(32,16)\halfdotline
\put(10,16){\circle{4}}
   \put(30,16){\circle{4}}

\put(0,6)\halfdotline
  \put(12,6)\dotline
  \put(32,6)\halfdotline
\put(10,6){\circle{4}}
   \put(30,6){\circle{4}}
\put(18,8){\small{${}_{1}$}}

\end{picture}

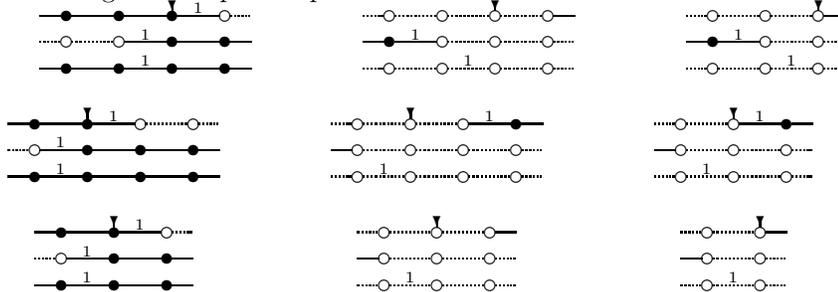
{\small Figure 10. From top to bottom, they represent the
derivation of $\tilde\EE\st_\eta$ in case $\eta$ is in
$\eta\in\tilde\bD_2-\tilde\bD_0$, in $\tilde\bD_0-\tilde\bD_2$ and
is away from $\tilde\bD_0\cup\tilde\bD_2$.}

\vsp

Exactly as in the case before, we can contract all those rational
curves that are immediately to the right of the marked nodes in
$\tilde\cX_\eta$ for all $\eta\in\tilde\bE_b$ simultaneously and
obtain a new family $\tilde\cX\st$ that has smooth total space.
Let $\varphi\mh \tilde\cX\to\tilde\cX\st$ be the stabilization,
then $\tilde\EE\st=\varphi\lsta\tilde\EE$ is a locally free and
the sheaves $\tilde\EE\st_\eta$ for $\eta\in\bE_b$ are exactly the
sheaves shown in the right column in figures 9 and 10.

\begin{lemm}
All stable modifications $\tilde\EE\st_\eta$ derived so far are
$0$-stable.
\end{lemm}

\begin{proof}
We will only sketch the proof here, since the details are exactly
the same as in the proof of Lemma \ref{4.3}. First,
$\tilde\EE_\eta$ could be $0$-unstable only if its associated GPB
were split. This is possible only when it is of type $I_b^{-4}$.
Because $\eta\in\tilde\bD_0\cap\tilde\bD_2$, the split of its
associated GPB implies that the tangent in $T_\xi\tilde\bM_1$
associated to $\eta$ lies in $T_\xi\bE_b$, which is impossible.
Thus the associated GPB is irreducible and hence all stable
modification derived are $0$-stable.
\end{proof}


We now consider the case where $\EE\st_\xi$ is of type $I_b^{+1}$.
Let $\eta\in\tilde\bE_b$ be a lift of $\xi$, let $\cX/\tilde\bV$
be an analytic neighborhood of $\eta$ as before and let
$\tilde\cX\to\cX$ be the small resolution constructed according to
a similar rule. We will have similar quotient family
$p\sta\EE\st\to p\sta\iota\lsta\LL$, and we will take the
$\tilde\EE$ be the kernel of this homomorphism. We still need to
determine the types of members in $\tilde\EE$. As before, we let
$\cN_0,\cdots,\cN_3$ be the loci of singular points of the fibers
of $\cW/\tilde\bU$, let $\cD_0,\cdots,\cD_3$ be their respective
images in $\tilde\bU$ and let $\tilde\cD_0,\cdots,\tilde\cD_3$ be
their proper transforms in $\tilde\bV$. The resulting type of the
stable modification will depend on whether $\eta$ is in
$\tilde\bD_0\cap\tilde\bD_2$, in $\tilde\bD_2-\tilde\bD_0$, in
$\tilde\bD_0-\tilde\bD_2$ or is away from
$\tilde\bD_0\cup\tilde\bD_2$. We will show their respective stable
modification by providing their associated graphes as before:

\begin{picture}(127,40)(0,0)

\put(70,30){\vector(0,-1){2}}

\put(0,26){\line(1,0){8}}
   \put(12,26){\line(1,0){16}}
   \put(32,26){\line(1,0){16}}
   \put(52,26){\line(1,0){16}}
   \put(72,26)\solid
   \put(92,26)\dotline
   \put(112,26)\halfdotline
\put(10,26){\circle*{4}}
   \put(30,26){\circle*{4}}
   \put(50,26){\circle*{4}}
   \put(70,26){\circle*{4}}
   \put(90,26){\circle{4}}
   \put(110,26){\circle{4}}
\put(78,28){\small{${}_{1}$}}

\put(0,16){\qbezier[5](0,0)(4,0)(7,0)}
  \put(12,16){\dotline}
   \put(32,16)\solid
   \put(52,16){\line(1,0){16}}
   \put(72,16){\line(1,0){16}}
   \put(92,16)\solid
   \put(112,16)\halfsolid
\put(10,16){\circle{4}}
   \put(30,16){\circle{4}}
   \put(50,16){\circle*{4}}
   \put(70,16){\circle*{4}}
   \put(90,16){\circle*{4}}
   \put(110,16){\circle*{4}}
\put(38,18){\small{${}_{1}$}}

\put(0,6){\line(1,0){8}}
  \put(12,6){\line(1,0){16}}
  \put(32,6){\line(1,0){16}}
  \put(52,6){\line(1,0){16}}
  \put(72,6){\line(1,0){16}}
  \put(92,6)\solid
  \put(112,6)\halfsolid
\put(10,6){\circle*{4}}
   \put(30,6){\circle*{4}}
   \put(50,6){\circle*{4}}
   \put(70,6){\circle*{4}}
   \put(90,6){\circle*{4}}
   \put(110,6){\circle*{4}}
\put(58,8){\small{${}_{1}$}}

\end{picture}
\begin{picture}(127,40)(0,0)

\put(70,30){\vector(0,-1){2}}

\put(0,26)\halfdotline
   \put(12,26)\dotline
   \put(32,26)\dotline
   \put(52,26)\dotline
   \put(72,26)\dotline
   \put(92,26)\solid
   \put(112,26)\halfsolid
\put(10,26){\circle{4}}
   \put(30,26){\circle{4}}
   \put(50,26){\circle{4}}
   \put(70,26){\circle{4}}
   \put(90,26){\circle{4}}
   \put(110,26){\circle*{4}}
\put(78,28){\small{${}_{1}$}}

\put(0,16)\halfsolid
  \put(12,16)\solid
   \put(32,16)\dotline
   \put(52,16)\dotline
   \put(72,16)\dotline
   \put(92,16)\dotline
   \put(112,16)\halfdotline
\put(10,16){\circle*{4}}
   \put(30,16){\circle{4}}
   \put(50,16){\circle{4}}
   \put(70,16){\circle{4}}
   \put(90,16){\circle{4}}
   \put(110,16){\circle{4}}
\put(18,18){\small{${}_{1}$}}

\put(0,6)\halfdotline
  \put(12,6)\dotline
  \put(32,6)\dotline
  \put(52,6)\dotline
  \put(72,6)\dotline
  \put(92,6)\dotline
   \put(112,6)\halfdotline
\put(10,6){\circle{4}}
   \put(30,6){\circle{4}}
   \put(50,6){\circle{4}}
   \put(70,6){\circle{4}}
   \put(90,6){\circle{4}}
   \put(110,6){\circle*{4}}
\put(58,8){\small{${}_{1}$}}

\end{picture}
\begin{picture}(120,40)(0,0)

\put(50,30){\vector(0,-1){2}}

\put(0,26)\halfdotline
   \put(12,26)\dotline
   \put(32,26)\dotline
   \put(52,26)\solid
   \put(72,26)\halfsolid
\put(10,26){\circle{4}}
   \put(30,26){\circle{4}}
   \put(50,26){\circle{4}}
   \put(70,26){\circle*{4}}
\put(58,28){\small{${}_{1}$}}

\put(0,16)\halfsolid
  \put(12,16)\solid
   \put(32,16)\dotline
   \put(52,16)\dotline
   \put(72,16)\halfdotline
\put(10,16){\circle*{4}}
   \put(30,16){\circle{4}}
   \put(50,16){\circle{4}}
   \put(70,16){\circle{4}}
\put(18,18){\small{${}_{1}$}}

\put(0,6)\halfdotline
  \put(12,6)\dotline
  \put(32,6)\dotline
  \put(52,6)\dotline
  \put(72,6)\halfdotline
\put(10,6){\circle{4}}
   \put(30,6){\circle{4}}
   \put(50,6){\circle{4}}
   \put(70,6){\circle{4}}
\put(38,8){\small{${}_{1}$}}

\end{picture}
\begin{picture}(130,40)(0,0)

\put(70,30){\vector(0,-1){2}}

\put(0,26){\line(1,0){8}}
   \put(12,26){\line(1,0){16}}
   \put(32,26){\line(1,0){16}}
   \put(52,26){\line(1,0){16}}
   \put(72,26)\solid
   \put(92,26)\halfdotline
\put(10,26){\circle*{4}}
   \put(30,26){\circle*{4}}
   \put(50,26){\circle*{4}}
   \put(70,26){\circle*{4}}
   \put(90,26){\circle{4}}
\put(78,28){\small{${}_{1}$}}

\put(0,16){\qbezier[5](0,0)(4,0)(7,0)}
  \put(12,16){\dotline}
   \put(32,16)\solid
   \put(52,16){\line(1,0){16}}
   \put(72,16){\line(1,0){16}}
   \put(92,16)\halfsolid
\put(10,16){\circle{4}}
   \put(30,16){\circle{4}}
   \put(50,16){\circle*{4}}
   \put(70,16){\circle*{4}}
   \put(90,16){\circle*{4}}
\put(38,18){\small{${}_{1}$}}

\put(0,6){\line(1,0){8}}
  \put(12,6){\line(1,0){16}}
  \put(32,6){\line(1,0){16}}
  \put(52,6){\line(1,0){16}}
  \put(72,6){\line(1,0){16}}
  \put(92,6)\halfsolid
\put(10,6){\circle*{4}}
   \put(30,6){\circle*{4}}
   \put(50,6){\circle*{4}}
   \put(70,6){\circle*{4}}
   \put(90,6){\circle*{4}}
\put(58,8){\small{${}_{1}$}}

\end{picture}
\begin{picture}(127,40)(0,0)

\put(70,30){\vector(0,-1){2}}

\put(0,26)\halfdotline
   \put(12,26)\dotline
   \put(32,26)\dotline
   \put(52,26)\dotline
   \put(72,26)\dotline
   \put(92,26)\halfsolid
\put(10,26){\circle{4}}
   \put(30,26){\circle{4}}
   \put(50,26){\circle{4}}
   \put(70,26){\circle{4}}
   \put(90,26){\circle{4}}

\put(0,16)\halfsolid
  \put(12,16)\solid
   \put(32,16)\dotline
   \put(52,16)\dotline
   \put(72,16)\dotline
   \put(92,16)\halfdotline
\put(10,16){\circle*{4}}
   \put(30,16){\circle{4}}
   \put(50,16){\circle{4}}
   \put(70,16){\circle{4}}
   \put(90,16){\circle{4}}
\put(18,18){\small{${}_{1}$}}

\put(0,6)\halfdotline
  \put(12,6)\dotline
  \put(32,6)\dotline
  \put(52,6)\dotline
  \put(72,6)\dotline
  \put(92,6)\halfdotline
\put(10,6){\circle{4}}
   \put(30,6){\circle{4}}
   \put(50,6){\circle{4}}
   \put(70,6){\circle{4}}
   \put(90,6){\circle{4}}
\put(58,8){\small{${}_{1}$}}

\end{picture}
\begin{picture}(120,30)(0,0)

\put(0,26)\halfdotline
   \put(12,26)\dotline
   \put(32,26)\dotline
   \put(52,26)\halfsolid
\put(10,26){\circle{4}}
   \put(30,26){\circle{4}}
   \put(50,26){\circle{4}}
\put(50,30){\vector(0,-1){2}}

\put(0,16)\halfsolid
  \put(12,16)\solid
   \put(32,16)\dotline
   \put(52,16)\halfdotline
\put(10,16){\circle*{4}}
   \put(30,16){\circle{4}}
   \put(50,16){\circle{4}}
\put(18,18){\small{${}_{1}$}}

\put(0,6)\halfdotline
  \put(12,6)\dotline
  \put(32,6)\dotline
  \put(52,6)\halfdotline
\put(10,6){\circle{4}}
   \put(30,6){\circle{4}}
   \put(50,6){\circle{4}}
\put(38,8){\small{${}_{1}$}}

\end{picture}
\begin{picture}(130,40)(0,0)

\put(50,30){\vector(0,-1){2}}

\put(0,26){\line(1,0){8}}
   \put(12,26){\line(1,0){16}}
   \put(32,26){\line(1,0){16}}
   \put(52,26){\line(1,0){16}}
   \put(72,26)\dotline
   \put(92,26)\halfdotline
\put(10,26){\circle*{4}}
   \put(30,26){\circle*{4}}
   \put(50,26){\circle*{4}}
   \put(70,26){\circle{4}}
   \put(90,26){\circle{4}}
\put(58,28){\small{${}_{1}$}}

\put(0,16){\qbezier[5](0,0)(4,0)(7,0)}
  \put(12,16)\solid
   \put(32,16)\solid
   \put(52,16){\line(1,0){16}}
   \put(72,16){\line(1,0){16}}
   \put(92,16)\halfsolid
\put(10,16){\circle{4}}
   \put(30,16){\circle*{4}}
   \put(50,16){\circle*{4}}
   \put(70,16){\circle*{4}}
   \put(90,16){\circle*{4}}
\put(18,18){\small{${}_{1}$}}

\put(0,6){\line(1,0){8}}
  \put(12,6){\line(1,0){16}}
  \put(32,6){\line(1,0){16}}
  \put(52,6){\line(1,0){16}}
  \put(72,6){\line(1,0){16}}
  \put(92,6)\halfsolid
\put(10,6){\circle*{4}}
   \put(30,6){\circle*{4}}
   \put(50,6){\circle*{4}}
   \put(70,6){\circle*{4}}
   \put(90,6){\circle*{4}}
\put(38,8){\small{${}_{1}$}}

\end{picture}
\begin{picture}(127,40)(0,0)

\put(50,30){\vector(0,-1){2}}

\put(0,26)\halfdotline
   \put(12,26)\dotline
   \put(32,26)\dotline
   \put(52,26)\dotline
   \put(72,26)\solid
   \put(92,26)\halfsolid
\put(10,26){\circle{4}}
   \put(30,26){\circle{4}}
   \put(50,26){\circle{4}}
   \put(70,26){\circle{4}}
   \put(90,26){\circle*{4}}
\put(78,28){\small{${}_{1}$}}

\put(0,16)\halfsolid
  \put(12,16)\dotline
   \put(32,16)\dotline
   \put(52,16)\dotline
   \put(72,16)\dotline
   \put(92,16)\halfdotline
\put(10,16){\circle{4}}
   \put(30,16){\circle{4}}
   \put(50,16){\circle{4}}
   \put(70,16){\circle{4}}
   \put(90,16){\circle{4}}

\put(0,6)\halfdotline
  \put(12,6)\dotline
  \put(32,6)\dotline
  \put(52,6)\dotline
  \put(72,6)\dotline
  \put(92,6)\halfdotline
\put(10,6){\circle{4}}
   \put(30,6){\circle{4}}
   \put(50,6){\circle{4}}
   \put(70,6){\circle{4}}
   \put(90,6){\circle{4}}
\put(38,8){\small{${}_{1}$}}

\end{picture}
\begin{picture}(120,40)(0,0)

\put(30,30){\vector(0,-1){2}}

\put(0,26)\halfdotline
   \put(12,26)\dotline
   \put(32,26)\solid
   \put(52,26)\halfsolid
\put(10,26){\circle{4}}
   \put(30,26){\circle{4}}
   \put(50,26){\circle*{4}}
\put(38,28){\small{${}_{1}$}}

\put(0,16)\halfsolid
  \put(12,16)\dotline
   \put(32,16)\dotline
   \put(52,16)\halfdotline
\put(10,16){\circle{4}}
   \put(30,16){\circle{4}}
   \put(50,16){\circle{4}}

\put(0,6)\halfdotline
  \put(12,6)\dotline
  \put(32,6)\dotline
  \put(52,6)\halfdotline
\put(10,6){\circle{4}}
   \put(30,6){\circle{4}}
   \put(50,6){\circle{4}}
\put(18,8){\small{${}_{1}$}}

\end{picture}
\begin{picture}(130,40)(-10,0)

\put(50,30){\vector(0,-1){2}}

\put(0,26){\line(1,0){8}}
   \put(12,26){\line(1,0){16}}
   \put(32,26){\line(1,0){16}}
   \put(52,26){\line(1,0){16}}
   \put(72,26)\halfdotline
\put(10,26){\circle*{4}}
   \put(30,26){\circle*{4}}
   \put(50,26){\circle*{4}}
   \put(70,26){\circle{4}}
\put(58,28){\small{${}_{1}$}}

\put(0,16){\qbezier[5](0,0)(4,0)(7,0)}
  \put(12,16)\solid
   \put(32,16)\solid
   \put(52,16){\line(1,0){16}}
   \put(72,16)\halfsolid
\put(10,16){\circle{4}}
   \put(30,16){\circle*{4}}
   \put(50,16){\circle*{4}}
   \put(70,16){\circle*{4}}
\put(18,18){\small{${}_{1}$}}

\put(0,6){\line(1,0){8}}
  \put(12,6){\line(1,0){16}}
  \put(32,6){\line(1,0){16}}
  \put(52,6){\line(1,0){16}}
  \put(72,6)\halfsolid
\put(10,6){\circle*{4}}
   \put(30,6){\circle*{4}}
   \put(50,6){\circle*{4}}
   \put(70,6){\circle*{4}}
\put(38,8){\small{${}_{1}$}}

\end{picture}
\begin{picture}(127,40)(-10,0)

\put(50,30){\vector(0,-1){2}}

\put(0,26)\halfdotline
   \put(12,26)\dotline
   \put(32,26)\dotline
   \put(52,26)\dotline
   \put(72,26)\halfsolid
\put(10,26){\circle{4}}
   \put(30,26){\circle{4}}
   \put(50,26){\circle{4}}
   \put(70,26){\circle{4}}

\put(0,16)\halfsolid
  \put(12,16)\dotline
   \put(32,16)\dotline
   \put(52,16)\dotline
   \put(72,16)\halfdotline
\put(10,16){\circle{4}}
   \put(30,16){\circle{4}}
   \put(50,16){\circle{4}}
   \put(70,16){\circle{4}}

\put(0,6)\halfdotline
  \put(12,6)\dotline
  \put(32,6)\dotline
  \put(52,6)\dotline
  \put(72,6)\halfdotline
\put(10,6){\circle{4}}
   \put(30,6){\circle{4}}
   \put(50,6){\circle{4}}
   \put(70,6){\circle{4}}
\put(38,8){\small{${}_{1}$}}

\end{picture}
\begin{picture}(120,40)(-10,0)

\put(30,30){\vector(0,-1){2}}

\put(0,26)\halfdotline
   \put(12,26)\dotline
   \put(32,26)\halfsolid
\put(10,26){\circle{4}}
   \put(30,26){\circle{4}}

\put(0,16)\halfsolid
  \put(12,16)\dotline
   \put(32,16)\halfdotline
\put(10,16){\circle{4}}
   \put(30,16){\circle{4}}

\put(0,6)\halfdotline
  \put(12,6)\dotline
  \put(32,6)\halfdotline
\put(10,6){\circle{4}}
   \put(30,6){\circle{4}}
\put(18,8){\small{${}_{1}$}}

\end{picture}

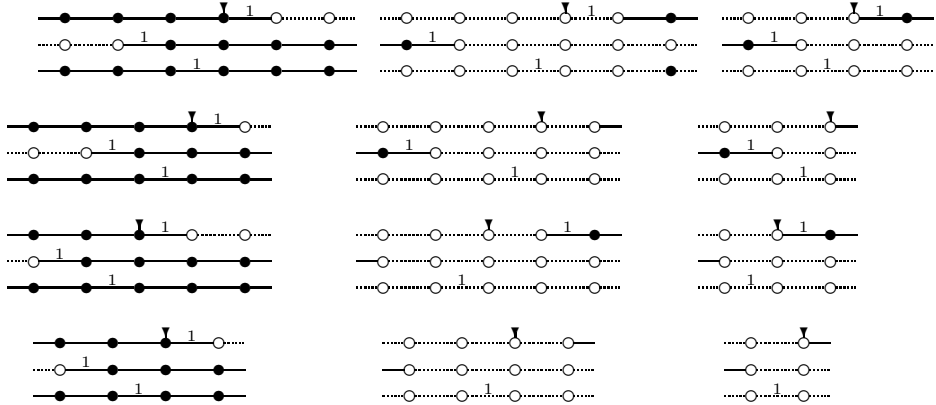
{\small Figure 11. The graphes represent the stable modifications
of type $I_b^{+1}$ sheaves over $\tilde\bD_0\cap\tilde\bD_2$, over
$\tilde\bD_2-\tilde\bD_0$, over $\tilde\bD_0-\tilde\bD_2$ or away
from $\tilde\bD_0\cup\tilde\bD_2$.}

\vsp

\begin{lemm}
Let $\tilde\EE$ over $\tilde\cX$ be the result of the stable
elementary modification of $p\sta\EE\st$, then all its members are
weakly $0$-stable and hence their stabilization are $0$-stable.
\end{lemm}

\begin{proof}
The proof is similar to the argument before and will be omitted.
\end{proof}

The family $\tilde\EE\st$ over each $\tilde\bV$ induces a morphism
$\tilde\bV\to\bM^0$ that is the local extension of the birational
$\tilde\bM\,-\!\to\bM^0$. Since both $\tilde\bM$ and $\bM^0$ are
smooth, the local extension patch together to form a morphism
$\tilde\bM\to\bM^0$, as desired.


\subsection{The existence of the descent $\tilde\bM\to \bM^0$}

In this subsection, we will show that the morphism $\Psi\mh
\tilde\bM\to\bM^0$ descends to a morphism $\bM_2\to\bM^0$.

We begin with a brief outline of our strategy. First, we know that
the flipped loci $\tilde\II_b$ is $\PP W_+$ and is the result of
$\bE_b$ after contracting $\bE_a\cap\bE_b$, and that the
exceptional divisor of $\tilde\bM_2\to\bM_1$ over $\tilde\II_b$ is
$\PP W_+\times_{\PFT}\PP W_-$. Since $\bM_2$ is a flip of $\bM_1$,
the projection $\Psi_4$ is the result of contracting all fibers of
$\PP W_+$. On the other hand, by our description of the
contraction $\tilde\bM_1\to\bM_1$, the exceptional divisor of
$\tilde\bM\to \tilde\bM_1$ over $\Psi\upmo_3(\tilde\II_b)$ is
$\bE_b\times_{\bM_1}\PP W_-$. Since both $\bM_2$ and $\bM^0$ are
smooth and since $\bE_b$ is proper, to show that $\Psi$ descends
to an morphism $\bM_2\to\bM^0$ it suffices to show that there is
an open subset $\bU\sub\bE_b$ and an open $\bV\sub \PP W^-$ so
that the restriction of $\Psi$ to $\bU\times_{\bM_1} \bV$ is a
composition of the second projection $\bU\times_{\bM_1}\bV\to\bV$
with a morphism $\bV\to \bM^0$.

To prove the last statement, we need a description of the normal
bundle $N_{\PP W_+/\bM_1}$ that relates directly to the elementary
modification we shall perform. It is expressed in terms of
relative extension sheaves; hence a tautological family on $\bM_1$
is required. There is one more technical difficulty: the space
$\bM_1$ is not a moduli space per se, thus we can not use
deformation theory to derive its tangent bundle. Nevertheless,
$\bM_1$ is birational to $\bM\uhalf$, and thus over a dense open
subset its tangent bundle is given by the deformation theory of
sheaves.

Our first step is to construct a tautological family over an open
subset of $\PP W_+\sub \bM_1$. Since $\PP W_+$ is a projective
bundle over $\PFT$, we shall content ourselves with constructing
such a family over an open subset of the fiber $\PP W_{+\eta}$ of
$\PP W_+$ over a general $\eta\in\PFT$. Such family will be the
universal extension of a sheaf $\cF_\eta$ by another sheaf
$\cL_\eta$ over $X_{1,1}\shar$.

We begin with constructing $\cF_\eta$ and $\cL_\eta$. Let
$\eta_0\in\Xi_1$ be any point associated to a pair of GPBs
$(F^G,L^G)$ with $F^0=\ell\oplus F|_{p_2}$ for a line $\ell\sub
F|_{p_1}$, and let $\beta_1\mh F\to F|_{p_1}/\ell=\kk(p_1)$ be the
induced homomorphism. Let $\tilde X_{\Ao}$ be the blowing up of
$(p_1,0)\in X\times\Ao$, let $\varphi\mh \tilde X_{\Ao} \to X$ be
the projection, and define
$$\tilde F=\ker\{\varphi\sta F\mapright{\varphi\sta\beta_1}\varphi\sta\kk(p_1)\}
|_{\tilde X_0},\quad \text{where}\ \tilde X_0=\tilde
X_{\Ao}\times_{\Ao}0.
$$
The sheaf $\tilde F$ is a locally free sheaf on $X_{1,0}\cong
\tilde X_0$ whose restriction to the unique rational curve
$D_1\sub X_{1,0}$ is isomorphic to $\cO\oplus\cO(1)$. To obtain a
sheaf on $X_{1,1}$ we consider the map $\varphi_2\mh X_{1,1}\to
X_{1,0}$ contracting the rational curve $D_2\sub X_{1,1}$ attached
to $p_2$. By abuse of notation, we still denote by $D_1$ the other
rational curve in $X_{1,1}$. The sheaf we intend to construct is
the direct sum
\begin{equation}
\lab{4.111} F\pri=\varphi_2\sta \tilde F\oplus \cO_{D_1}\oplus
\cO_{D_2}(1),
\end{equation}
where $\cO_{D_i}=\iota_{i\ast}\cO_{D_i}$ with $\iota_i\mh D_i\to
X_{1,1}$ the inclusion.

Our next step is to glue $F\pri$ along the two marked points of
$X_{1,1}$ using a lift $\eta\in \PFT$ of $\eta_0\in\Xi_1$ that is
defined by a homomorphism $\beta_2\mh F\to \kk(p_2)$. Let $\psi\mh
X_{1,1}\to X_{1,1}\shar$ be the obvious morphism, let
$q_-=D_1\cap\psi\upmo(q\shar)$ and $q_+= D_2\cap
\psi\upmo(q\shar)$ be the two marked points of $X_{1,1}$, and let
$q_0=D_1\cap X$ and $q_2=D_2\cap X$. We consider the space
$K_-=\Hom_{X_{1,1}}(\cO_{D_1},F\pri)$ and the subspaces of
$F\pri|_{q_-}$:
$$V_{-,1}=\{ f(q_-)\mid f\in K_-, f(q_0)=0\}\and
V_{-,2}=\{ f(q_-)\mid f\in K_-\}.
$$
Obviously,
$$0\ne V_{-,1}\varsubsetneq V_{-,2}\varsubsetneq F\pri|_{q_-}
$$
form a filtration that depends only on $F^G$. Similarly, we define
a filtration
$$0\ne V_{+,1}\varsubsetneq V_{+,2}\varsubsetneq F\pri|_{q_+}
$$
via
$$V_{+,1}=\ker\{\varphi_2\sta\tilde F|_{q_+}\mapright{\beta_2}
\varphi_2\sta\kk(p_2)|_{q_+}\} \and V_{+,2}=V_{+,1}\oplus
\cO_{D_2}(1)|_{q_+}.
$$
This filtration depends on $\eta$. After that, we pick an
isomorphism $h\mh F\pri|_{q_-}\to F\pri|_{q_+}$ that preserves
these two filtrations, and identify $F\pri$ along the two marked
points $q_-$ and $q_+$ via this isomorphism to form a vector
bundle on $X_{1,1}\shar$. We denote the resulting vector bundle by
$\cF_h$.

Given two such homomorphisms $h_1$ and $h_2$, we say
$\cF_{h_1}\sim \cF_{h_2}$ if there is an automorphism $\sigma$ of
$X_{1,1}\shar$ and an isomorphism $\sigma\sta \cF_{h_2}\cong
\cF_{h_1}$.

\begin{lemm}
\lab{4.112} Let $\eta\in\PFT$ be any element and let $\cF_h$ be
the sheaf on $X_{1,1}\shar$ so constructed. Then $\cF_h$,
modulo the equivalence relation so defined, is independent of the
choice of $h$. Let $\cF_\eta$ be a representative of this
equivalence class. Then for any $\eta\pri\in\PFT$,
$\cF_\eta\sim\cF_{\eta\pri}$ if and only if $\eta=\eta\pri$.
\end{lemm}

\begin{proof}
Let $G$ be the group of pairs $(v,\sigma)$ where $\sigma$ is an
automorphism of the pointed curve $(X_{1,1}, q_\pm)$ and $v$ is an
isomorphism $\sigma\sta F\pri\mapright{v} F\pri$. It is direct to
check that the tautological homomorphism $G\to \Aut(F\pri|_{p_2})$
preserves the filtration $V_{+,\bullet}$ while the image of $G\to
\Aut(F\pri|_{p_1})$ is exactly the subgroup of automorphisms that
preserve the filtration $V_{-,\bullet}$. It follows that the
equivalence class of the sheaf $\cF_h$ is independent of the
choice of $h$.

The proof of the second part is straightforward and will be
omitted.
\end{proof}

Finally, let $\iota_0\mh X\to X_{1,1}\shar$ be the tautological
inclusion and let $\cL_\eta=\iota_{0\ast}L(-p_1-p_2)$.

Next, we will construct a vector space $\cW_\eta$ and a family of
sheaves over $\PP\cW_\eta$. Later we will show that $\cW_\eta$ is
canonically isomorphic to $W_{+\eta}$ and the family over
$\PP\cW_\eta$ and the tautological family over $\PP W_{+\eta}$
coincide over a dense open subset of $\PP \cW_\eta\cong\PP
W_{+\eta}$.

The vector space $\cW_\eta$ is the kernel of the canonical
homomorphism
$$\Ext^1_{X_{1,1}\shar}(\cF_\eta,\cL_\eta)\lra H^0(\ext^1(\cF_\eta,\cL_\eta))
\equiv H^0(\kk(q_0)\oplus \kk(q_2))\mapright{\text{proj}}
\kk(q_2).
$$

Next we construct a family of curves over $\PP\cW_\eta$. Let
$\cD_1=\Po\times\PP\cW_\eta$, let $\cO(1)$ be the degree one line
bundle over $\PP\cW_\eta$ and let $\cD_2=\PP(\cO\oplus\cO(1))$ be
the associated projective bundle over $\PP\cW_\eta$. We fix two
sections
$$\cQ_-=0\times \PP\cW_\eta\and \cQ_0=\infty\times\PP\cW_\eta
$$
of $\cD_1$  and pick two sections
$$\cQ_2=\PP(0\oplus\cO(1))\and \cQ_+=\PP(\cO\oplus 0)
$$
of $\cD_2$. We then glue $\cD_1$ to $X\times\PP\cW_\eta$ by
identifying $\cQ_0$ with $p_1\times\PP\cW_\eta$, and then glue
$\cD_2$ to $\cX_{1,0}$ by identifying $\cQ_2$ with
$p_2\times\PP\cW_\eta$. We denote the family from the first gluing
by $\cX_{1,0}$ and denote the family resulting from both gluing by
$\cX_{1,1}$. The first is a constant family of $X_{1,0}$ and the
second is a non-constant family of $X_{1,1}$ over $\PP\cW_\eta$.
Let $\cX_{1,1}\shar$ be the result of gluing the two sections
$\cQ_-$ and $\cQ_+$ of $\cX_{1,1}$. It is a family of
$X_{1,1}\shar$ over $\PP \cW_\eta$.

We now construct a sheaf $\FF$ over $\cX_{1,1}\shar$. Let
$\cX_{1,1}\to\cX_{1,0}$ be the contraction of the component
$\cD_2$, and let $\Phi: \cX_{1,1}\to X_{1,1}$ be the composition
of the contraction $\cX_{1,1}\to\cX_{1,0}$ with the projection
$\cX_{1,0}\to X_{1,0}\sub X_{1,1}$. We consider the sheaf of
$\cO_{\cX_{1,1}}$-modules
\begin{equation}
\lab{4.25} \FF\pri=\Phi\sta \varphi\sta_2\tilde F\oplus
\cO_{\cD_1}\oplus\cO_{\cD_2}(\cQ_2),
\end{equation}
according to the convention of (\ref{4.111}). Clearly, there are
canonical isomorphisms
$$\FF\pri|_{\cQ_\pm}\cong
F\pri|_{q_\pm}\otimes_{\kk}\cO_{\cQ_\pm}.
$$
Hence any isomorphism $h\mh F\pri|_{q_-}\cong F\pri|_{q_+}$
induces a canonical isomorphism $\tilde h\mh \FF\pri|_{\cQ_-}\cong
\FF\pri|_{\cQ_+}$. Using an isomorphism $h$ that preserves the two
filtrations as defined in Lemma \ref{4.112}, we can glue $\FF\pri$
along the two marked sections $\cQ_-$ and $\cQ_+$ to obtain a new
sheaf $\FF$ over $\cX_{1,1}\shar$. Clearly, restricting to each
fiber of $\cX_{1,1}$ over $\PP\cW_\eta$, the sheaf $\FF\pri$ is
merely the sheaf $F\pri$ constructed before, and the sheaf $\FF$
is isomorphic to the $\cF$ constructed in Lemma \ref{4.112}.

We are ready to construct the desired tautological family $\EE$
over $\cX_{1,1}\shar$. First, recall that $\FF\pri$ has a direct
summand $\cO_{\cD_2}(\cQ_2)$. Since
$\cO_{\cD_2}(\cQ_2)|_{\cQ_2}\cong\cO_{\cQ_2}(-1)$\footnote{In case
$\pi\mh Z\to\PP\cW_\eta$ is a family and $\cE$ is a sheaf of
$\cO_Z$-modules, we will use $\cE(1)$ to denote the sheaf
$\cE\otimes\pi\sta\cO_{\PP \cW_\eta}(1)$.}, the inclusion
$\cO_{\cD_2}(\cQ_2)\sub \FF\pri$ defines a subsheaf
$$g_1\mh\cO_{\cQ_2}(-1)\to \FF\pri|_{\cQ_2}\equiv \FF|_{\cQ_2}.
$$
Let $\LL$ be the sheaf of $\cO_{\cX_{1,1}\shar}$-modules
$\LL=\bar\iota\lsta L(-p_1-p_2)$, where $\bar\iota\mh
X\times\PP\cW_\eta\to \cX_{1,1}$ is the tautological inclusion.
Next let $\pi\mh \cX_{1,1}\shar\to\PP\cW_\eta$ be the projection
and consider the following two relative extension sheaves and the
natural homomorphism between them:
\begin{equation}
\lab{4.113} \ext_{\cX_{1,1}\shar/\PP\cW_\eta}^1(\LL(-1), \FF)\lra
\ext_{\hat\cQ_2}^1(\hat\LL(-1),\hat\FF),
\end{equation}
where the latter is canonically isomorphic to
\begin{equation}
\lab{4.23} \Ext^1_{W[2]_0\shar}(\cL_\eta,\cF_\eta)\otimes
\cO_{\PP\cW_\eta}(-1).
\end{equation}
Here $\hat\cQ_2$ is the formal completion of $\cX_{1,1}\shar$
along $\cQ_2$ while
$\hat\LL=\LL\otimes_{\cO_{\cX_{1,1}\shar}}\!\cO_{\hat\cQ_2}$, etc.
Because the support of $\LL$ only intersects $\cD_2$ along $\cQ_2$
and the restriction of $\FF$ to
$X_{1,0}\times\PP\cW_\eta\sub\cX_{1,1}\shar$ is the pull-back of
$F\pri|_{X_{1,0}}$, the kernel of the above homomorphism is
canonically isomorphic to $\cW_\eta\otimes\cO_{\PP\cW_\eta}(1)$.
Let $\eps$ be a tautological section of
$\cW_\eta\otimes\cO_{\PP\cW_\eta}(1)$. It can be viewed as a
section of the relative extension sheaf, thus defines an extension
sheaf $\EE\pri$ fitting into the exact sequence
$$0\lra \FF\lra \EE\pri\lra \LL(-1)\lra 0.
$$
Because $\eps$ is in the kernel of (\ref{4.113}), $\EE\pri$ is not
locally free along $\cD_2$. Not only that, there is a homomorphism
$\mu_1\mh \EE\pri\lra \cO_{\cQ_2}(-1)$ so that the composite
$\FF\to\EE\pri\to\cO_{\cQ_2}(-1)$ is identical to the composite
$\FF\to \cO_{\cD_2}(\cD_2)|_{\cQ_2}\equiv \cO_{\cQ_2}(-1)$ induced
by (\ref{4.25}). Let $\mu_2\mh \EE\pri\to \cO_{\cQ_2}(-1)$ be
induced by $\EE\pri\to\LL(-1)|_{\cQ_2}$ and define $\EE\ppri$ by
the exact sequence
$$0 \lra \EE\ppri \lra \EE\pri\mapright{(\mu_1,\mu_2)} \cO_{\cQ_2}(-1) \lra 0.
$$
The resulting sheaf $\EE\ppri$ is locally free along $\cQ_2$.

The sheaf $\EE\pri$ is still non-locally free along some points of
$\cQ_0$. Indeed, the section $\eps$ composed with the homomorphism
$$\ext_{\cX_{1,1}\shar/\PP\cW_\eta}^1(\LL(-1),\FF)
\lra \ext_{\hat\cQ_0}^1(\hat\LL(-1),\hat\FF) \cong\cO_{\cQ_0}(1)
$$
defines a section of $\cO_{\cQ_0}(1)$ whose vanishing locus is
exactly where $\EE\ppri$ is not locally free. Let $s\in
H^0(\cO_{\cQ_0}(1))$ be this section. Before we proceed, we need
to resolve the non-locally freeness of $\EE\ppri$. We first glue
$\cD_2$ to $X\times\PP\cW_\eta$ by identifying $\cQ_2$ with
$p_2\times\PP\cW_\eta$, and then glue $\cD_1$ to the resulting
family by identifying $\cQ_-$ with $\cQ_+$ in the obvious way. Let
$\varphi_1\mh\cX\pri\to\cX_{1,1}\shar$ be the projection, which is
the smoothing of $\cX_{1,1}\shar$ along $\cQ_0\sub\cX_{1,1}\shar$.
Clearly, $\cX\pri$ is a locally constant family of $X_{0,2}$'s. We
next blow up $\cX\pri$ along $p_1\times s\upmo(0)$, where $s$ is
the section of $\cO_{\PP\cW_\eta}(1)$ mentioned before. We denote
the blowing up by $\tilde\cX$, and let $\cQ_0\pri\sub \tilde\cX$
be the proper transform of $p_1\times\PP\cW_\eta\sub\cX\pri$.
Lastly, we construct a new family $\cX_{1,1}\shar$ by identifying
(gluing) the two sections $\cQ_0\pri$ and $\cQ_0$ of $\cX\pri$ in
the obvious way, and we keep the section $\cQ_-=\cQ_+\sub\cX\shar$
as its marked section. This way $\cX\shar$ is a family whose
members are either $X_{1,1}\shar$ or $X_{2,1}\shar$. Let
$\varphi\mh \cX\shar\to \cX_{1,1}\shar$ be the tautological
projection.

\begin{lemm}
There is a unique family of locally free sheaves $\tilde\EE$ over
$\cX\shar$ so that $\varphi\lsta\tilde\EE=\EE\dpri$ and
$R^1\varphi\lsta\tilde\EE=0$.
\end{lemm}

\begin{proof}
The proof is straightforward and will be omitted.
\end{proof}

For $\xi\in\PP\cW_\eta$ we denote by $\tilde\EE_\xi$ the
restriction of $\tilde\EE$ to the fiber $\cX\shar_\xi$ of
$\cX\shar$ over $\xi$.

\begin{lemm}\lab{43.1}
There is a line $\Sigma_\eta\sub\PP\cW_\eta$ so that for each
$\xi\in\PP\cW_\eta-\Sigma_\eta$ the sheaf $\tilde\EE_\eta$ is
$\half$-stable.
\end{lemm}

\begin{proof}
We will postpone the proof until the next subsection.
\end{proof}

Since for $\xi\in\PP\cW_\eta-\Sigma_\eta$ the sheaf
$\tilde\EE_\xi$ is $\half$-stable, by the universal property of
$\bM\uhalf$ the family $\tilde\EE$ induces a canonical morphism
$u\mh \PP\cW_\eta-\Sigma_\eta \to \bM\uhalf$. Further, by the
construction of the family, it is clear that $u$ factor through
$\II_b-\II_a$, and hence to $\tilde\II_b-\tilde\II_a\sub\bM_1$.
Since $\tilde\II_b=\PP W_+$, the morphism $u$ induces a morphism
$$\tilde u: \PP\cW_\eta-\Sigma_\eta\lra \PP W_+.
$$
On the other hand, the construction of the family $\tilde\EE$
ensures that the composition of $\tilde u$ with the projection
$\PP W_+\to \PP\cF|_{p_2}\dual$ maps $\PP\cW_\eta-\Sigma_\eta$ to
the point $\eta\in\PP\cF|_{p_2}\dual$. Thus $\tilde u$ factor
through
$$u_\eta: \PP\cW_\eta-\Sigma_\eta \lra \PP W_{+\eta}.
$$

\begin{lemm}
The morphism $u_\eta$ extends to an isomorphism $\PP\cW_\eta\to
\PP W_{+\eta}$.
\end{lemm}

\begin{proof}
By construction $u_\eta$ is one-one. Since $\dim\PP\cW_\eta
=\dim\PP W_{+\eta}\geq 3$, $u_\eta$ is an isomorphism away from a
line. It is direct to check that $u_\eta$ maps lines in
$\PP\cW_\eta-\Sigma_\eta$ to lines in $\PP W_{+\eta}$. Hence
$u_\eta$ automatically extends to an isomorphism.
\end{proof}

By the argument in Section 3, there is a vector bundle $W_-$ over
$\PP\cF|_{p_2}\dual$ so that the normal bundle
$$N_{\tilde\II_b/\bM_1}\cong\varphi\sta W_-\otimes\cO_{\PP W_+}(-1),
$$
where $\varphi\mh \PP W_+\to \PP\cF|_{p_2}\dual$ is the
projection.

\begin{lemm}\lab{43.2}
The normal bundle $N_{\tilde\II_b/\bM_1}|_{\PP\cW_\eta}$ is
canonically isomorphic to
$$\Ext^1_{W[2]\shar}(\cL_\eta,\cF_\eta)\otimes \cO_{\PP \cW_\eta}(-1).
$$
\end{lemm}

\begin{lemm}\lab{43.3}
The restriction of the morphism $\Psi$ to the preimage of $\PP
W_{+\eta}$, say $\phi_\eta\mh
\PP\cW_\eta\times\PP\Ext^1_{W[2]\shar}(\cL_\eta,\cF_\eta)\to
\bM^0$ is the composite of the second projection with a morphism
$h\mh \PP\Ext^1_{W[2]\shar}(\cL_\eta,\cF_\eta)\to \bM^0$.
\end{lemm}

\begin{proof}
We will postpone the proofs of these two Lemmas until the next
subsection.
\end{proof}

Since $\bM_2$ is a flip of $\bM_1$ along
$\tilde\II_b\cup\tilde\II_c$, the restriction of $\Psi_3$ to the
exceptional divisor over $\tilde\II_b$, which is $\PP W_+\times\PP
W_-$, is the composite of the second projection with the morphism
$\PP W_-\to\bM_2$. In particular, this proves that

\begin{lemm}
For any $z\in \image\phi_\eta$, the image set $(\Psi_4\circ \tilde
\Psi_2)(\phi_\eta\upmo (z))$ is a single point set in $\bM_2$.
\end{lemm}

Since $\eta\in\PP\cF_\eta|_{p_2}\dual$ is an arbitrary point, this
proves

\begin{lemm}
For any closed $z\in\bM_2$ there is a unique point $z\pri\in\bM^0$
so that $(\Psi_4\circ \tilde \Psi_2)\upmo(z)=\Psi\upmo(z\pri)$.
\end{lemm}

As a corollary, this proves the equivalence result we set out to
prove:

\begin{prop}
The induced birational map $\bM^0\sim \bM_2$ is an isomorphism of
varieties.
\end{prop}

\subsection{The proof of Lemmas \ref{43.1}, \ref{43.2} and \ref{43.3}}

In this subsection, we will give the proof of Lemmas \ref{43.1},
\ref{43.2} and \ref{43.3}.

\begin{proof}[Proof of Lemma \ref{43.1}]
We need to investigate when the sheaf $\tilde\EE_\xi$ is
$\half$-stable. For the moment, we assume $\xi$ is away from the
vanishing locus $s\upmo(0)$. Then $\tilde\EE_\xi$ is of type
$\II_b^{+0}$ that fits into the exact sequence
\begin{equation}
\lab{4.21} 0\lra \cF_\eta\lra \tilde\EE_\xi\lra \cL_\eta\lra 0.
\end{equation}
Following the discussion in Section 3, $\tilde\EE_\xi$ is not
$\half$-stable if and only if there is a sheaf $\bar\cL_\eta$ that
is locally free away from the marked node $q_1=q\shar$ so that it
fits into the diagram with exact rows:
$$\begin{CD}
0 @>>> \cF_\eta @>>> \tilde\EE_\xi @>>> \cL_\eta @>>> 0\\
@.@AA{\subset}A @AAA @|\\
0@>>> \cO_{D_1}(-1)\oplus\cO_{D_2}(-1) @>>> \bar\cL_\eta @>>>
\cL_\eta @>>> 0.
\end{CD}
$$
Because $\bar\cL_\eta$ is unique and the left square is a
push-out, $\bar\EE_\xi$ is uniquely determined by the left
vertical inclusion. On the other hand, the subsheaf $\cO_{D_2}(-1)
\hookrightarrow \cF_\eta$ is unique and there is a $\Po$ family of
subsheaves $\cO_{D_1}(-1) \hookrightarrow\cF_\eta$. Hence there is
a $\Po$ family of extensions (\ref{4.21}) that are derived from
the diagram above. Further, it is easy to see that there is one
choice of $\cO_{D_1}(-1) \hookrightarrow\cF_\eta$ so that the
associated sheaf $\tilde \EE_\xi$ is not locally free at $q_0$. A
quick reasoning shows that this corresponds exactly to the case
where $\xi\in s\upmo(0)$. Combined, this shows that there is a
line $\Sigma_\eta\sub\PP\cW_\eta$ so that for all $\xi\in
\PP\cW_\eta-\Sigma_\eta\cup s\upmo(0)$ the associated sheaves
$\tilde\EE_\xi$ are $\half$-stable.

In case $\xi\in s\upmo(0)$, a similar argument shows that
$\tilde\EE_\xi$ is $\half$-stable unless $\xi\in s\upmo(0)\cap
\Sigma_\eta$. This proves the lemma.
\end{proof}

\begin{proof}[Proof of Lemma \ref{43.2}]
We now prove that
\begin{equation}
\lab{44.1} N_{\tilde\II_b/\bM_1}|_{\PP\cW_\eta}\cong
\Ext^1_{W[2]\shar}(\cL_\eta,\cF_\eta)\otimes \cO_{\PP
\cW_\eta}(-1).
\end{equation}
Let $W[2]\shar/\At$ be the family of marked curves containing
$W[2]_0\shar \cong X_{1,1}\shar$ as its central fiber. Then
$\tilde\EE_\xi$, where
$\xi\in\PP\cW_\eta-s\upmo(0)\cup\Sigma_\eta$, is a sheaf over
$W[2]\shar_0$. It is known that the first order deformations of
$\tilde\EE_\xi$ as sheaves of $\cO_{w[2]\shar}$-modules is
$\Ext_{W[2]\shar}^1(\tilde\EE_\xi,\tilde\EE_\xi)$, which fits into
the diagram
$$
\begin{CD}
\Ext_{W[2]\shar}^1(\cL_\eta,\cF_\eta) @. @. \Ext_{W[2]\shar}^1(\cL_\eta,\cL_\eta)\\
@VVV @. @VVV\\
\Ext_{W[2]\shar}^1(\tilde\EE_\xi,\cF_\eta) @>>>
\Ext_{W[2]\shar}^1(\tilde\EE_\xi,\tilde\EE_\xi)@>{\phi_1}>>
\Ext_{W[2]\shar}^1(\tilde\EE_\xi,\cL_\eta) \\
@VVV @. @VV{\phi_2}V\\
\Ext_{W[2]\shar}^1(\cF_\eta,\cF_\eta) @. @. \Ext_{W[2]\shar}^1(\cF_\eta,\cL_\eta)\\
\end{CD}
$$
Because the standard $(\CC\sta)^{\times 2}$ action on $\At$ lifts
to an action on $W[2]\shar\to \At$, it induces a homomorphism
$\CC^{\oplus 2}\equiv T_0\At\to
\Ext_{W[2]\shar}^1(\tilde\EE_\xi,\tilde\EE_\xi)$. Since
$\tilde\EE_\xi$ is $\half$-stable, $[\tilde\EE_\xi] \in\bM\uhalf$,
and hence lies in $\II_b^+$. Then the tangent space
$T_{[\tilde\EE_\xi]}\bM\uo$ at $[\tilde\EE_\xi]$ is canonically
isomorphic to
$\Ext_{W[2]\shar}^1(\tilde\EE_\xi,\tilde\EE_\xi)/T_0\At$.

We now claim that the kernel of $\phi=\phi_2\circ\phi_1$ contains
the image of $T_0\At\to
\Ext_{W[2]\shar}^1(\tilde\EE_\xi,\tilde\EE_\xi)$ and the tangent
space of $\II^+_b$ at $[\tilde\EE_\xi]$ is the quotient
$\ker(\phi)/ T_0\At$. Indeed, the groups
$\Ext_{W[2]\shar}^1(\cL_\eta,\cL_\eta)$,
$\Ext_{W[2]\shar}^1(\cF_\eta,\cF_\eta)$ and
$\Ext_{W[2]\shar}^1(\cL_\eta,\cF_\eta)$ parameterize the first
order deformations of $\cL_\eta$, of $\cF_\eta$ and the space of
extensions of $\cL_\eta$ by $\cF_\eta$. It is direct to check that
the kernel of $\phi$ is the tangent space at $[\tilde\EE_\xi]$ of
the space of all sheaves of type $\II^+_b$. Because the
$(\CC\sta)^{\times 2}$ action preserves this space, we have
$\image(T_0\At)\sub\ker(\phi)$ and hence $\ker(\phi)/T_0\At\equiv
T_{[\tilde\EE_\xi]}\II^+_b$.

We now show that $\phi$ is surjective. Once this is established,
then the normal vector space to $\II_b^+$ at $\tilde\EE_\eta$ is
canonically isomorphic to
\begin{equation}
\lab{4.22} N_{\II_b^+/\bM^{1/2}}|_{\tilde\EE_\xi}\equiv
\Ext^1_{W[2]\shar}(\cF_\eta,\cL_\eta).
\end{equation}
First, since $T_{[\tilde\EE_\phi]}\II^+_a=\ker(\xi)/T_0\At$, the
image of $\phi$ is the normal vector space to $\II^+_a$ in
$\bM\uo$ at $[\tilde\EE_\xi]$. By (\ref{3.21}), we know that the
normal vector space has dimension $2g-1$. Thus to prove the lemma
it suffices to show that
$$\dim \Ext_{W[2]\shar}^1(\cF_\eta,\cL_\eta)=2g-1.
$$
Recall that $q\shar=q_1$. By a direct computation we have the
exact sequence
$$
0\lra \kk(q_0)\oplus \kk(q_2)\lra
\ext_{W[2]\shar}^1(\cF_\eta,\cL_\eta)\lra
\hom(\cF_\eta,\cL_\eta)\otimes_{\kk} T_0\At \lra 0.
$$
We claim that $H^0(\hom(\cF_\eta,\cL_\eta))=0$. First, let
$F^G=(F,F^0)$ and $L^G=(L,0)$ be the associated GPB vector bundles
of $[\tilde\EE_\xi]$ in $\bG^{1/3}_{2,4,3}$ and
$\bG^{1/3}_{1,3,0}$, respectively, that is the image of the
morphism $\II^+_a\to \bG^{1/3}_{2,4,3}\times\bG^{1/3}_{1,3,0}$
introduced in section 3.2. By abuse of notation, we let $j\mh X\to
W[2]_0\shar$ be the main irreducible component. Then we have exact
sequences of sheaves of $\cO_X$-modules
$$0\lra \cF_\eta\otimes_{\cO_{W[2]\shar}}\cO_X\lra F\lra \kk(p_1)\lra 0
$$
and
$$0\lra \cL_\eta\otimes_{\cO_{W[2]\shar}}\cO_X\lra L\lra \kk(p_1)\oplus\kk(p_2)\lra 0.
$$
Further, a direct check shows that
$$\hom(\cF_\eta,\cL_\eta)=j\lsta\ker\{F\dual\otimes L\lra \kk(p_1)\otimes L|_{p_1}
\oplus F\dual\otimes L|_{p_2}\}.
$$
Hence any nontrivial homomorphism $F\to L$ in
$H^0(\hom(\cF_\eta,\cL_\eta))$ is a homomorphism $F^G\to L^G$ of
GPBs. But both $F^G$ and $L^G$ are $\frac13$-stable GPBs and thus
there is no nontrivial homomorphisms between them.
This proves that $H^0(\hom(\cF_\eta,\cL_\eta))=0$.
Combined with the exact sequence
$$
0\lra H^1(\hom(\cF_\eta,\cL_\eta))\lra
\Ext^1_{W[2]\shar}(\cF_\eta,\cL_\eta)\lra
H^0(\ext^1_{W[2]\shar}(\cF_\eta,\cL_\eta))\lra 0,
$$
we obtain
$$\dim\Ext^1_{W[2]\shar}(\cF_\eta,\cL_\eta)=2+h^1(\hom(\cF_\eta,\cL_\eta)
=2g-1.
$$
This proves that the arrow $\phi$ is surjective.

Because the isomorphism (\ref{4.22}) is canonical, restricting to
$\PP\cW_\eta-s\upmo(0)\cup\Sigma_\eta$ we have canonical
isomorphism
$$N_{\II_b^+/\bM^{1/2}}|_{\PP\cW_\eta-s\upmo(0)\cup\Sigma_\eta}
\cong
\ext_{\cX_{1,1}\shar/\PP\cW_\eta}^1(\LL(-1),\FF)|_{\PP\cW_\eta-s\upmo(0)\cup\Sigma_\eta}.
$$
Because $\codim \Sigma_\eta\geq 2$,  $s\upmo(0)$ is a hypersurface
and  both $N_{\II_b^+/\bM^{1/2}}$ and
$\ext_{\cX_{1,1}\shar/\PP\cW_\eta}^1(\LL(-1),\FF)$ are of the
forms $\CC^{2g-1}\otimes \cO_{\PP\cW_\eta}(-1)$, the above
isomorphism must extend to an isomorphism (\ref{44.1}), as
desired. This proves the lemma.
\end{proof}

\begin{proof}[Proof of Lemma \ref{43.3}]
We pick an element $(\xi,v)\in
\PP\cW_\eta\times\PP\Ext^1_{W[2]\shar}(\cL_\eta,\cF_\eta)$. We
assume that the sheaf $\tilde\EE_\xi$ is a sheaf over
$W[2]\shar_0$, which is $\half$-stable and the image of $v$ in
$$H^0\bl\ext^1_{W[2]\shar}(\cF_\eta,\cL_\eta)\br =\kk(q_0)\oplus
\kk(q_2)
$$
is not contained in either $\kk(q_0)$ or in $\kk(q_2)$. Now let
$B=\spec\kk[u]/(u^2)$ and let $B_0\sub B$ be the closed point.
Then a lift $\tilde
v\in\PP\Ext^1_{W[2]\shar}(\tilde\EE_\xi,\tilde\EE_\xi)$ of $v$,
namely $\phi(\tilde v)=v$, defines a sheaf of
$\cO_{W[2]\shar\times B}$-modules $\tilde\EE_\xi(\tilde v)$ that
is the extension of $\tilde\EE_\xi$ by $\tilde\EE_\xi\otimes \cI$
defined by the class $\tilde v$. Here $\cI$ is the ideal sheaf of
$W[2]\shar\times B_0\sub W[2]\shar\times B$. Then a direct local
calculation of extension sheaves shows that there is an embedding
$\omega\mh B\to \At$ that does not lie in the two coordinate lines
of $\At$, and so $\tilde\EE_\xi(\tilde v)$ is a rank three locally
free sheaf of $\cO_{W[2]\shar\times_{\At}B}$-modules.

We now suppose $\Psi_2\upmo(\xi)$ is a single point. Then
$(\xi,v)\in \tilde\bM_1$ lifts to a unique element in $\tilde\bM$,
which we denote by $(\xi,v)$ as well. Then following the
discussion in subsection 4.1, the image $\Psi((\xi,v))\in\bM^0$ is
the point associated to the sheaf $\tilde\FF_\xi(\tilde v)$ that
was constructed by first taking the kernel of the composite
$$\FF_\xi(\tilde v)=\ker\bigl\{\tilde\EE_\xi(\tilde v)\lra \tilde\EE_\xi(\tilde v)
\otimes_{\cO_{W[2]\shar\times_{\At}B}}\cO_{W[2]\shar_0}\cong
\tilde \EE_\xi\lra\cL_\eta\bigr\}
$$
and then restricting to the closed fiber $W[2]\shar_0$
$$\tilde\FF_\xi(\tilde v)\equiv
\FF_\xi(\tilde
v)\otimes_{\cO_{W[2]\shar\times_{\At}B}}\cO_{W[2]_0\shar}.
$$
First of all, since $B\to\At$ does not lie in any of the two
coordinate lines, to perform the elementary modification we do not
need to modify any of the nodes in $W[2]_0\shar$
and the sheaf $\tilde\FF_\xi(\tilde v)$ is locally free. On the
other hand, $\tilde\FF_\xi(\tilde v)$ is the cokernel of the
composite $\cF_\eta\equiv \cF_\eta\otimes\cI\lra \FF_\xi(\tilde
v)$ that is the unique lifting of
$$\cF_\eta\otimes_{\cO_{W[2]_0\shar}}\cI\lra \tilde\EE_\xi\otimes_{\cO_{W[2]_0\shar}}
\cI\lra \tilde\FF_\xi(\tilde v).
$$
Hence $\tilde\FF_\xi(\tilde v)$ fits into the exact sequence
$$0\lra \cL_\eta\lra \tilde\FF_\xi(\tilde v)\lra \cF_\eta\lra 0
$$
and the extension class of this exact sequence is a multiple of
$$v\in\PP\Ext^1_{W[2]\shar}(\cL_\eta,\cF_\eta)
$$
we started with. In particular the image $\Psi((\xi,v))$ depends
only on $v$. Now we pick an (analytic) open subset $U_\eta$ of
$\xi\in \PP\cW_\eta$ and $V_\eta\sub
\PP\Ext^1_{W[2]\shar}(\cL_\eta,\cF_\eta)$ so that
$$\tilde
\Psi_2|_{\tilde\Psi_2\upmo(U_\eta\times V_\eta)}\mh
\tilde\Psi_2\upmo(U_\eta\times V_\eta)\to U_\eta\times V_\eta
$$
is an isomorphism. Then the fact that $\Psi((\xi,v))$ depends only
on $v$ implies that the restriction of $\Psi$ to
$\tilde\Psi_2\upmo(U_\eta\times V_\eta)$ is the composite of the
second projection $U_\eta\times V_\eta\to V_\eta$ with a morphism
$V_\eta\to \bM^0$. Since $\Psi$ is morphism and
$$\tilde\Psi_2\upmo(\PP W_+\times_{\PFT}\PP W_-)\cong
\bE_b\times_{\PFT}\PP W_-,
$$
its restriction to $\tilde\Psi_2\upmo(\PP W_+\times_{\PFT}\PP
W_-)$ must be the composite of the second projection with a
morphism $\PP W_-\to\bM^0$. This proves the lemma.
\end{proof}



\def\CC{\mathbb{C} }
\def\PP{\mathbb{P} }
\def\RR{\mathbb{R} }
\def\Hom{\mathrm{Hom} }
\def\Eo{E|_{p_1} }
\def\Et{E|_{p_2} }

\section{The vanishing results}

The purpose of this section is to prove the main theorem of this
paper.
\begin{theo}
Let $M_{3,\chi}(Y)$ be the moduli space of stable vector bundles
over a smooth irreducible curve $Y$ of genus $g$ for $\chi\equiv
1, 2$ mod 3. Then we have $$c_i(M_{3,\chi}(Y))=0\ \ \ \text{for }
i>6g-5.$$
\end{theo}

Since $M_{3,1}(Y)\cong M_{3,2}(Y)$ by $E\to E^{\vee}\otimes L$ for
a fixed line bundle $L$ of degree $1$, we may assume $\chi\equiv
1$ mod 3, say $\chi=4$.

Our proof is induction on the genus $g$. When $g=1$,
$M_{3,4}(Y)\cong Y$ by Atiyah's theorem and hence we have the
vanishing result. We assume from now on that $g\ge 2$.

In the previous sections, we established the following diagram:
$$\xymatrix{
& &{\bM^0}\ar@{<..>}[r]_{flips}\ar[dl]_{normalization}&{\bM^{1/2}}
\ar@{<..>}[r]_{flips}&{\bM^1}\ar[dd]^{\overline{GL(3)}}\\
M_{3,4}(Y)\ar@{~>}[r]_{degeneration}&M_{3,4}(X_0)\\
& & & & {M_{3,7}(X)} }
$$

Suppose $c_i(M_{3,4}(X))=0$ for $i>6g-11$. We want to show that
$c_i(M_{3,4}(Y))=0$ for $i>6g-5$.

\subsection{Chern classes of $\bM^0$}

Let $S_0=M_{3,7}(X)$ and $E\to S_0\times X$ be a universal bundle.
Recall that a vector bundle on $X_n$ is $\alpha$-stable if and
only if its associated GPB $(V,V^0)$ is $\alpha$-stable. When
$\alpha=1^-$,  this is equivalent to $V$ being stable. Hence
$\bM^{1^-}$ is a fiber bundle over $S_0$ obtained by blowing up
$\bG^{1^-}=Gr(3,E|_{p_1+p_2})$, the Grassmannian bundle over
$S_0$. Let $\pi_1:S_1=\PP \Hom(E|_{p_1},E|_{p_2})\to S_0$ be the
projectivization of the bundle $\Hom(\Eo,\Et)$.

We blow up $S_1$ along the locus of rank 1 homomorphisms
$$B:=\PP E|_{p_1}^{\vee}\times_{S_0} \PP E|_{p_2}$$ and let
$\pi_2:S_2\to S_1$ be the blow-up map.

The exceptional divisor $\Delta_1:=\tilde B$ of $\pi_2:S_2\to S_1$
and the proper transform $\Delta_2$ of the locus of rank 2
homomorphisms in $S_1$ are normal crossing divisors. Let
$\Delta=\Delta_1+\Delta_2$.

\begin{lemm}
\label{vanlem} $c_i(\Omega_{S_2/S_0}(\log \Delta))=0$ for $i>6$.
Consequently, if $c_i(\Omega_{S_0})=0$ for $i>6g-11$, then
$c_i(\Omega_{S_2}(\log \Delta))=0$ for $i>6g-5$.
\end{lemm}

Since $\Omega_{S_2/S_0}(\log\Delta)$ is locally free of rank $8$,
it suffices to check that the 7th and 8th Chern classes vanish.
The proof is a lengthy computation. See the Appendix.

\vspace{.5cm}

Let $S_3$ be the result of two blow-ups of $\PP \big(\Hom
(\Eo,\Et)\oplus \cO_{S_0}\big)$ first along the section $\PP
\cO_{S_0}$ and then along the closure of the locus of rank 1
homomorphisms in $\Hom (\Eo,\Et)\subset \PP \big(\Hom
(\Eo,\Et)\oplus \cO_{S_0}\big)$. The obvious rational map $$\PP
\big(\Hom (\Eo,\Et)\oplus \cO_{S_0}\big)\dashrightarrow \PP \Hom
(\Eo,\Et)=S_1$$ becomes a $\PP^1$-bundle after the above first
blow-up and the preimage of $B$ is the center of the second
blow-up. Hence we get a $\PP^1$-bundle projection
$$\pi_3:S_3\to S_2$$
and $S_2$ naturally embeds into $S_3$.

Next, we blow up $S_3$ along $\Delta_1=\tilde B\subset S_2$ which
lies in $S_3$ as a codimension 2 subvariety. Let $\pi_4:S_4\to
S_3$ be the blow-up. Then by local computation, $S_4$ is the same
as the result of the blow-ups of $\PP \big(\Hom (\Eo,\Et)\oplus
\cO_{S_0}\big)$, first along $\PP \cO_{S_0}$, second along $B$
which lies in $\PP \Hom (\Eo,\Et)$ and finally along the proper
transform of the closure of the locus rank 1 homomorphisms in
$\Hom (\Eo,\Et)\subset \PP \big(\Hom (\Eo,\Et)\oplus
\cO_{S_0}\big)$. So if we finally blow up $S_4$ along the proper
transform of $\Delta_2\subset S_2$ which lies in $S_4$ as a
codimension 2 subvariety, then we obtain the moduli space
$\bM^{1^-}$ of $1^-$-stable bundles which we also denote by $S_5$
and the last blow-up is denoted by $\pi_5:S_5\to S_4.$ Recall that
$\bM^{1^-}$ has six divisors
$\tilde\bY_0,\tilde\bY_1,\tilde\bY_2,\tilde\bZ_0,\tilde\bZ_1,\tilde\bZ_2$.
Let $D$ be the sum of these.

\begin{lemm}\label{vanlem2}
$c_i(\Omega_{\bM^1}(\log D))=0$ for $i>6(g-1)$ if and only if
$c_i(\Omega_{S_2}(\log \Delta))=0$ for $i>6(g-1)$.
\end{lemm}
\begin{proof}
The proof is due to Gieseker \cite{Gie}. Notice that we have four
divisor
$\widetilde{\Delta}_0,\widetilde{\Delta}_1,\widetilde{\Delta_2},\widetilde{\Delta}_3$
in $S_3$ which are the images of
$\tilde\bZ_0,\tilde\bY_1,\tilde\bY_2,\tilde\bY_0$ respectively.
Let $\widetilde{\Delta}$ be the sum of $\widetilde{\Delta}_i$'s.
Notice that $\widetilde{\Delta}_1=\pi_3^{-1}(\Delta_1)$ and
$\widetilde{\Delta}_2=\pi_3^{-1}(\Delta_2)$. Hence we have an
exact sequence
$$0\to \pi_3^*\Omega_{S_2}(\log \Delta)\to \Omega_{S_3}(\log
\widetilde{\Delta})\to \Omega_{S_3/S_2}(\log
(\widetilde{\Delta}_0+\widetilde{\Delta}_3))\to 0.$$

But the line bundle $\Omega_{S_3/S_2}(\log
(\widetilde{\Delta}_0+\widetilde{\Delta}_3))$ is trivial since we
can find a nowhere vanishing section as follows: as $\pi_3$ is a
$\PP^1$-bundle there is an open covering $\{U_i\}$ of $S_2$ and
rational functions $z_i$ on $\pi_3^{-1}(U_i)$ with a simple pole
at $\widetilde{\Delta}_0$ and a simple zero at
$\widetilde{\Delta}_3$. Because $z_i=z_jf_{ij}$ for a nowhere
vanishing function $f_{ij}$ on $U_i\cap U_j$, $dz_i/z_i$ gives a
well-defined section of
$\Omega_{S_3/S_2}(\log(\widetilde{\Delta}_0+\widetilde{\Delta}_3))$.

Next, $S_4$ was obtained by blowing up along the intersection of
two divisors $\widetilde{\Delta}_0$ and $\widetilde{\Delta}_1$ in
$S_3$. Let $Z_1'$ be the exceptional divisor of $\pi_4$. By a
local computation, we get $$\pi_4^*\Omega_{S_3}(\log
\widetilde{\Delta})\cong \Omega_{S_4}(\log
(\widetilde{\Delta}+Z_1')).$$ By the same argument, we see that
$$\pi_5^*\Omega_{S_4}(\log
(\widetilde{\Delta}+Z_1'))\cong \Omega_{\bM^1}(\log D).$$ The
lemma now follows immediately.
\end{proof}

By Lemmas \ref{vanlem} and \ref{vanlem2}, we deduce the vanishing
of Chern classes for $\bM^{1^-}$.
\begin{coro} $c_i(\Omega_{\bM^{1^-}}(\log D))=0$ for $i>6g-5$
\end{coro}

\subsection{From ${\bf M}^{1^-}$ to ${\bf M}^0$}
The goal of this subsection is to show the following.

\begin{prop}
$c_i(\Omega_{\bM^{1^-}}(\log D))=0$ for $i>6g-5$ iff
$c_i(\Omega_{\bM^0}(\log D))=0$ for $i>6g-5$.
\end{prop}

Recall that $\bM^0$ is obtained from $\bM^{1^-}$ by a sequence of
flips along subvarieties each of which lies in the intersection of
two of the six divisors and the blow-up center is not contained in
any other divisor.

We use a lemma from \cite{Gie}. Let $J$ be the base of a flip. In
other words, there are two vector bundles $E$ and $F$ over $J$ and
a variety $S$ into which $\PP E$ is embedded. And the normal
bundle to $Z=\PP E$ is the pull-back of $F$ tensored with
$\cO_{\PP E}(-1)$. Let $\tilde S$ be the blow-up of $S$ along $Z$
and $S'$ be the blow-down of $\tilde S$ along the $\PP
E$-direction. Then $\tilde S$ is the blow-up of $S'$ along $Z'=\PP
F$. Suppose that there are normally crossing smooth divisors $D_i$
in $S$ such that $Z$ is contained in $D_1\cap D_2$ as a smooth
subvariety but no other divisor contains $Z$. Let $D=\sum D_i$ and
$D'$ be its proper transform in $S'$. The following lemma is from
\cite[\S12]{Gie}.\footnote{Gieseker assumed that $Z\subset D_1\cap
D_2$ and $Z\cap D_k=\emptyset$ for $k\ne 1,2$. But the same proof
works as long as $Z\cap D$ is a smooth divisor in $Z$.}

\begin{lemm}
Suppose the top $k$ Chern classes of $J$ vanish. Then
$c_i(\Omega_S(\log D))=0$ for $i>\dim S-k-1$ iff
$c_i(\Omega_{S'}(\log D'))=0$ for $i>\dim S-k-1$.
\end{lemm}

The flip bases for $\alpha=2/3$ are as follows: The moduli spaces
$\bP^{2/3}_{2,5,p_i}$ of stable parabolic bundles with parabolic
weight 1/3 and quasi-parabolic structure at $p_i$ for $i=1,2$ lie
in the moduli of GPBs $\bG^{2/3}_{2,5,1}$. Let
$\widetilde{\bG}^{2/3}_{2,5,1}$ be the blow-up of
$\bG^{2/3}_{2,5,1}$ along $\bP^{2/3}_{2,5,p_1}\cup
\bP^{2/3}_{2,5,p_2}$. Then the flip bases are
\begin{itemize} \item $\widetilde{\bG}^{2/3}_{2,5,1}\times Jac(X)$
 \item a Jacobian of $X$ times a
$\PP^1$ bundle over $\bP^{2/3}_{2,5,p_1}$ \item a Jacobian of $X$
times a $\PP^1$ bundle over $\bP^{2/3}_{2,5,p_2}$.
\end{itemize}

Because the underlying vector bundle of a parabolic bundle or a
GPB above is stable, all these three moduli spaces are fiber
bundles over $M_{2,5}(X)$. By Gieseker's theorem \cite{Gie}, we
know that the top $2g-3$ Chern classes of $M_{2,5}(X)$ vanish.
Hence the top $3g-4$ Chern classes of the flip bases for
$\alpha=2/3$ vanish. From the above lemma, we have
$$c_i(\Omega_{\bM^{1/2}}(\log D))=0\ \ \ \text{ for }i>6g-5.$$

We can similarly deal with the flip bases for $\alpha=1/3$: the
moduli spaces $\bP^{1/3}_{2,4,p_i}$ of stable parabolic bundles
with parabolic weight 2/3 and quasi-parabolic structure at $p_i$
for $i=1,2$ lie in the moduli of GPBs\footnote{The choice of 1
dimensional subspace $V_1$ of $E|_{p_1}$ gives rise to the 3
dimensional subspace $V=V_1+E|_{p_2}$. This is a GPB in
$\bG^{1/3}_{2,4,3}$.} $\bG^{1/3}_{2,4,3}$. Let
$\widetilde{\bG}^{1/3}_{2,4,3}$ be the blow-up of
$\bG^{1/3}_{2,4,3}$ along $\bP^{1/3}_{2,4,p_1}\cup
\bP^{1/3}_{2,4,p_2}$. Then the flip bases are
\begin{itemize}
\item $\widetilde{\bG}^{1/3}_{2,4,3}\times Jac(X)$
 \item a Jacobian of $X$ times a
$\PP^1$ bundle over $\bP^{1/3}_{2,4,p_1}$
 \item a Jacobian of $X$ times a
$\PP^1$ bundle over $\bP^{1/3}_{2,4,p_2}$.
\end{itemize}

Because the underlying vector bundle of a parabolic bundle above
is stable, the moduli spaces $\bP^{1/3}_{2,4,p_1}$ and
$\bP^{1/3}_{2,4,p_2}$ are fiber bundles over $M_{2,3}(X)$. By
Gieseker's theorem \cite{Gie}, we know that the top $2g-3$ Chern
classes of $M_{2,3}(X)$ vanish.

The moduli space $\widetilde{\bG}^{1/3}_{2,4,3}$ is not a fiber
bundle over $M_{2,3}(X)$ but this is isomorphic to a divisor in
Gieseker's moduli space: consider the universal family $\cF$ over
$\widetilde{\bG}^{1/3}_{2,4,3}\times X$. Blow up this space along
$\bP^{1/3}_{2,4,p_1}\times p_1$ and $\bP^{1/3}_{2,4,p_2}\times
p_2$. Perform elementary modifications as in \S2.4 so that we get
a family of curves over $\widetilde{\bG}^{1/3}_{2,4,3}$ and a
vector bundle on the family of curves. The restriction of this
vector bundle to the proper transforms of
$\widetilde{\bG}^{1/3}_{2,4,3}\times p_1$ and
$\widetilde{\bG}^{1/3}_{2,4,3}\times p_2$ is equipped with a
choice of basis and we can glue the rank 2 bundle $\cO\oplus
\cO(1)$ over a rational curve $\PP^1$ to get a vector bundle over
the family of nodal genus $g$ curves. It is elementary to check
that this is a family of bundles in the Gieseker's moduli space
$\bM^{1/3}_{2,3}$ for the rank 2 case and so we get a morphism
$$\widetilde{\bG}^{1/3}_{2,4,3}\to \bM^{1/3}_{2,3}.$$ It is now an
easy matter to check that this morphism is bijective onto a
divisor of rank 1 locus in the Gieseker's moduli space. Hence,
$\widetilde{\bG}^{1/3}_{2,4,3}$ becomes a fiber bundle over
$M_{2,1}(X)$ after a flip whose base is the product of two
Jacobians over $X$. By Gieseker's lemma again, we deduce that the
top $2g-3$ Chern classes of $\widetilde{\bG}^{1/3}_{2,4,3}$ vanish
and hence the top $3g-4$ Chern classes of all the flip bases for
$\alpha=1/3$ vanish. From Gieseker's lemma, we have
$$c_i(\Omega_{\bM^{0}}(\log D))=0\ \ \ \text{ for }i>6g-5.$$

The argument at the end of \S13 in \cite{Gie} enables us to deduce
the vanishing Chern classes of the general member of the family
$\bM_{3,4}(\mathfrak W)$ from the vanishing of the Chern classes
of $\Omega_{\bM^{0}}(\log D)$. So we conclude that
$$c_i(M_{3,4}(Y))=0\ \ \ \text{ for } i>6g-5.$$


\section{Appendix}

{\small\small

The purpose of this appendix is to prove Lemma \ref{vanlem}.

Notice that $E|_{p_1}$ is isomorphic to $E|_{p_2}$. Let
$a_1,a_2,a_3$ be the Chern roots of $E|_{p_1}\cong E|_{p_2}$ and
let $\xi=c_1(\cO_{S_1}(-1))$. Then the Chern roots of
$\Hom(\Eo,\Et)$ are $0,0,0,\pm(a_1-a_2),\pm(a_2-a_3),\pm(a_3-a_1)$
and thus the cohomology ring $H^*(S_1)$ is  the polynomial algebra
$H^*(S_0)[\xi]$ over $H^*(S_0)$ modulo the relation
\begin{equation}\label{r1}\xi^3(\xi^2-(a_1-a_2)^2)
(\xi^2-(a_2-a_3)^2)(\xi^2-(a_3-a_1)^2).\end{equation} From the
exact sequence
$$0\to \Omega_{S_1/S_0}\to \pi_1^*\Hom(\Eo,\Et)\otimes
\cO_{S_1}(-1)\to \cO\to 0$$ we deduce that the total Chern class
of the relative cotangent bundle of $S_1$ over $S_0$ is
$$
c(\Omega_{S_1/S_0})=(1+\xi)^3((1+\xi)^2-(a_1-a_2)^2)
((1+\xi)^2-(a_2-a_3)^2)((1+\xi)^2-(a_3-a_1)^2)).$$

Similarly, we can describe the cohomology rings and the total
Chern classes of the relative tangent bundles of $\PP
E|_{p_1}^{\vee}$ and $\PP E|_{p_2}$ over $S_0$. Let
$$u=c_1(\cO_{\PP E|_{p_1}^{\vee}}(-1))+\frac13 c_1(E|_{p_1})$$
$$v=c_1(\cO_{\PP E|_{p_2}}(-1))-\frac13 c_1(E|_{p_1}).$$ We
intentionally shifted the generators to make our computation
simpler. Then, we have
$$H^*(\PP E|_{p_1}^{\vee})=H^*(S_0)[u]/\langle u^3+\alpha u+\beta \rangle$$
$$H^*(\PP E|_{p_2})=H^*(S_0)[v]/\langle v^3+\alpha v-\beta \rangle$$
where
$$\alpha=c_2(E|_{p_1})-\frac13 c_1^2(E|_{p_1})$$
$$\beta=c_3(E|_{p_1})-\frac13
c_1(E|_{p_1})c_2(E|_{p_1})+\frac2{27}c_1^3(E|_{p_1}).$$ Also we
have
$$c(T_{\PP E|_{p_1}^{\vee}/S_0})=(1-u)^3+\alpha
(1-u)-\beta=1-3u+3u^2+\alpha$$
$$c(T_{\PP E|_{p_2}/S_0})=(1-v)^3+\alpha
(1-v)+\beta=1-3v+3v^2+\alpha.$$

Using $\alpha$ and $\beta$, we can rewrite
$$c(\Omega_{S_1/S_0})=(1+\xi)^3((1+\xi)^6+6\alpha
(1+\xi)^4+9\alpha^2(1+\xi)^2+4\alpha^3+27\beta^2)$$ as one can
check by direct computation.

From \cite{Ful}, we get the exact sequences
$$0\to \cO_{\tilde B}(-1)\to g^*N_{B/S_1}\to F\to 0$$
$$0\to T_{S_2}\to \pi_2^*T_{S_1}\to \jmath_*F\to 0$$
where $g:\tilde B\to B$ is the restriction of $\pi_2$ to the
exceptional divisor $\tilde B$ and $\jmath$ is the inclusion of
$\tilde B$. Therefore, we have
\begin{equation}
c(T_{S_2})=\pi_2^*c(T_{S_1})/c(\jmath_*F)\label{chexseq2}\end{equation}
\begin{equation}\label{chexseq}
c(F)=g^*c(N_{B/S_1})/c(\cO_{\tilde B}(-1)).\end{equation}

Since $\cO_{S_1}(-1)$ restricts to $\cO_{\PP
E|_{p_1}^*}(-1)\boxtimes \cO_{\PP E|_{p_2}}(-1)$, $\xi$ restricts
to $$c_1(\cO_{\PP E|_{p_1}^*}(-1))+c_1(\cO_{\PP
E|_{p_2}}(-1))=u+v.$$ The restriction of the relative tangent
bundle $T_{S_1/S_0}$ to $B$ has total Chern class
$$(1-u-v)^3((1-u-v)^6+6\alpha
(1-u-v)^4+9\alpha^2(1-u-v)^2+4\alpha^3+27\beta^2)$$ while the
total Chern class of the relative tangent bundle $T_{B/S_0}$ is
$$(1-3u+3u^2+\alpha)(1-3v+3v^2+\alpha).$$
Hence the total Chern class of the normal bundle $N_{B/S_1}$ is
$$\frac{(1-u-v)^3((1-u-v)^6+6\alpha
(1-u-v)^4+9\alpha^2(1-u-v)^2+4\alpha^3+
27\beta^2)}{(1-3u+3u^2+\alpha)(1-3v+3v^2+\alpha)}
$$
which is by direct computation equal to
$$1-6\xi+(15\xi^2+4\alpha-3uv)-(15\xi^3+9\alpha
\xi)+(6\xi^4+6\alpha \xi^2)$$ with $\xi|_{B}=u+v$ understood by
abuse of notations.

Let $\eta=c_1(\cO_{S_2}(\tilde B))\in H^2(S_2)$. Then
$\eta|_{\tilde B}=c_1(\cO_{\tilde B}(-1))$ since $\cO_{S_2}(\tilde
B)|_{\tilde B}\cong \cO_{\tilde B}(-1)$. So we have
$$H^*(\tilde B)\cong H^*(B)[\eta]/
\langle
\eta^4+6\xi\eta^3+(15\xi^2+4\alpha-3uv)\eta^2+(15\xi^3+9\alpha
\xi)\eta+(6\xi^4+6\alpha \xi^2)\rangle.$$ In fact, it is easy to
see that the above relation lifts to a relation
\begin{equation}\label{r2}
\eta^4+6\xi\eta^3+(15\xi^2+4\alpha)\eta^2-3\jmath_*(uv)\eta+(15\xi^3+9\alpha
\xi)\eta+(6\xi^4+6\alpha \xi^2)=0
\end{equation}
in $H^*(S_2)$.

From (\ref{chexseq}), we have
$$\begin{array}{ll}
&c(F)=c(N_{B/S_1})/(1+\eta)=\\
&1-(6\xi+\eta)+(15\xi^2+4\alpha-3uv +6\xi\eta+\eta^2)\\&
-(15\xi^3+9\alpha\xi+15\xi^2\eta+4\alpha\eta+6\xi\eta^2+\eta^3-3uv\eta).
\end{array}$$

By local computation we have
$$c(\cO_{S_2}(\Delta_1))=1+\eta,\ \ \
c(\cO_{S_2}(\Delta_2))=1-3\xi-2\eta.$$

Hence we have\small
\begin{equation}\label{relchcl1}
c(\Omega_{S_2/S_0}(\log \Delta))=\frac{(1+\xi)^3((1+\xi)^6+6\alpha
(1+\xi)^4+9\alpha^2(1+\xi)^2+4\alpha^3+27\beta^2)}{(1-\eta)
(1+3\xi+2\eta)}c(\Omega_{S_2/S_1}).\end{equation}
 For $c(\Omega_{S_2/S_1})$,
we compute
\begin{equation}\label{relchcl2}
c(T_{S_2/S_1})=c(T_{S_2})/\pi_2^*c(T_{S_1})=\frac1{c(\jmath_*F)}\end{equation}
and change the signs of the terms of degree $\equiv 2$ (mod $4$).

It is a consequence of the Grothendieck-Riemann-Roch theorem that
\begin{equation}c(\jmath_*F)=1-\jmath_*\left(\frac1{\eta}\big(1-\prod
 \frac{1+b_i}{1+b_i-\eta}\big)\right)\label{relchcl3}\end{equation}
where $b_i$ are the Chern roots of $F$, i.e.
$$\begin{array}{ll}&\prod (1+b_i) = 1-(6\xi+\eta)+(15\xi^2+4\alpha-3uv
+6\xi\eta+\eta^2)\\&
-(15\xi^3+9\alpha\xi+15\xi^2\eta+4\alpha\eta+6\xi\eta^2+\eta^3-3uv\eta).
\end{array}$$
By expanding, we see that
$$\begin{array}{ll}
&\prod (1+b_i-\eta) =1-(6\xi+4\eta)+(15\xi^2+4\alpha-3uv
+18\xi\eta+6\eta^2)\\&
-(15\xi^3+9\alpha\xi+30\xi^2\eta+8\alpha\eta+18\xi\eta^2+4\eta^3-6uv\eta).
\end{array}
$$
Hence, we have 
\begin{equation}\lab{50.1}
\frac1{\eta}\big(1-\prod
\frac{1+b_i}{1+b_i-\eta}\big)=\frac{F}{1-D}.
\end{equation}
Here \small
$$A=(6\xi+4\eta)-(15\xi^2+4\alpha+18\xi\eta+6\eta^2)
+(15\xi^3+9\alpha\xi+30\xi^2\eta+8\alpha\eta+18\xi\eta^2+4\eta^3-6uv\eta),
$$
$$B=-3+(12\xi+5\eta)
-(15\xi^2+12\xi\eta+3\eta^2+4\alpha),
$$
$D=A+3uv$ and $F=B+3uv$. Then by expanding\footnote{We used Maple
7 for the computations in this section.} (\ref{50.1}) and
collecting all terms of degrees up to 14, we obtain \small
$$
\begin{array}{ll} &27u^{3}v^{3} + 9u^{2}v^{2}+ 3 uv + 18BA^{5}uv +
12BA^{3}uv+ 3Auv\\&+ 3A^{2}uv+18Au^{2}v^{2}+
3A^{3}uv+27A^{2}u^{2}v^{2}+ 81Au^{3}v^{3}+ 3A^{4}uv+ 36A^{3} u^{2}
v^{2}\\&
 + 162A^{2} u^{3}  v^{3}
 + BA^{3} + BA^{4} + B + 3A^{5}uv + 45A^{4}u^{2}v^{2}\\&
 + 270A^{3}u^{3}v^{3}
 + 54A^{5} u^{2} v^{2} + 405A^{4} u^{3} v^{3} + 3uvA^{7}
 + 567u^{3}v^{3}A^{5} \\&+ 63u^{2}v^{2}A^{6}+ 27Bu^{3}v^{3}
 + 9Bu^{2}v^{2}+ 3Buv+ 3A^{6}uv+ BA^{2} \\&+ BA
 + BA^{5}+ BA^{6} + BA^{7} + 6BAuv+ 9BA^{2}uv\\&+ 27BAu^{2}v^{2}
 + 54BA^{2}u^{2}v^{2}+ 108BAu^{3}v^{3}+ 15BA^{4}uv +
90BA^{3}u^{2}v^{2}\\&+ 270BA^{2}u^{3}v^{3} + 135BA^{4}u^{2} v^{2}+
540BA^{3}u^{3}v^{3} + 189BA^{5}u^{2}v^{2} \\&+ 945BA^{4}u^{3}v^{3}
+ 21BA^{6}uv.\end{array}$$

By the projection formula and $\jmath_*1=\eta$,
$$\jmath_*\left(\frac1{\eta}\big(1-\prod
 \frac{1+b_i}{1+b_i-\eta}\big)\right)=\jmath_* \left(F/(1-D)\right)$$
is, up to degree 16, equal to\small
$$\begin{array}{ll} &27\jmath_*(u^{3}v^{3}) + 9\jmath_*(u^{2}v^{2})
+ 3 \jmath_*(uv) + 18BA^{5}\jmath_*(uv) + 12BA^{3}\jmath_*(uv)+
3A\jmath_*(uv)\\&+ 3A^{2}\jmath_*(uv)+18A\jmath_*(u^{2}v^{2})+
3A^{3}\jmath_*(uv)+27A^{2}\jmath_*(u^{2}v^{2})\\ &
+81A\jmath_*(u^{3}v^{3})+ 3A^{4}\jmath_*(uv)+ 36A^{3}
\jmath_*(u^{2} v^{2})\\&
 + 162A^{2} \jmath_*(u^{3}  v^{3})
 + BA^{3}\eta + BA^{4}\eta + B\eta + 3A^{5}\jmath_*(uv) + 45A^{4}\jmath_*(u^{2}v^{2})\\&
 + 270A^{3}\jmath_*(u^{3}v^{3})
 + 54A^{5} \jmath_*(u^{2} v^{2}) + 405A^{4} \jmath_*(u^{3} v^{3}) + 3\jmath_*(uv)A^{7}
 + 567\jmath_*(u^{3}v^{3})A^{5} \\&+ 63\jmath_*(u^{2}v^{2})A^{6}+
 27B\jmath_*(u^{3}v^{3})
 + 9B\jmath_*(u^{2}v^{2})+ 3B\jmath_*(uv)+ 3A^{6}\jmath_*(uv)+ BA^{2}\eta \\&+
 BA\eta + BA^{5}\eta+ BA^{6}\eta + BA^{7}\eta + 6BA\jmath_*(uv)+ 9BA^{2}\jmath_*(uv)\\
 &+ 27BA\jmath_*(u^{2}v^{2})
 + 54BA^{2}\jmath_*(u^{2}v^{2})+ 108BA\jmath_*(u^{3}v^{3})+ 15BA^{4}\jmath_*(uv) +
90BA^{3}\jmath_*(u^{2}v^{2})\\&+ 270BA^{2}\jmath_*(u^{3}v^{3}) +
135BA^{4}\jmath_*(u^{2} v^{2})+ 540BA^{3}\jmath_*(u^{3}v^{3}) +
189BA^{5}\jmath_*(u^{2}v^{2})
\\&+ 945BA^{4}\jmath_*(u^{3}v^{3}) + 21BA^{6}\jmath_*(uv).\end{array}$$

Substitute the above expression into (\ref{relchcl3}) and expand
(\ref{relchcl2}) up to degree 16. Change the signs of the terms of
degree $\equiv 2$ mod 4 and plug it into (\ref{relchcl1}).

Now we can compute the Chern classes by direct computation from
(\ref{relchcl1}). The 7th Chern class is, up to sign, equal
to\small
$$\begin{array}{ll}
&- 378\xi^{5}\eta^{2} + 36\alpha\xi\jmath_*(uv)\eta -
72\alpha\xi^{2}\eta^{3} - 138\alpha^{2}\xi\eta^{2}-
492\alpha^{2}\xi^{2}\eta - 516\alpha\xi^{3}\eta^{2}\\
&- 1056\alpha\xi^{4}\eta - 540\xi^{2}\jmath_*(uv)\eta^{2}
 + 558\alpha\xi^{2}\jmath_*(uv) + 108\alpha\jmath_*(u^2v^2) -
810\xi^{2}\jmath_*(u^2v^2) \\
&- 630\xi^{3}\jmath(uv)\eta - 564\xi^{6}\eta - 252\xi^{7 } -
504\alpha\xi^{5}  - 252\alpha^2\xi^{3}  - 72\xi^{4}\eta^{3}\\& -
54\alpha^{2}\jmath_*(uv)
 + 72\xi^{4}\jmath_*(uv)
- 126\xi\eta^{3}\jmath_*(uv)- 54\eta\beta^2 - 54\jmath_*(u^3v^3)
\end{array}$$
and the 8th Chern class is\small
$$\begin{array}{ll}
&1332\jmath_*(uv)\eta^2\alpha\xi+2178\jmath_*(uv)\xi^2\alpha\eta-
2835\jmath_*(u^2v^2)\xi^3
-1143\jmath_*(u^2v^2)\eta^3+2214\jmath_*(uv)\xi^5\\
&-6939\jmath_*(u^2v^2)\xi^2\eta
-4968*\jmath_*(u^2v^2)\xi\eta^2+1485\jmath_*(u^2v^2)\alpha\xi
+774\jmath_*(u^2v^2)\alpha\eta\\
&-72\xi^5\eta^3+132\xi^6\eta^2
-336\alpha^2\xi\eta^3-588\alpha^2\xi^2\eta^2-408\alpha\xi^3\eta^3
-456\alpha\xi^4\eta^2+1152\jmath_*(u^2v^2)\eta^3\\
&+6372\xi^2\jmath_*(u^2v^2)\eta
-828\alpha\jmath_*(u^2v^2)\eta+3942\eta^2\jmath_*(u^2v^2)\xi
+846\jmath_*(uv)\xi^4\eta+30\jmath_*(uv)\alpha^2\eta\\
&+1008\jmath_*(uv)\eta^2\xi^3+72\jmath_*(uv)\xi^2\eta^3
-6\jmath_*(uv)\alpha\eta^3+3150\jmath_*(uv)\alpha\xi^3\\
&-828\jmath_*(uv)\alpha^2\xi-567\eta\jmath_*(u^3v^3)-648\alpha\xi^5\eta
-498\alpha^2\xi^3\eta+132\alpha^3\xi\eta-18\xi^7\eta+153\xi^8\\
&-810\xi\jmath_*(u^3v^3)+54\beta^2\eta^2+165\alpha^2\xi^4
+6\xi^2\alpha^3+81\xi^2\beta^2+312\alpha\xi^6+594\jmath_*(u^2v^2)\jmath_*(uv).
\end{array}$$

Notice that the 7th Chern class is the image by $\jmath_*$
of\small \begin{equation} \begin{array}{ll} & - 378\xi^{5}\eta+
36\alpha uv\xi\eta- 72\alpha\xi^{2}\eta^{2}- 138\alpha^2\xi\eta-
492\alpha^2\xi^{2}- 516\alpha\xi^{3}\eta - 1056\alpha\xi^{4}\\& -
540\xi^{2}\eta^{2}uv+ 558\alpha uv\xi^{2}+ 108\alpha u^{2}v^{2}-
810\xi^{2}u^{2}v^{2}  - 630\xi^{3} uv\eta- 564\xi^{6}\\&-
72\xi^{4}\eta^{2}- 54\alpha^2 u v+ 72\xi^{4} uv- 126\xi\eta^{3}uv
 - 54\beta^2 - 54u^{3}v^{3} \\&+ 42\alpha\xi
(\eta^{3} + 6\xi\eta^{2} + (15\xi^{2} + 4\alpha )\eta - 3uv\eta +
15\xi^{3} +
9\alpha\xi) \\
&+ 42\xi^{3}(\eta^{3} + 6\xi\eta^{2} + (15\xi^{2} + 4\alpha )\eta
- 3uv\eta + 15\xi^{3} + 9\alpha\xi)
\end{array}\label{chk7}\end{equation}
with $\xi^4+\alpha\xi^2$ replaced by
$$\jmath_*[-\left(\eta^3+6\xi\eta^2+(15\xi^2+4\alpha)\eta
-3uv\eta+15\xi^3+9\alpha\xi\right)/6]$$ from the relation
(\ref{r2}).

 This is a class in $H^*(\widetilde B)$ which is a polynomial algebra
 over $H^*(S_0)$ generated by $u$, $v$, and $\eta$ with the relations
$\xi=u+v$, $$u^3+\alpha u+\beta=0$$ $$v^3+\alpha v-\beta=0$$
$$\eta^4+6\xi\eta^3+(15\xi^2+4\alpha-3uv)\eta^2+(15\xi^3+9\alpha
\xi)\eta+(6\xi^4+6\alpha \xi^2)=0.$$ Using Gr\"obner package, one
can check that the class (\ref{chk7}) is zero. Therefore we proved
that the 7th Chern class vanishes.

We apply the same strategy for the 8th Chern class. The only term
we cannot express as the image of $\jmath_*$ in the above fashion
using (\ref{r2}) is the term $81\xi^2\beta^2$. Let
$$\mu=(\xi^4+\alpha\xi^2)(\xi^2+\alpha)(\xi^2+4\alpha)+27\xi^2\beta^2.$$
This is exactly the relation (\ref{r1}) divided by $\xi$ and thus
we have $$\xi\mu=0.$$ Then by the relation (\ref{r2}) as above
$c_8(\Omega_{S_2/S_0}(\log \Delta))-3\mu$ is the image by
$\jmath_*$ of\small
$$\begin{array}{ll}
&- 72\xi^{5}\eta^{2} - 336\alpha^2\xi\eta^{2} -
588\alpha^2\xi^{2}\eta
 - 408\alpha\xi^{3}\eta^{2} - 456\alpha\xi^{4}\eta + 132\xi^{6}\eta - 18\xi^{7} \\
& - 648\alpha\xi^{5} + 1008\xi^{3}uv\eta^{2} - 498\alpha^2\xi^{3}
+ 54\eta\beta^2 + 132\xi\alpha^3\\
& -  \frac{ 51}{2} \xi^{4}(y^{3} + 6\xi\eta^{2} + (15\xi^{2} +
4\alpha - 3uv)\eta + 15\xi^{3} + 9\alpha\xi) \\
&- \alpha^2 (y^{3} + 6\xi\eta^{2} + (15\xi^{2} + 4\alpha -
3uv)\eta + 15\xi^{3} + 9\alpha\xi) + 2214\xi^{5}uv -
2835\xi^{3}u^{2}v^{2}\\
& - \frac {53}{2} \alpha\xi^{2}(y^{3} + 6\xi\eta^{2} + (15\xi^{2}
+ 4\alpha - 3uv)\eta + 15\xi^{3} + 9\alpha\xi) - 810 u^{3}v^{3}\xi
+ 27u^{3}v^{3}\eta\\
& + 3150\alpha uv\xi^{3} + 1332\alpha uv\xi\eta^{2}
 + 2178\alpha u v\xi^{2}\eta - 54\alpha u^{2} v^{2}\eta \\&+
9\eta^{3} u^{2} v^{2} + 30\alpha^2 u v\eta + 1485\alpha u^{2}
v^{2}\xi - 828 \alpha^2 u v\xi + 846\xi^{4} u v\eta \\ &-
1026\eta^{2} u^{2} v^{2}\xi  - 567\xi^{2} u^{2} v^{2}\eta +
72\xi^{2}\eta^{3} u v - 6\eta^{3}\alpha uv \\&  +\frac {1}{2}
(\xi^{2} + \alpha)(\xi^{2}
 + 4\alpha)(\eta^{3} + 6\xi\eta^{2} + (15\xi^{2} + 4\alpha - 3u
v)\eta + 15\xi^{3} + 9\alpha\xi)
\end{array}
$$
If we simplify this expression using the Gr\"obner package for the
ring $H^*(\widetilde B)$, we get
\begin{equation}\label{ch8van}\begin{array}{ll} & 3\eta^{3} u^{2} v^{2} +
\alpha u^{2}\eta^{3} + 2 \eta^{3}\alpha u v - 3\eta^{3} u\beta +
\alpha v^{2}\eta^{3}
 + 3\eta^{3} v\beta + a^{4}\eta^{3} \\
& + 6\alpha u^{2} v^{2}\eta + 9\eta u^{2} v\beta + 4\alpha^2\eta
u^{2}
- 9\eta u v^{2}\beta + 2\alpha^2 u v\eta - 6\eta\alpha u\beta \\
& + 4\alpha^2\eta v^{2} + 6\eta\alpha v\beta + 4\alpha^3\eta +
9\beta^2\eta. \end{array}\end{equation}

If we multiply $\eta$ to this expression (\ref{ch8van}), we get
zero! Hence $c_8(\Omega_{S_2/S_0}(\log \Delta))-3\mu$ lies in
$H^*(S_1)$ because its restriction to $\widetilde B$ is exactly
the above expression multiplied by $\eta$ and is equal to zero.

If we multiply $\xi=u+v$ to (\ref{ch8van}) and simplify using the
Gr\"obner package, we get zero! By the projection formula, this
implies that $c_8(\Omega_{S_2/S_0}(\log \Delta))-3\mu$ lies in the
kernel of multiplication by $\xi$
$$\xi:H^*(S_1)\to H^{*+2}(S_1)$$
 which is exactly $H^*(S_0)\mu$. Therefore, we deduce that
$$c_8(\Omega_{S_2/S_0}(\log \Delta))-3\mu= c\mu$$
for some rational number $c$. To compute $c$, we restrict the
image by $\jmath_*$ of (\ref{ch8van}) to a fiber of
$\pi_1\circ\pi_2:S_2\to S_0$ so that $\alpha=0$ and $\beta=0$.
Using the explicit relations it is now an elementary exercise to
check that $c=-3$. Hence, we conclude that
$c_8(\Omega_{S_2/S_0}(\log \Delta))=0$. }


\end{document}